\documentclass[letter]{amsart}

\usepackage{amssymb,amsmath}
\usepackage{stmaryrd,paralist}
\usepackage{bbm}
\usepackage{graphics,graphicx}

\usepackage[latin1]{inputenc}
\usepackage{subfigure}
\usepackage[margin=1.4 in]{geometry}
\usepackage{hyperref,color}
\usepackage[dvipsnames]{xcolor}
\graphicspath{{Figures/}} 
\graphicspath{{./Figures/}}

\usepackage{url}

\newcommand{\addllncs}[1]{}

\newtheorem{theo}{Theorem}
\numberwithin{theo}{section}
\newtheorem{thm}[theo]{Theorem}
\newtheorem{lem}[theo]{Lemma}
\newtheorem{lemma}[theo]{Lemma}
\newtheorem{prop}[theo]{Proposition}
\newtheorem{proposition}[theo]{Proposition}

\newtheorem{definition}[theo]{Definition}

\theoremstyle{remark}
\newtheorem{remark}[theo]{Remark}
\newtheorem{example}[theo]{Example}

{\medbreak}

\def\ni{\noindent}

\newcommand{\cacher}[1]{}

\def\NN{\mathbb{N}}
\def\RR{\mathbb{R}}

\def\br{\mathbb{R}}

\def\cA{\mathcal{A}}

\def\cO{\mathcal{O}}

\def\cT{\mathcal{T}}
\def\cS{\mathcal{S}}

\def\cW{\mathcal{W}}
\def\cL{\mathcal{L}}

\newcommand{\al}{\alpha}

\newcommand{\eps}{\epsilon}

\newcommand{\ov}[1]{\overline{#1}}
\newcommand{\wt}[1]{\widetilde{#1}}

\newcommand{\fig}[3]{\begin{figure}[h!]\begin{center}\includegraphics[#1]{#2.pdf}\end{center}\caption{#3}\label{fig:#2}\end{figure}}

\def\set5{[1:5]}

\newcommand{\dd}{\textbf{d}}
\newcommand{\ee}{\textbf{e}}
\newcommand{\ff}{\textbf{f}}
\newcommand{\vv}{\textbf{v}}

\newcommand{\sups}[1]{\textsuperscript{#1}}

\newcommand{\linit}{\textrm{left-init}}
\newcommand{\rinit}{\textrm{right-init}}
\newcommand{\lterm}{\textrm{left-term}}
\newcommand{\rterm}{\textrm{right-term}}
\newcommand{\wtA}{\widetilde{A}}

\newcommand{\aA}{\overrightarrow{A}}

\newcommand{\ovAo}{\overline{A_o}}
\newcommand{\ovAe}{\overline{A_e}}
\newcommand{\ovAos}{\overline{A_o^*}}
\newcommand{\ovAes}{\overline{A_e^*}}

\newcommand{\four}{\{1,2,3,4\}}

\newcommand{\conv}{\textrm{converging}}
\newcommand{\dive}{\textrm{diverging}}
\newcommand{\miss}{\textrm{missing}}

\def\Pro{P_{\mathrm{red}}^{\mathrm{out}}}
\def\Pbo{P_{\mathrm{blue}}^{\mathrm{out}}}
\def\Pri{P_{\mathrm{red}}^{\mathrm{in}}}
\def\Pbi{P_{\mathrm{blue}}^{\mathrm{in}}}
\def\Pr{P_{\mathrm{red}}}
\def\Pb{P_{\mathrm{blue}}}

\newcommand{\vinfty}{v_\infty}
\newcommand{\Gbox}{{G^{\boxtimes}}}

\newcommand{\plm}{planar map }
\newcommand{\plmm}{planar map}
\newcommand{\plms}{planar maps }
\newcommand{\plmms}{planar maps}

\begin{document}

\author[Olivier Bernardi, \'Eric Fusy, and Shizhe Liang]{Olivier Bernardi$^{*}$ \and \'{E}ric Fusy$^{\dagger}$ \and Shizhe Liang$^{+}$}
\thanks{$^{*}$Department of Mathematics, Brandeis University, Waltham MA, USA,
bernardi@brandeis.edu.\\
$^{\dagger}$LIGM/CNRS, Universit\'e Gustave Eiffel, Champs-sur-Marne, France, eric.fusy@univ-eiffel.fr.\\
$^{+}$Department of Mathematics, Brandeis University, Waltham MA, USA, Shizhe1011@brandeis.edu.\\
}

\title[A census of graph-drawing algorithms]{A census of graph-drawing algorithms based on generalized transversal structures}
\date{\today}

%

\begin{abstract}
We present two graph drawing algorithms based on the recently defined \emph{grand-Schnyder woods}, which are a far-reaching generalization of the classical Schnyder woods. The first is a straight-line drawing algorithm for plane graphs with faces of degree 3 and 4 with no separating 3-cycle, while the second is a rectangular drawing algorithm for the dual of such plane graphs.

In our algorithms, the coordinates of the vertices are defined in a global manner, based on the underlying grand-Schnyder woods. The grand-Schnyder woods and drawings are computed in linear time. 

When specializing our algorithms to special classes of plane graphs, we recover the following known algorithms:
\begin{compactitem} 
\item He's algorithm for rectangular drawing of 3-valent plane graphs, based on transversal structures,  
\item Fusy's algorithm for the straight-line drawing of triangulations of the square, based on transversal structures,
\item Bernardi and Fusy's algorithm for the orthogonal drawing of 4-valent plane graphs, based on 2-orientations,
\item Barri\`ere and Huemer's algorithm for the straight-line drawing of quadrangulations, based on separating decompositions.	 
\end{compactitem}
Our contributions therefore provide a unifying perspective on a large family of graph drawing algorithms that were originally defined on different classes of plane graphs and were based on seemingly different combinatorial structures.
\end{abstract}

\maketitle

\vspace{-7mm}

\section{Introduction}\label{sec:intro}



In the rich literature on graph drawing algorithms (see for instance~\cite{BattistaETT99,NishizekiR04} and references therein), two of the most classical and popular types of drawings are \emph{straight-line drawings} (in which edges are represented by straight-line  segments) and \emph{orthogonal drawings} (where edges are represented by a sequence of axis-aligned segments). 


\fig{width=.55\linewidth}{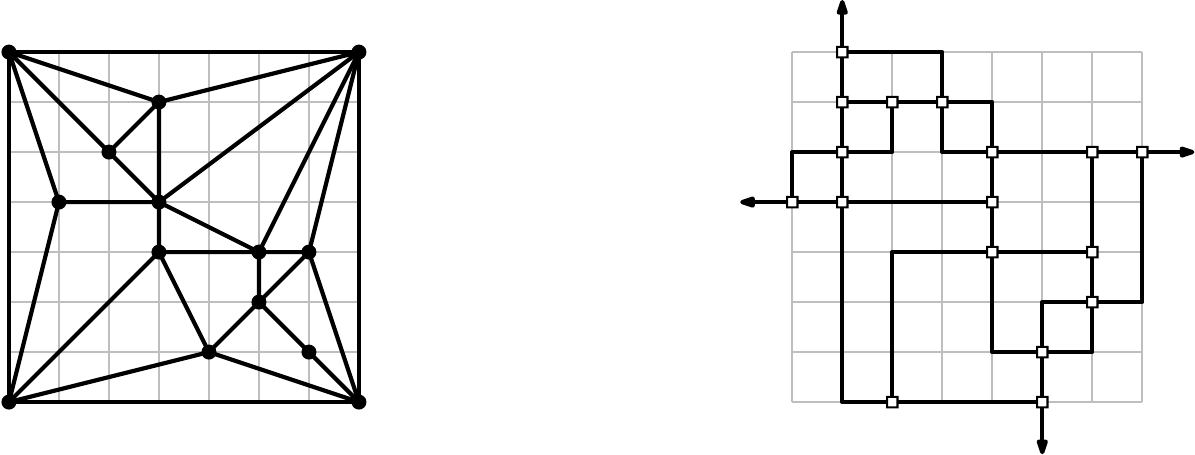}{Straight-line drawing of a plane graph with faces of degree~3 and~4 (left), and rectangular drawing of its dual (right).} 

\subsection{New drawing algorithms based on grand-Schnyder woods}
In this article we present two new drawing algorithms: a straight-line drawing algorithm for simple plane graphs with faces of degree 3 or 4 with no separating cycle of length $3$, and an orthogonal drawing algorithm for their duals. Graph drawings obtained using these algorithms are shown in Figure~\ref{fig:intro-pictures}. In these algorithms, the coordinates of the vertices are defined globally, either in terms of face-counts in certain regions or in terms of increasing functions over some bipolar orientations. The combinatorial structures and the drawings can be computed in linear time.

Our drawing algorithms are based on a recently discovered class of combinatorial structures~\cite{OB-EF-SL:Grand-Schnyder}. These structures are called \emph{grand-Schnyder structures} in honor of their relation to the classical \emph{Schnyder woods} underlying a most famous drawing algorithm~\cite{Schnyder:wood2}. Precisely, for every integer $d\geq 3$, a class of combinatorial structures called \emph{$d$-grand-Schnyder} (\emph{$d$-GS} for short), was defined for \emph{$d$-maps}, that is, connected plane graphs such that the outer face has degree $d$ and inner faces have degree at most $d$. It was shown in~\cite{OB-EF-SL:Grand-Schnyder} that a $d$-map admits a $d$-GS structure if and only if its non-facial cycles have length at least $d$. The classical Schnyder woods correspond to 3-GS structures on triangulations. 

In this article, the algorithms are based on 4-GS structures. 
Interestingly, when specializing 4-GS structures to the case of the triangulations of the square, one recovers a combinatorial structure used in another classical graph drawing algorithm, namely the \emph{regular edge labelings} used by He to define a rectangular drawing algorithm for 3-valent plane graphs~\cite{He93:reg-edge-labeling}. Since then, regular edge labelings have been rediscovered and popularized under the name \emph{transversal structures} which we shall adopt in the present article. In a different direction, when specializing 4-GS structures to quadrangulations one recovers the separating decompositions (or equivalently, 2-orientations) which have also found applications to graph drawing algorithms: straight-line drawings and orthogonal drawings~\cite{Barriere-Huemer:4-Labelings-quadrangulation,OB-EF:Schnyder}, as well as visibility representations and segment contact representations~\cite{RosenstiehlT86,TamassiaT86,TamassiaTollis89}. 

As we will later explain, our orthogonal  drawing algorithm interpolates between He's algorithm~\cite{He93:reg-edge-labeling} (for 3-valent plane graphs) and Bernardi and Fusy's algorithm~\cite{OB-EF:Schnyder} (for 4-valent plane graphs), while our straight-line drawing algorithm interpolates between Fusy's algorithm ~\cite{Fu07b} (for triangulations) and Barri\`ere and Huemer's algorithm~\cite{Barriere-Huemer:4-Labelings-quadrangulation} (for quadrangulations). This reveals the deep hidden commonality between these 4 algorithms (in addition to relating them to other Schnyder-wood based algorithms~\cite{Schnyder:wood2, Felsner:woods}). Our proofs are also uniform, and make use of a simple local planarity criteria which could be of independent interest (Lemma~\ref{lem:planar-criteria}).

Before describing our results further, we give a quick review of the 4 drawing algorithms in~\cite{Barriere-Huemer:4-Labelings-quadrangulation,OB-EF:Schnyder,Fu07b,He93:reg-edge-labeling}.

\subsection{A review of four classical graph-drawing algorithms}
In this subsection we review the following graph drawing algorithms:
\begin{itemize}
\item the straight-line drawing algorithm for irreducible triangulations of the square due to Fusy~\cite{Fu07b}, and the (bendless) orthogonal drawing algorithm for their duals due to He~\cite{He93:reg-edge-labeling},
\item the straight-line drawing algorithm for simple quadrangulations due to Barri\`ere and Huemer~\cite{Barriere-Huemer:4-Labelings-quadrangulation}, and the orthogonal drawing algorithm for their duals due to Bernardi and Fusy~\cite{OB-EF:Schnyder}. 
\end{itemize}
%
%
We start by defining the combinatorial structures underlying these algorithms: transversal structures and separating decompositions.

\fig{width=0.9\linewidth}{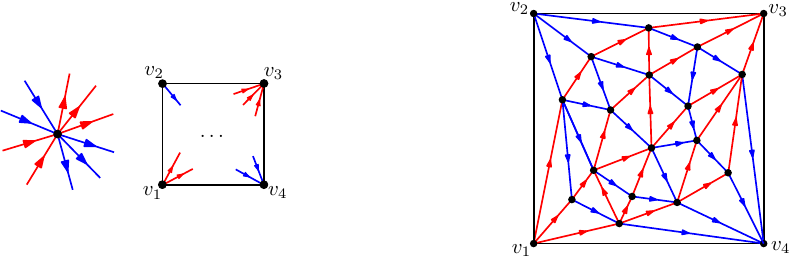}{Left: The local conditions (at inner vertices and at outer vertices) for transversal structures. Right: Example of a transversal structure on a triangulation of the square.}

We call \emph{\plm} a planar drawing of a connected planar graph, considered up to continuous deformation.
A \emph{triangulation of the square} is a \plm such that every inner face has degree 3, while the outer face has degree 4. Such a \plm is called \emph{irreducible} if it has no cycle of length less than 4 except for the contour of the inner faces.
It was shown by He~\cite{He93:reg-edge-labeling} that irreducible triangulations of the square can be endowed with some combinatorial structures that he named \emph{regular edge-labelings}. These structures were later rediscovered by Fusy~\cite{Fu07b}, who used the name \emph{transversal structures} that we adopt below.
Consider an irreducible triangulation of the square $G$, with its outer vertices labeled $v_1,v_2,v_3,v_4$ in clockwise order.
A \emph{transversal structure} for $G$ is an orientation and bicoloration of its inner edges --- in blue and red --- satisfying the local conditions shown in Figure~\ref{fig:transversal}: 
\begin{enumerate}
\item
At each inner vertex the incident edges form 4 groups in clockwise order: outgoing red, outgoing blue, incoming red, incoming blue. 
\item
The edges incident to the outer vertices $v_1,v_2,v_3,v_4$ are respectively outgoing red, outgoing blue, incoming red, and incoming blue.
\end{enumerate}
It is known that these local conditions imply the global condition that the spanning subgraph $G_r$ formed by the red edges and the 4 outer edges directed from $v_1$ to $v_3$ is an acyclic orientation with a unique source $v_1$ and a unique sink $v_3$; and similarly the spanning subgraph $G_b$ formed by the blue edges and the 4 outer edges directed from $v_2$ to $v_4$ is an acyclic orientation 
with a unique source $v_2$ and a unique sink $v_4$. See Figure~\ref{fig:transversal} for an example.

\fig{width=\linewidth}{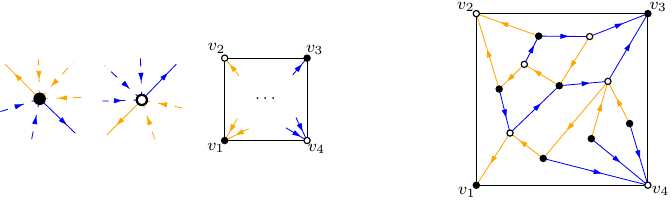}{Left: The local conditions, at inner vertices and at outer vertices, for separating decompositions. Right: Example of a separating decomposition on a simple quadrangulation.}

A \emph{(simple) quadrangulation} is a \plm (without loops nor double edges) such that every face has degree 4. 
It was shown in~\cite{de1995bipolar} that simple quadrangulations can be endowed with combinatorial structures called \emph{separating decompositions}. 
Consider a simple quadrangulation $G$ with its outer vertices labeled $v_1,v_2,v_3,v_4$ in clockwise order, and consider the proper bicoloration of its vertices in black and white such that $v_1$ is black (quadrangulations are bipartite). 
A \emph{separating decomposition} for $G$ is an orientation and bicoloration of its inner edges --- in blue and orange --- satisfying the local conditions shown in Figure~\ref{fig:separating_decomp}:
\begin{enumerate}
\item
At every black (resp. white) inner vertex there is a single outgoing blue edge $e$ and a single outgoing orange edge $e'$, and all the incident edges between $e$ and $e'$ in clockwise order are blue (resp. orange) while all the incident edges between $e'$ and $e$ are orange (resp. blue).
\item 
The edges incident to the outer vertices $v_1,v_2,v_3,v_4$ are all incoming; they are orange at $v_1,v_2$ and blue at $v_3,v_4$.
\end{enumerate} 
It is shown in~\cite{de1995bipolar} that the subgraph $T_o$ formed by the inner vertices, the orange edges, and the outer edge $(v_1,v_2)$ is a tree, with edges oriented toward $\{v_1,v_2\}$; by convention it is considered as rooted at $v_2$. Similarly the subgraph $T_b$ formed by the inner vertices, the blue edges and the outer edge $(v_3,v_4)$ is a tree, rooted at $v_4$.

The two combinatorial structures described above have been used to define straight-line drawing algorithms using face-counting operations. These algorithms were introduced respectively in~\cite{Fu07b} for transversal structures, and in~\cite{Barriere-Huemer:4-Labelings-quadrangulation} for separating decompositions. We will now describe these algorithms before turning to the dual picture.

\fig{width=\linewidth}{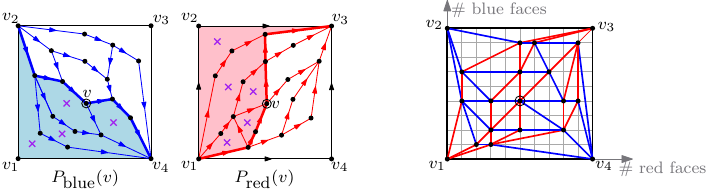}{The straight-line drawing obtained using the transversal structure shown in Figure~\ref{fig:transversal}. The blue and red separating paths are shown for the circled vertex $v$. For this vertex, the blue path $P_b(v)$ has 4 faces on its right in the blue map, and the red path has 5 faces on its left in the red map. Accordingly, $v$ has coordinates $(5,4)$ in the drawing. Placing similarly all the inner vertices yields the planar straight-line drawing shown on the right. Since the red and blue maps have 10 and 9 inner faces respectively, the grid has size $10\times 9$.}

We start by describing the drawing algorithm from~\cite{Fu07b} which is illustrated in Figure~\ref{fig:draw_transversal}. Let $G$ be an irreducible triangulation of the square endowed with a transversal structure. For every inner vertex $v$ of $G$ we will define some blue and red \emph{separating paths}, represented in Figure~\ref{fig:draw_transversal}, which are used to determine the placement of $v$.
The \emph{leftmost outgoing blue path} of an inner vertex $v$ is the path $\Pbo(v)$ starting at $v$ and taking the leftmost outgoing blue edge at each step until reaching $v_4$, while the \emph{rightmost incoming blue path} is the path $\Pbi(v)$ starting at $v$ and taking the rightmost incoming blue edge at each step until reaching $v_2$. We use the notation $P^-$ for the reversal of a directed path $P$, and define the \emph{blue separating path} at $v$ as the blue path $\Pb(v)=\Pbo(v)\cup \Pbi(v)^-$ going from $v_2$ to $v_4$. 
We define similarly the \emph{leftmost outgoing red path} $\Pro(v)$ (going from $v$ to $v_3$), the \emph{rightmost incoming red path} $\Pri(v)$ (going from $v$ to $v_1$), and the \emph{red separating path} $\Pr(v)=\Pro(v)\cup \Pri(v)^-$.
In the straight-line drawing of $G$ given by~\cite{Fu07b}, the vertices are placed on a grid of size $f_r\times f_b$, where $f_r,f_b$ are the numbers of inner faces in the red and blue maps $G_r,G_b$:
\begin{compactitem}
\item the outer vertices of $G$ are placed at the 4 corners of the grid with $v_1$ at the origin, 
\item each inner vertex $v$ is placed at $(x_v,y_v)$, where $x_v$ is the number of inner faces of $G_r$ on the left of $\Pr(v)$, and $y_v$ is the number of inner faces of $G_b$ on the right of $\Pb(v)$. 
\end{compactitem}
An example is given in Figure~\ref{fig:draw_transversal}.
 
 \fig{width=\linewidth}{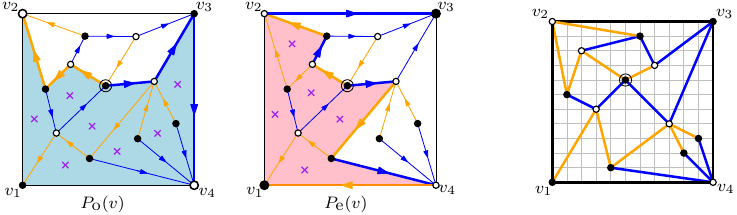}{The straight-line drawing obtained using the separating decomposition shown in Figure~\ref{fig:separating_decomp}. 
 The even (resp. odd) separating path of the circled vertex is shown,
with 5 faces on its left (resp. $7$ faces on its right). Accordingly the circled vertex has coordinates $(5,7)$ in the drawing. Placing 
similarly all the inner vertices yields the planar straight-line drawing shown on the right. Since there are 11 inner faces, the grid has size $11\times 11$.}

Next we describe the algorithm from~\cite{Barriere-Huemer:4-Labelings-quadrangulation}, which is illustrated in Figure~\ref{fig:str_dr_separating}. 
Let $G$ be a simple quadrangulation endowed with a separating decomposition.
For each inner vertex $v$ we define an odd and an even separating path as illustrated in Figure~\ref{fig:str_dr_separating}. For all $i\in\four$ we define $P_i(v)$ as the directed path from $v$ to $v_{i+3}$ obtained by following the outgoing arcs as follows: 
\begin{compactitem}
\item for $i=1$ (resp. $i=3$) we follow the outgoing blue (resp. orange) arc at every vertex, 
\item for $i=2$ (resp. $i=4$) we follow the outgoing blue (resp. orange) arc at every black vertex, and the outgoing orange (resp. blue) arc at every white vertex. 
\end{compactitem}
In the above definition of $P_i(v)$ the arcs $(v_2,v_3)$ and $(v_3,v_4)$ are blue, and the arcs $(v_4,v_1)$ and $(v_1,v_2)$ are orange. It is shown in~\cite{FeHuKa} that with these conventions, the path $P_i(v)$ is always well-defined and ends at $v_{i+3}$ for all $i$~\cite{FeHuKa}.
The \emph{odd separating path} at $v$ is the path $P_o(v):=P_1(v)\cup P_3(v)^-$ going from $v_2$ to $v_4$, while the \emph{even separating path} at $v$ is the path $P_e(v):=P_4(v)\cup P_2(v)^-$ going from $v_1$ to $v_3$. The separating paths are always simple paths as shown in~\cite{FeHuKa}.
In the straight-line drawing of $G$ given by~\cite{Barriere-Huemer:4-Labelings-quadrangulation}, the vertices are placed on a grid of size $f\times f$, where $f$ is the number of inner faces:
\begin{compactitem}
\item the outer vertices of $G$ are placed at the 4 corners of the grid with $v_1$ at the origin, 
\item each inner vertex $v$ is placed at $(x_v,y_v)$, where $x_v$ is the number of inner faces of $G$ on the left of $P_e(v)$, and $y_v$ is the number of inner faces of $G$ on the right of $P_o(v)$. 
\end{compactitem}
An example is given in Figure~\ref{fig:str_dr_separating}. 

Let us highlight the similarity between the algorithms in~\cite{Fu07b} and~\cite{Barriere-Huemer:4-Labelings-quadrangulation}. As already noted in~\cite{Barriere-Huemer:4-Labelings-quadrangulation} (building upon the \emph{0,1-labelings} of~\cite{FeHuKa}), there is a natural way to label the inner corners of a separating decomposition with labels in $\four$ using the rule represented in Figure~\ref{fig:labeling-separating}, and the paths $P_1(v),\ldots,P_4(v)$ associated to an inner vertex $v$ can be defined in terms of this labeling: at each inner vertex the path $P_i(v)$ takes the outgoing edge separating the labels $i$ and $i+1$. Similarly, there is a natural way to label the inner corners of a transversal structure using the rules represented in Figure~\ref{fig:labeling-separating}, and the blue and red separating paths can be defined in terms of these labelings: at each inner vertex the path $\Pbo(v)$ (resp. $\Pri(v)$, $\Pbi(v)$, $\Pro(v)$) takes the edge separating the labels 1 and 2 (resp. 2 and 3, 3 and 4, 4 and 1).

\fig{width=\linewidth}{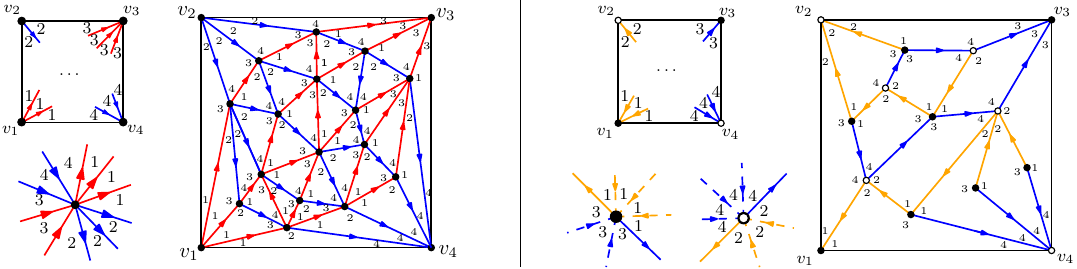}{Left: Corner labeling associated to a transversal structure (rule for defining the labels around outer and inner vertices, and an example). Right: Corner labeling associated to a separating decomposition (rule for defining the labels around outer vertices and around black and white inner vertices, and an example).}

We now turn our attention to the dual picture, and review two orthogonal drawing algorithms based on the same pair of combinatorial structures~\cite{OB-EF:Schnyder,He93:reg-edge-labeling}. An \emph{orthogonal drawing} of a planar graph $G'$ is a drawing such that each edge is a sequence of axis-aligned segments, a bend being a turning point on an edge. An orthogonal drawing is called \emph{one-bend} if every edge has at most one bend. All orthogonal drawing algorithms to be reviewed and introduced in the paper produce one-bend drawings  
(upon not drawing the root vertex, see below). 

 \fig{width=\linewidth}{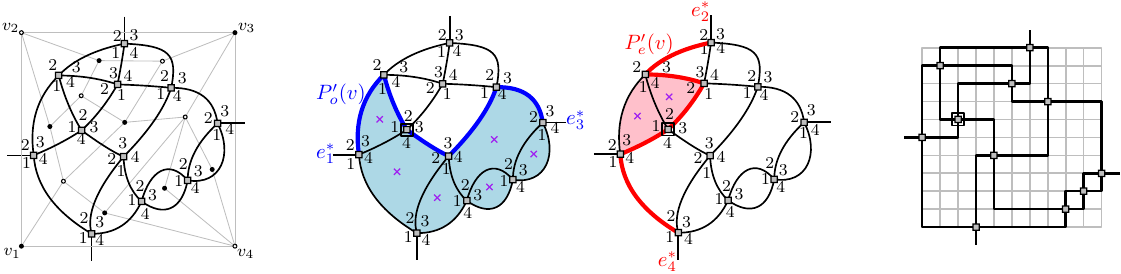}{Orthogonal drawing from~\cite{OB-EF:Schnyder} obtained using the (dual of the) separating decomposition shown in Figure~\ref{fig:separating_decomp}. 
 The odd and even separating paths of the circled vertex $v$ are indicated. The even path $P_e'(v)$ has 2 faces on its left, while the odd path $P_o'(v)$ has $6$ faces on its right, hence the circled vertex has coordinates $(2,6)$ in the drawing. Placing 
similarly all the dual vertices (excluding $\vinfty$) yields the orthogonal drawing shown on the right.}

We first describe orthogonal drawing algorithm from~\cite{OB-EF:Schnyder} for the dual of simple quadrangulations.  
Let $G$ be a simple quadrangulation, let $G^*$ be the dual of $G$, and let $\vinfty$ be the vertex of $G^*$ dual to the outer face of $G$. We call $\vinfty$ the \emph{root vertex} of $G^*$. The drawing of $G':=G^*\setminus \vinfty$ is again based on a face-counting procedure.
For $i\in\four$ we define $e_i^*$ as the edge of $G^*$ dual to the outer edge $\{v_i,v_{i+1}\}$ of $G$, and we let $v_i^*$ be the non-root endpoint of $e_i^*$.
Recall that, given a separating decomposition of $G$, one can associate some labels to the inner corners of $G$ as indicated in Figure~\ref{fig:labeling-separating}. By duality, this gives a labeling of the corners of $G^*$ not incident to $\vinfty$. 
This labeling is such that the labels $1,2,3,4$ appear in clockwise order around each non-root vertex. See Figure~\ref{fig:BeFu} for an example. 
For $i\in\{1,2,3,4\}$, let $P_i'(v)$ be the path starting at $v$ and taking the edge separating the corners of labels $i$ and $i+1$ at each vertex, until reaching $v_i^*$. 
The \emph{even separating path} at $v$ is the simple directed path $P_e'(v):=P_2'(v)\cup P_4'(v)^-$ from $v_4^*$ to $v_2^*$ (where we again use the notation $P^-$ to indicate the path obtained by reversing $P$). The  \emph{odd separating path} at $v$ is the path $P_o'(v):=P_3'(v)\cup P_1'(v)^-$ from $v_1^*$ to $v_3^*$.
In the orthogonal drawing of $G'$ from~\cite{OB-EF:Schnyder}, the coordinates $(x,y)$ of each non-root vertex $v$ are such that $x$ (resp. $y$) is the number of inner faces of $G'$ on the left of $P_e'(v)$ (resp. on the right of $P_o'(v)$). 
Each edge is then drawn as the union of two axis-aligned segments in the unique way such that, at each vertex~$v$, the path $P_1'(v)$ (resp. $P_2'(v)$, $P_3'(v)$, $P_4'(v)$) leaves $v$ toward the left (resp. top, right, bottom). Equivalently, the corner of label 1 (resp. 2,3,4) occupies the south-west (resp. north-west, north-east, south-east) quadrant at $v$.
If $G'$ has $n$ inner faces, then the grid size is $n\times n$. An example is given in Figure~\ref{fig:BeFu}.

 \fig{width=\linewidth}{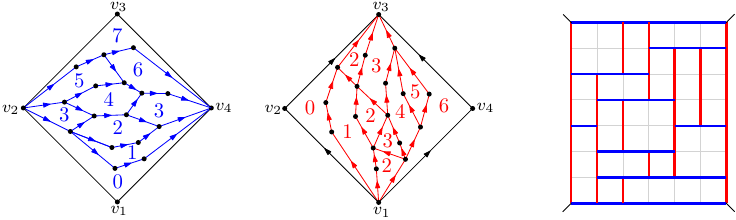}{Rectangular drawing obtained using increasing functions on the transversal structure shown in Figure~\ref{fig:transversal}. 
 The increasing functions of $G_b$ and $G_r$ are indicated on the left, and the resulting rectangular drawing of $G'$ is shown on the right.}

Last but not least, we recall the rectangular drawing algorithm from~\cite{He93:reg-edge-labeling} based on transversal structures.
Recall that a \emph{rectangular drawing} is an orthogonal drawing, such that edges have no bend, and all faces (including the 
outer face) have a rectangular shape. Let $G$ be an irreducible triangulation of the square, let $G^*$ be its dual, and let $\vinfty$ be the vertex of $G^*$ dual to the outer face of $G$. We consider a transversal structure of $G$, and define a rectangular drawing of $G':=G^*\setminus \vinfty$ using some ``increasing functions'' on the faces of the blue and red maps, as represented in Figure~\ref{fig:He_drawing}.
An \emph{increasing function} for $G_r$ (resp. $G_b$) is an assignment of values in $\mathbb{R}$ to the inner faces such that every inner edge has a smaller (resp. larger) label on the left side than on the right side. 
Every vertex $v\in G'$ is dual to a triangular inner face $f$ of $G$ which is contained in a face $f_r(v)$ of $G_r$ and in a face $f_b(v)$ of $G_b$.
In the rectangular drawing of $G'$ from~\cite{He93:reg-edge-labeling}, the coordinates $(x,y)$ of $v$ are such that $x$ (resp. $y$) is the label of the face $f_r(v)$ (resp. $f_b(v)$). If face-labels are in $\mathbb{N}$, then one obtains a rectangular drawing on a grid. 
The choice of increasing function for $G_r$ (resp. $G_b$) yielding the smallest grid-width (resp. grid-height) is the one where the label of the leftmost (resp. rightmost) inner face is $0$, and the label of each other inner face is $1$ plus the maximum of the labels of the faces that are adjacent from the left (resp. right); equivalently, each face is labeled by the length of a longest path in the dual bipolar orientation. This optimal choice of increasing functions is used in Figure~\ref{fig:He_drawing}.

\subsection{Outline of the article}
In this article we present two new drawing algorithms that make use of the \emph{grand-Schnyder woods} introduced in~\cite{OB-EF-SL:Grand-Schnyder}.
Our straight line drawing algorithm applies to the class of \emph{4-adapted maps}, which are the \plms such that the outer face has degree 4, the inner faces have degree at most 4 and the non-facial cycles have length at least 4. Our orthogonal drawing algorithm applies to the dual of 4-adapted maps.


It was shown in~\cite{OB-EF-SL:Grand-Schnyder} that any 4-adapted map admits a \emph{4-grand-Schnyder wood} (or\emph{4-GS wood} for short), which can be computed in linear time. The definition of 4-GS woods will be given in Section~\ref{sec:combinatorics}.
Note that the irreducible triangulations of the square are the 4-adapted maps having only inner faces of degree $3$, while the simple quadrangulations are the 4-adapted maps having inner faces of degree $4$. Hence, irreducible triangulations of the square and simple quadrangulations form two subclasses of 4-adapted maps. 
Moreover, as we will explain in Section~\ref{sec:4GS1}, the transversal structures of a triangulation of the square are in bijection with the $4$-GS woods of that triangulation. 
As for separating decompositions, they correspond bijectively to a special class of $4$-GS woods of quadrangulations satisfying some parity conditions (the \emph{even 4-GS woods}). So the concept of 4-GS woods allows one to interpolate between the concept of transversal structures used in~\cite{Fu07b,He93:reg-edge-labeling} and the concept of separating decompositions used in~\cite{Barriere-Huemer:4-Labelings-quadrangulation,OB-EF:Schnyder}.

Regarding straight-line drawings, we present a ``face-counting'' algorithm for 4-adapted maps, where the coordinates of each vertex are determined as the number of faces in certain ``regions'' associated to that vertex. 
When specialized to irreducible triangulations of the square or to simple quadrangulations, the face counting algorithm yields the two algorithms from~\cite{Fu07b} and from~\cite{Barriere-Huemer:4-Labelings-quadrangulation} respectively.\footnote{For triangulations of the square, we more precisely recover the algorithm from~\cite{Fu07b} in its compact version (where every grid-line is occupied by at least one vertex).} 
We then present a more general straight-line drawing algorithm for 4-adapted maps, where the coordinates of each vertex are obtained using \emph{increasing functions} on some bipolar orientations associated to the 4-GS wood. 
We show that the increasing function algorithm can be performed in linear time, and that it leads to a grid size that is less or equal to the grid size of the face-counting algorithm. We also discuss some strategies for further decreasing the grid-size. 

Regarding orthogonal drawings, we first give a ``face-counting'' algorithm for the dual of 4-adapted maps (where the coordinates of each vertex are defined as the number of faces in certain ``regions'' associated to that vertex). We then present a more general ``increasing function" drawing algorithm. The face-counting algorithm recovers the algorithm from~\cite{OB-EF:Schnyder} when applied to the dual of quadrangulations, while the increasing-function algorithm recovers the algorithm from~\cite{He93:reg-edge-labeling} when applied to the dual of triangulations.
We show that the increasing function algorithm can be performed in linear time, and we then discuss strategies to minimize the number of bends and to decrease the grid-size.
The orthogonal drawings produced by our algorithm are \emph{right-chiral}: an edge starting at a 3-valent vertex is either straight or has a right-bend (in particular edges between 3-valent vertices are straight). We show that the dual of a 4-map admits a right-chiral drawing if and only if it is 4-adapted, in which case we give an algorithm to construct a right-chiral drawing achieving the minimum number of bends possible.


The article is organized as follows. In Section~\ref{sec:combinatorics} we define the combinatorial structures relevant for our drawing algorithms: 4-GS labelings, 4-GS woods, and some associated bipolar orientations (in the primal and in the dual). In Section~\ref{sec:primal_algo} we present our straight-line drawing algorithms for 4-adapted maps: we start by presenting a face-counting algorithm (Section~\ref{sec:face-counting-primal}) before presenting an increasing-function algorithm (Section~\ref{sec:increasing-primal}). We then prove the correctness of the algorithms (based on a local correctness criterion, Lemma~\ref{lem:planar-criteria}), analyze the time-complexity and grid size, and present an optimization for reducing the grid size (Sections~\ref{sec:proof-planar-primal} to~\ref{sec:optimization}). In Section~\ref{sec:dual_algo} we present our orthogonal drawing algorithms for the dual of 4-adapted maps. We start by presenting a face-counting algorithm (Section~\ref{sec:face-counting-dual}), before presenting an increasing-function algorithm (Section~\ref{sec:increasing_dual}). We prove the correctness of the algorithms in Section~\ref{sec:planarity_ortho}, before discussing the relation with He's algorithm, discussing additional properties of the drawing, explaining how to deal with vertices of degree 2, analyzing the time-complexity and the grid size, and presenting optimizations for reducing the grid size (Sections~\ref{sec:comparison-He} to~\ref{sec:reduce-grid-dual}).
Finally in Section~\ref{sec:conclusion} we present some variations on the face-counting algorithms (Section~\ref{sec:variant-face-count}), we discuss a way to triangulate 4-adapted maps while preserving 4-adaptedness (Section~\ref{sec:triangulating}), and we give an algorithm for minimizing the number of bends when drawing the duals of 4-adapted maps (Section~\ref{sec:bend}).





\section{Combinatorial structures}\label{sec:combinatorics}

In this section we lay out some basic definitions about \plmms, 
and then define the combinatorial structures which underlie our drawing algorithms.
These combinatorial structures have been introduced in \cite{OB-EF-SL:Grand-Schnyder} under the name \emph{grand-Schnyder structures}. Precisely, we will be using 4-grand-Schnyder structures and their duals.

\subsection{Planar maps}
Our \emph{graphs} are finite and undirected. 
A \emph{simple graph} is a graph without loops nor multiple edges. A \emph{\plm} is a connected planar graph drawn in the plane without edge crossing, considered up to continuous deformation.  A \emph{corner} of a \plm is the sector delimited by 2 consecutive edges around a vertex. A vertex or face has degree $k$ if it is incident to $k$ corners.
An \emph{arc} of a graph is an edge together with a chosen orientation of this edge (so each edge gives rise to two \emph{opposite} arcs). The endpoints of an arc $a$ are called the \emph{initial} and \emph{terminal} vertices of $a$. We write $a=(u,v)$ to indicate that $u$ and $v$ are the initial and terminal vertices of the arc $a$ respectively. The infinite face of a \plm is called the \emph{outer face}, while the others are called \emph{inner faces}. We also call \emph{outer} the vertices, edges, arcs and corners incident to the outer face, and \emph{inner} the others. 

A \emph{3,4-angulation of the square} is a \plm such that the inner faces have degree 3 or~4, and the outer face has degree 4 and is incident to 4 distinct outer vertices. We canonically denote the outer vertices of such a \plm by $v_1,v_2,v_3,v_4$, in clockwise order around the outer face.
A 3,4-angulation of the square is called \emph{adapted} if any simple cycle which is not the contour of a face has length at least 4. Note that adapted 3,4-angulations are necessarily simple. Note also that a \plm with every face of degree 4 is adapted if and only if its girth is 4 (where \emph{girth} means the minimal length of cycles).

Recall that the \emph{dual} $G^*$ of a \plm $G$ is obtained as follows: 
\begin{compactitem}
\item we draw a vertex $v_f$ of $G^*$ inside every face $f$ of $G$ and call $v_f$ the \emph{dual vertex} of $f$, 
\item and we draw an edge $e^*$ of $G^*$ connecting $v_f$ and $v_g$ across each edge $e$ separating the faces $f$ and $g$ of $G$, and call $e^*$ the \emph{dual edge} of the \emph{primal edge} $e$.
\end{compactitem}
The \emph{root vertex} of $G^*$ is the vertex $\vinfty$ dual to the outer face of $G$.

If $G$ is a 3,4-angulation of the square, then its dual graph $G^*$ is a \plm having vertices of degree 3 or~4, with a designated root vertex $\vinfty$ of degree~4. We call such a map a \emph{rooted 3,4-map}. A rooted 3,4-map $G^*$ is called \emph{dual-adapted} if $G$ is adapted. Equivalently, a rooted 3,4-map $G^*$ is dual-adapted if the deletion of less than 4 edges does not disconnect $G^*$ into 2 subgraphs having more than 1 vertex each. Note that dual-adaptedness is less demanding than \emph{4-edge connectivity} since vertices of degree 3 are allowed.

\subsection{Grand-Schnyder structures}\label{sec:4GS1}
In this subsection we recall some definitions and results from \cite{OB-EF-SL:Grand-Schnyder} about 4-GS labelings and 4-GS woods.
 
Let $G$ be a 3,4-angulation of the square. A \emph{4-labeling} of $G$ is an assignment of a \emph{label} in $\{1,2,3,4\}$ to each inner corner of $G$; see for example Figure~\ref{fig:4-grand-Schnyder}(a). 
A $4$-GS labeling of $G$ is a $4$-labeling satisfying certain local conditions represented in Figure~\ref{fig:def-4GS-labeling}.
These conditions are best expressed in terms of \emph{jumps}. Consider a $4$-labeling of $G$. Let $c$ and $c'$ be two inner corners of $G$, and let $i$ and $i'$ be their respective labels. The \emph{label jump} from the corner $c$ to the corner $c'$ is defined as the integer $\delta$ in $\{0,1,2,3\}$ such that $i+\delta\equiv i'$ modulo 4 (in other words, the label jump $\delta$ is $i'-i$ if $i'-i\geq 0$, and $i'-i+4$ otherwise).

For an inner vertex $v$ of $G$, the \emph{sum of clockwise jumps} around~$v$ is the sum of label jumps between consecutive corners in clockwise order around~$v$. Similarly, the \emph{sum of clockwise jumps} around a face $f$ is the sum of label jumps between consecutive corners in clockwise order around~$f$.

\begin{definition}\label{def:GS-labeling}
Let $G$ be a 3,4-angulation of the square. As usual, we assume that the outer vertices of $G$ are denoted by $v_1,\ldots,v_4$ in clockwise order around the outer face.
A \emph{$4$-GS labeling} of $G$ is a 4-labeling satisfying the following conditions.
\begin{itemize}
\item[(L0)] For all $i\in \{1,2,3,4\}$, all the inner corners incident to $v_i$ have label $i$.
\item[(L1)] For every inner vertex or face of $G$, the sum of clockwise jumps around this vertex or face is $4$.
\item[(L2)] The label jump from a corner to the next corner around a face are always strictly positive.
\item[(L3)] Let $e$ be an inner edge incident to a face $f$ of degree 3. Let $c$ and $c'$ be the consecutive corners incident to $e$ in clockwise order around~$f$, and let $v$ be the vertex incident to~$c'$. The label jump $\delta$ from $c$ to $c'$ and the label jump $\eps$ from $c'$ to the next corner in clockwise order around~$v$ satisfy $\delta+\eps\geq 2$. 
\end{itemize}
\end{definition}

\fig{width=\linewidth}{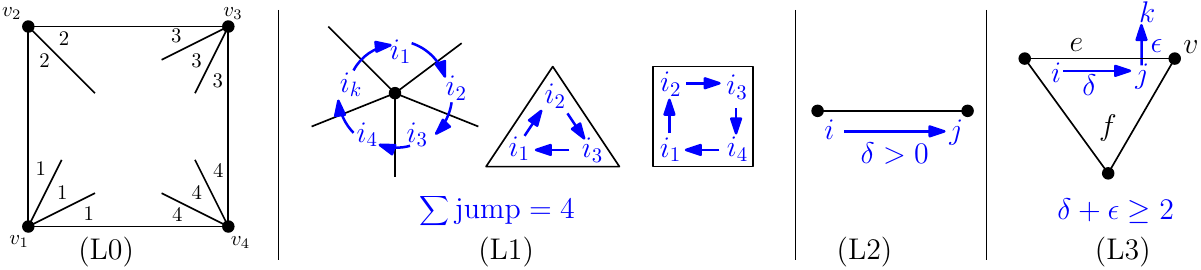}{Conditions defining 4-GS labelings.}

It was shown in \cite{OB-EF-SL:Grand-Schnyder} that a 3,4-angulation of the square $G$ admits a 4-GS labeling if and only if it is adapted. Moreover, in this case, a 4-GS labeling can be constructed in linear time (number of operations linear in the number of vertices of $G$). Let us state a few additional properties of 4-GS labelings.
\begin{remark}
\begin{compactitem}
\item Conditions (L1) and (L2) imply that around any inner face of degree 4 the labels are 1,2,3,4 in clockwise order. For an inner face of degree $3$, the labels are $i,i+1,i+2$ (modulo 4) in clockwise order for some $i\in\{1,2,3,4\}$. In particular the label jumps along edges are either 1 or 2.
\item The sum of label jumps counterclockwise around any inner edge $e$ is equal to~$4$. This property has been established in \cite{OB-EF-SL:Grand-Schnyder}.
Note that this implies that the label jumps around vertices are either 0, 1 or 2.
\end{compactitem}
\end{remark}

We now come back to the relation between transversal structures, separating decompositions, and 4-GS structures. As explained in the introduction, there is a natural way to associate a 4-labeling to a transversal structure, or to a separating decomposition. We start with the case of transversal structures. Let $G$ be an irreducible (i.e. adapted) triangulation of the square. 
Given a transversal structure $\cT$ of $G$, we define a labeling $\cL$ of the inner corners of $G$ by the rule represented in Figure \ref{fig:labeling-separating}(left):
\begin{compactitem}
\item for all $i\in\four$, the corners incident to the outer vertex $v_i$ are labeled $i$,
\item a corner incident to an inner vertex $v$ is labeled 1 (resp. 2,3,4) if the edge preceding $c$ clockwise around $v$ is outgoing red (resp. outgoing blue, incoming red, incoming blue). 
\end{compactitem}
It was shown in~\cite{OB-EF-SL:Grand-Schnyder} that this gives a bijection between the set of transversal structures of $G$ and the set of 4-GS labelings of $G$. 
We now turn to separating decompositions. Let $G$ be a simple (i.e. adapted) quadrangulation with its outer vertices labeled $v_1,\ldots,v_4$, and its inner vertices properly bicolored in black and white such that $v_1$ is black.
Given a separating decomposition $\cS$ of $G$, we define a labeling of the inner corners of $G$ by the rule represented in Figure \ref{fig:labeling-separating}(right):
\begin{compactitem}
\item for all $i\in\four$, the corners incident to the outer vertex $v_i$ are labeled $i$,
\item a corner incident to a black (resp. white) inner vertex $v$ is labeled 1 (resp. 4) if it is between the outgoing orange  arc and the outgoing blue arc clockwise around $v$, and it is labeled 3 (resp. 2) otherwise. 
\end{compactitem}
This way of defining a corner labeling was given in~\cite{Barriere-Huemer:4-Labelings-quadrangulation}. As shown in~\cite{OB-EF:Schnyder,OB-EF-SL:Grand-Schnyder}, it gives a bijection between the separating decompositions of $G$ and the \emph{even} 4-GS labelings of $G$, that is, the 4-GS labelings such that the labels around black vertices are odd and the labels around white vertices are even.

\fig{width=\linewidth}{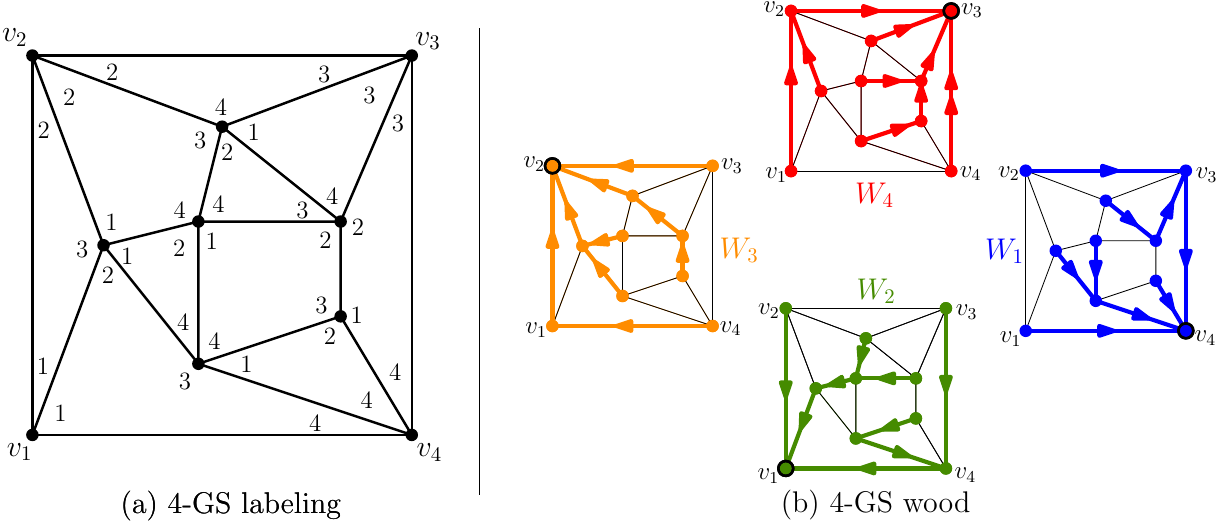}{(a) A 4-GS labeling $\cL$ of an adapted 3,4-angulation of the square. (b) The corresponding 4-GS wood $\cW=(W_1,W_2,W_3,W_4)$.}

Next, we define the 4-GS wood associated to a 4-GS labeling.
Let $G$ be a 3,4-angulation of the square. For each arc $a$ of $G$ there are four corners incident to $a$. 
The corner at the left of an arc $a$ around its initial vertex is called the \emph{left-initial corner of $a$}. The \emph{right-initial}, \emph{left-terminal}, and \emph{right-terminal} corners of $a$ are defined similarly. Given a 4-GS labeling $\cL$ of $G$, and an inner arc $a$, we denote by $\linit(a)$ (resp. $\rinit(a)$, $\lterm(a)$, $\rterm(a)$), the label of the left-initial (resp. right-initial, left-terminal, right-terminal) corner of $a$. This is represented in Figure~\ref{fig:labels-around-edge}(a).
For $i,j$ in $\four$, we denote by $[i:j[$ the set of integers $\{i,i+1,\ldots,j-1\}$ modulo 4. Concretely, if $i=j$ then $[i:j[=\emptyset$; if $i<j$ then $[i:j[=\{k\mid i\leq k<j\}$; and if $i>j$ then $[i:j[=\{k\mid i\leq k\leq 4\textrm{ or } 1\leq k<j\}$.

\fig{width=\linewidth}{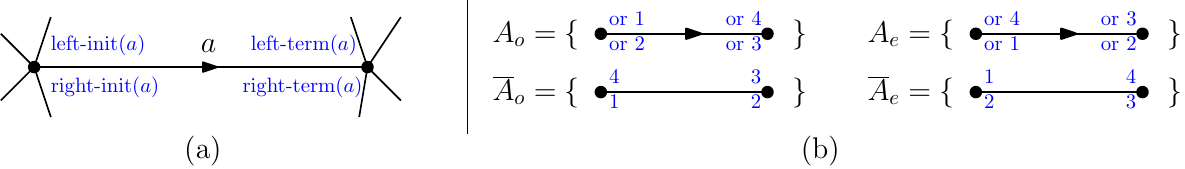}{(a) The labels around an arc. (b) The inner arcs in the acyclic orientation $A_o$ and $A_e$, and the inner edges not in their support.}

Given a 4-GS labeling $\cL$ of $G$, we define four oriented subgraphs $W_1,\ldots,W_4$ as follows. 
For all $k\in \four$, the subgraph $W_k$ has the same vertices as $G$ and set of arcs
$$W_k=\{a \textrm{ inner arc}\mid k\in[\linit(a),\rinit(a)[ \}\cup \{(v_{k+1},v_{k+2}),(v_{k+2},v_{k+3}),(v_{k},v_{k+3})\}.$$
The oriented subgraphs $W_1,\ldots,W_4$ corresponding to the 4-GS labeling in Figure~\ref{fig:4-grand-Schnyder}(a)
are represented in Figure~\ref{fig:4-grand-Schnyder}(b). We call $\cW=(W_1,\ldots,W_4)$ a \emph{4-GS wood} of $G$, by analogy to the woods defined by Schnyder for triangulations \cite{Schnyder:wood1}. 
It was shown in \cite{OB-EF-SL:Grand-Schnyder} that for all $k\in \four$, the subgraph $W_k$ is an oriented spanning tree rooted at $v_{k+3}$, with every arc in the tree oriented from child to parent. Moreover, in the rooted tree $W_k$, the vertices $v_k$ and $v_{k+1}$ are leaves (i.e. they are not incident to any inner edge of $G$ in $W_k$).\footnote{We mention that there is a small change of convention in the choice of the root of $W_k$ in the present article compared to \cite{OB-EF-SL:Grand-Schnyder}: in  \cite{OB-EF-SL:Grand-Schnyder} the root of $W_k$ was defined to be $v_k$.} 

\subsection{Bipolar orientations}\label{sec:4GS2}
In this subsection we explain how to use 4-GS structures to define several bipolar orientations which play a key role in the drawing algorithms. 
Recall that a \emph{bipolar orientation} is an orientation of a graph which is acyclic and has a unique \emph{source} (vertex only incident to outgoing arcs) and a unique \emph{sink} (vertex only incident to incoming arcs).
Given a 4-GS labeling $\cL$ of $G$, we define the \emph{odd} and \emph{even bipolar orientations} $B_o$ and $B_e$ as follows:
\begin{eqnarray*}
B_o= W_1\cup W_3^{-} \textrm{ and } B_e=W_4\cup W_2^{-}.
\end{eqnarray*}
In the above definition, the symbol $W_k^{-}$ indicates the oriented subgraph obtained from $W_k$ by reversing every arc. 
The bipolar orientations $B_o$ and $B_e$ are indicated in Figure~\ref{fig:bipolars-primal}(a). 
\fig{width=\linewidth}{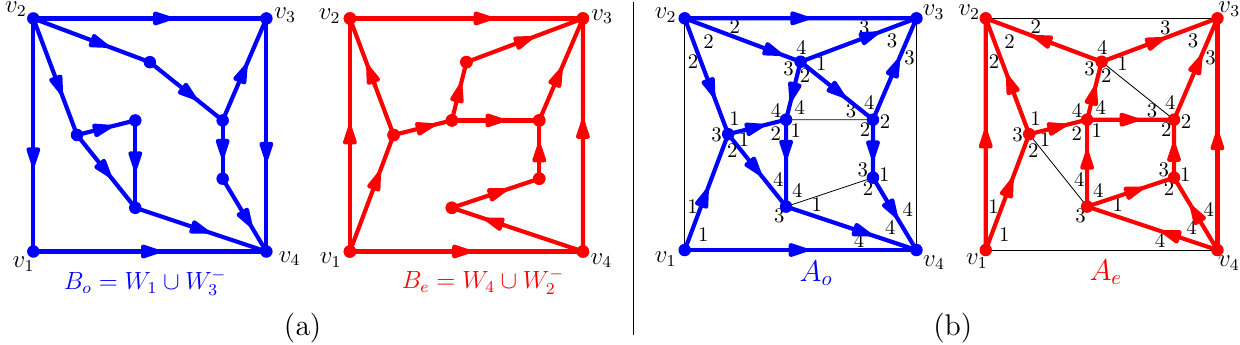}{(a) Odd and even bipolar orientations. (b) Odd and even acyclic orientations.}

It is easy to check that neither $B_o$ nor $B_e$ contains any edge oriented both ways, that is, if an arc $a$ belongs to $B_o$ (resp. $B_e$) then the opposite arc $-a$ does not belong to $B_o$ (resp. $B_e$). Hence $B_o$ and $B_e$ are oriented subgraphs of $G$. 
Furthermore we will prove below that $B_o$ and $B_e$ are not usurping their name, as these orientations are indeed bipolar.

Next we define the \emph{odd} and \emph{even acyclic orientations} $A_o$ and $A_e$.
These are the oriented subgraphs of $G$ having the same vertices as $G$ and the following sets of arcs:
\begin{eqnarray*}
A_o &\!\!\!=\!\!\!&\{a \textrm{ inner arc}\mid \linit(a)=1\textrm{ or }\rinit(a)=2\textrm{ or }\rterm(a)=3\textrm{ or }\lterm(a)=4\}\\
&&\cup\{(v_1,v_4),(v_2,v_3)\},\\
A_e &\!\!\!=\!\!\!&\{a \textrm{ inner arc}\mid \linit(a)=4\textrm{ or }\rinit(a)=1\textrm{ or }\rterm(a)=2\textrm{ or }\lterm(a)=3\}\\
&&\cup\{(v_1,v_2),(v_4,v_3)\}.
\end{eqnarray*}
This definition of $A_o$ and $A_e$ is illustrated in Figure~\ref{fig:labels-around-edge}(b).  An example is given in Figure~\ref{fig:bipolars-primal}(b).
It is easy to see that $A_o$ and $A_e$ do not contain any edge oriented both ways, hence they are oriented subgraphs of $G$. 
Observe also that $B_o$ is contained in $A_o\cup \{(v_2,v_1),(v_3,v_4)\}$, and $B_e$ is contained in $A_e\cup \{(v_1,v_4),(v_2,v_3)\}$. The situation is represented in Figure \ref{fig:label-to-bipolar}.

\begin{remark} \label{rk:dual-first-rk}
It is easy to see that an inner arc $a$ of $G$ is in $A_o$  if and only if $2\in [\rinit(a):\rterm(a)[$ or $4\in[\lterm(a):\linit(a)[$. This is to be compared to the definition of $B_e$: an inner arc $a$  is in $B_e$ if and only if  $2\in[\rterm(a):\lterm(a)[$ or $4\in [\linit(a):\rinit(a)[$.
Similarly, $a$ is in $A_e$ if $1\in [\rinit(a):\rterm(a)[$ or $3\in[\lterm(a):\linit(a)[$, while it is in $B_o$ if  $1\in [\linit(a):\rinit(a)[$ or $3\in[\rterm(a):\lterm(a)[$. We will comment further on this similarity in Section \ref{sec:dual4GS}.
\end{remark}

\fig{width=\linewidth}{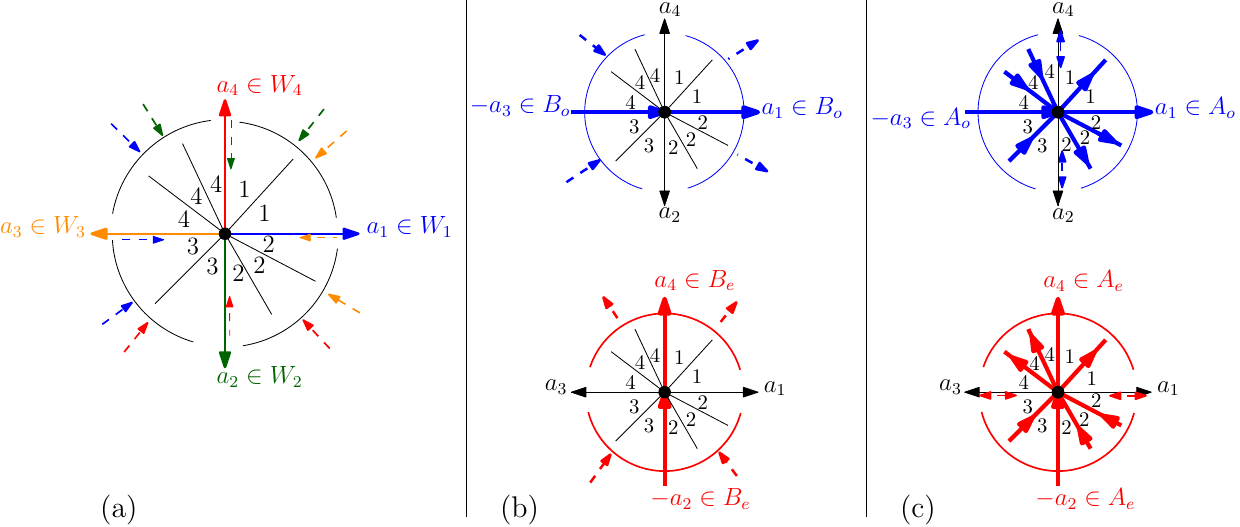}{(a) The arcs of $W_1,\ldots,W_4$ around a vertex $v$, superimposed with the corner labels around $v$. The arc $a_i$ is the outgoing arc of $W_i$ at $v$. The arcs $a_1,\ldots,a_4$ appear on this clockwise order around $v$, but we may have $a_i=a_{i+1}$. The possible incoming arcs belonging to the edges between $a_i$ and $a_{i+1}$ are indicated: the edges can bear incoming arcs from $W_{i+2}$ or $W_{i+3}$ (or both, or neither). Meanwhile the arc opposite to $a_i$ may belong to $W_{i+1}$ (b) The arcs of $B_o$ and $B_e$ around $v$. For instance the incoming arcs of $B_o$ at $v$ must appear strictly between $a_2$ and $a_4$ in clockwise order around $v$.  (c) The arcs of $A_o$ and $A_e$ around $v$.}

Let $\ovAo$ (resp. $\ovAe$) be the set of edges of $G$ which are not in the support of $A_o$ (resp. $A_e$). The set $\ovAo$ contains the outer edges $\{v_1,v_2\}$ and $\{v_3,v_4\}$ together with some inner edges, and it is easy to see that these inner edges are exactly those supporting the arcs $a$ satisfying 
$$\linit(a)=4~\textrm{ and }~\rinit(a)=1~\textrm{ and }~\rterm(a)=2~\textrm{ and }~\lterm(a)=3.$$
Furthermore it is easy to see that $\ovAo$ does not contain any cycle. Indeed, if there was such a cycle in $\ovAo$, then it would support a cycle of $W_4$, which is impossible since $W_4$ is a tree. Similarly $\ovAe$ has no cycle, hence both $\ovAo$ and $\ovAe$ are forests. One can therefore consider the oriented \plm $\wtA_o=A_o/\ovAo$ (resp. $\wtA_e=A_e/\ovAe$), that is, the oriented map obtained from $G$ by orienting the edges in $A_o$ (resp. $A_e$) and contracting the edges in $\ovAo$ (resp. $\ovAe$).

\begin{lemma}\label{lem:primal-bipolar}
Let $G$ be a 3,4-angulation of the square endowed with a 4-GS labeling. Let $B_o,B_e,A_o,A_e,\wtA_o,\wtA_e$ be the oriented \plms defined as above.
Then $\wtA_o$ and  $\wtA_e$ are bipolar orientations. Consequently, $A_o$ and  $A_e$ are acyclic orientations, and $B_o,B_e$ are bipolar orientations. 
\end{lemma}

\begin{proof}
We start by showing that $\wtA_o$ is acyclic. Suppose the contrary. For each simple directed cycle of $\wtA_o$ we consider the region of the plane it encloses. Let $C$ be a directed cycle of $\wtA_o$ whose enclosed region is minimal (not containing a smaller such region). By minimality of $C$, there is no directed path in $\wtA_o$ starting and ending on $C$ and staying strictly in the region enclosed by $C$. 
Since every inner vertex of $\wtA_o$ is incident to at least one incoming and at least one outgoing arc, the region enclosed by $C$ does not contain any arc (otherwise, one could extend this arc into a directed path in both direction, so as to construct a directed path starting and ending on $C$ and staying in the region enclosed by $C$). Hence the arcs of $C$ are all incident to a single inner face of $G$. However, it is easy to check that for any inner face $f$ of $G$, the set of arcs of $\wtA_o$ incident to $f$ does not form a directed cycle. Indeed any inner face $f$ of $G$ is incident to at least one arc of $A_o$ oriented clockwise around $f$ (since $f$ has a corner labeled 2 or a corner labeled 3, or both) and to at least one arc of $A_o$ oriented counterclockwise around $f$ (since $f$ has a corner labeled 1 or a corner labeled 4, or both).
This gives a contradiction and proves that $\wtA_o$ is acyclic. 
Furthermore, we claim the the unique source (resp. sink) of $\wtA_o$ is the vertex resulting from the contraction of the outer edge $\{v_1,v_2\}$ (resp. $\{v_3,v_4\}$). Indeed any inner vertex of $G$ is incident to at least one incoming arc and one outgoing arc (since the inner arcs of $B_o=W_1\cup W_3^{-}$ are in $A_o$), hence this also holds for any inner vertex of $\wtA_o$. Hence $\wtA_o$ is a bipolar orientation.

Since $\wtA_o$ is acyclic, the orientation $A_o$ is also acyclic (because if $A_o$ contained a directed cycle, then its contraction $\wtA_o$ would also contain one).
In fact $A_o\cup\{ (v_2,v_1),(v_3,v_4)\}$ is also acyclic (for the same reason). Hence the orientation $B_o$, which is contained in $A_o\cup\{ (v_2,v_1),(v_3,v_4)\}$, is also acyclic. Furthermore $B_o$ is bipolar since its unique sink (resp. source) is the outer vertex $v_2$ (resp. $v_4$).

Similarly, $\wtA_e$ is bipolar, $A_e$ is acyclic, and $B_e$ is bipolar (all the definitions coincide with the odd case up to increasing the indices by 1 modulo 4).
\end{proof}


\subsection{Dual Grand-Schnyder structures}\label{sec:dual4GS}
We now recall  from \cite{OB-EF-SL:Grand-Schnyder} some basic facts about \emph{dual 4-GS structures}. 

Let $G$ be a 3,4-angulation of the square. Let $G^*$ be the rooted 3,4-map dual to $G$, and let $\vinfty$ be its root vertex. The edges incident to $\vinfty$ are denoted by
$e_1^*,\ldots,e_4^*$, where $e_i^*$ is the edge of $G^*$ dual to the outer edge  $\{v_i,v_{i+1}\}$ of $G$. 
The edges $e_1^*,\ldots,e_4^*$ are called \emph{root edges}, and all the other edges are called \emph{non-root edges}. 
The faces $f_1^*,\ldots,f_4^*$ in counterclockwise order around $\vinfty$ are called \emph{root-faces}, where $f_i^*$ is dual to the outer vertex $v_i$ of $G$. All the other faces are called \emph{non-root faces}.



\begin{definition}\label{def:dual_label}
Let $G^*$ be a rooted 3,4-map. A \emph{dual 4-GS labeling} of $G^*$ is a labeling of its non-root corners in $\four$ satisfying the following conditions.
\begin{itemize}
    \item[(L0\sups{*})] For all $i \in \{1,2,3,4\}$, the corners incident to the root-face $f_i^*$ have label $i$.
    \item[(L1\sups{*})] Around every non-root face or non-root vertex, the sum of clockwise jumps is 4.
    \item[(L2\sups{*})] Consecutive corners around a vertex have distinct labels.
    \item[(L3\sups{*})] Let $a$ be an non-root arc whose initial vertex has degree 3. Let $c,c',c''$ be the left-initial, right initial and right terminal corners of $e$. The label jump $\delta$ from $c$ to $c'$ and the label jump $\delta$ from $c'$ to $c''$ satisfy $\delta+\eps\geq 2$. 
\end{itemize} 
\end{definition}

\fig{width=\linewidth}{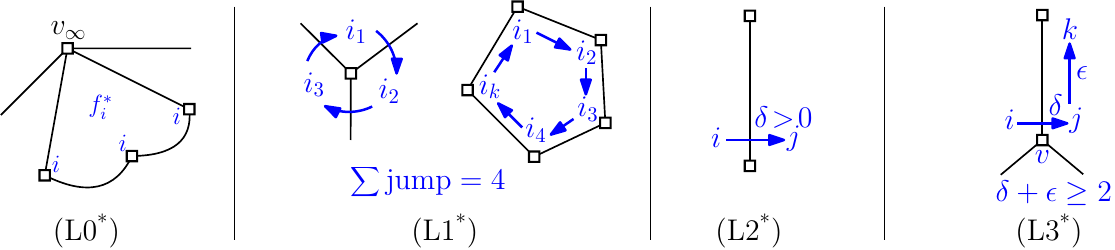}{Conditions defining dual 4-GS labelings.}

The conditions defining dual 4-GS labelings are illustrated in Figure \ref{fig:dual_label}. Observe that if $G$ is a 3,4-angulation, then  a 4-GS labeling of $G$ induces a  dual 4-GS labeling of $G^*$ as represented in Figure \ref{fig:4_dual_big}(a), since Conditions (L0\sups{*})--(L3\sups{*}) are dual to Conditions (L0)--(L3).

\fig{width=\linewidth}{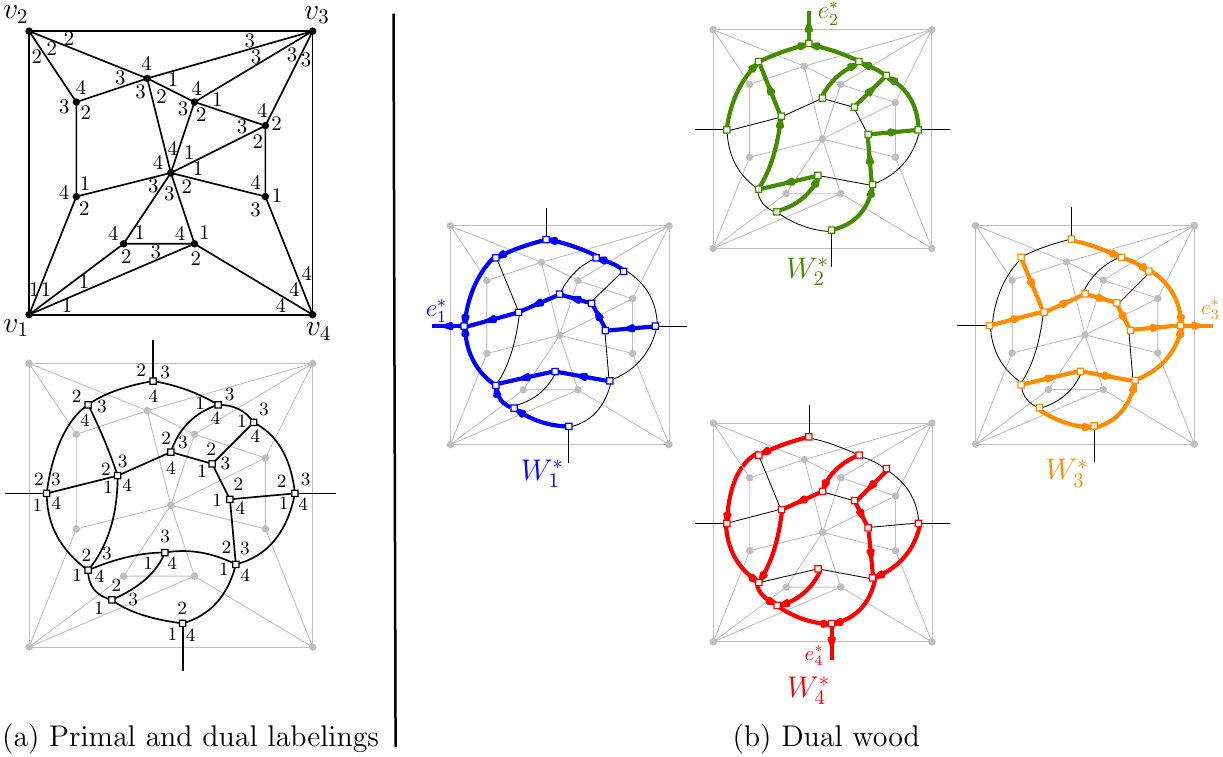}{(a) A 4-GS labeling $\cL$ (top) and the dual labeling $\cL^*$ (bottom). (b) The corresponding dual 4-GS wood $\mathcal{W^*} = (W_1^*,W_2^*,W_3^*,W_4^*)$.}

Given a dual 4-GS labeling $\cL^*$ one can define some subsets of arcs $W_1^*,...,W_4^*$ using a similar rule as in Section \ref{sec:4GS1}: for all $k\in \four$, 
$$W_k^*=\{a\textrm{ arc of }G^*\mid k\in[\linit(a),\rinit(a)[\}.$$ 
An example is represented in Figure \ref{fig:4_dual_big}(b).
It was shown in \cite{OB-EF-SL:Grand-Schnyder} that for all $k\in\four$, the oriented subgraph $W_k^*$ is a spanning tree oriented toward the root vertex $\vinfty$ and containing  the root-edge $e_k^*$. The tuple $\mathcal{W}^* = (W_1^*,...,W_4^*)$ is called a \emph{dual 4-GS wood}.

Given a dual 4-GS wood $\mathcal{W}^* = (W_1^*,...,W_4^*)$, we define some bipolar orientations and acyclic orientations. The \emph{odd} and \emph{even dual bipolar orientations} are defined as 
$$B_o^* = W_3^* \cup (W_1^*)^{-}~\textrm{ and }~B_e^* = W_2^* \cup (W_4^*)^{-}.$$ 
An example is represented in Figure \ref{fig:bipolars_dual_big}(a). 

\fig{width=\linewidth}{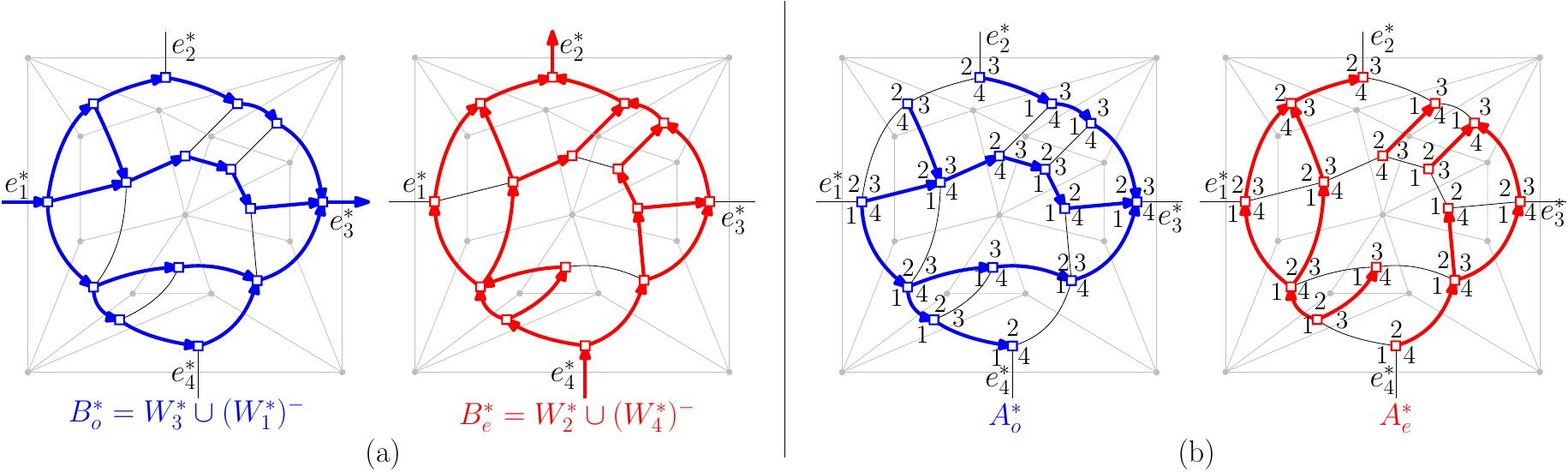}{(a) Odd and even dual bipolar orientations. (b) Odd and even dual acyclic orientations, which will be defined in Section \ref{sec:increasing_dual}.}

As we now explain, the orientations  $B_o^*$ and $B_e^*$ are actually the dual of the bipolar contractions $\wtA_e$ and $\wtA_o$ defined in Section \ref{sec:4GS1}. This is represented in Figure \ref{fig:dual-orientations}. For an arc $a$ of a \plm $G$, the \emph{dual arc} $a^*$ is the arc of $G^*$ crossing $a$ from left to right. The \emph{dual} of an oriented \plm $A$ is the oriented \plm consisting of all the dual arcs.
Let $\mathcal{L}$ be a 4-GS labeling of a 3,4-angulation of the square $G$, and let $\mathcal{L}^*$ be the corresponding dual $4$-GS labeling of $G^*$. It is easy to see, considering Remark \ref{rk:dual-first-rk}, that the arcs in $B_o^*$ (resp. $B_e^*$) are exactly the dual of the bipolar contractions $\wtA_e$ (resp. $\wtA_o^-$). 
Hence  $B_o^*$ (resp. $B_e^*$) is the dual of the oriented map $\wtA_e$ (resp. $\wtA_o^-$). The fact that $\wtA_o$ and $\wtA_e$ are bipolar orientations implies that $B_o^*$, $B_e^*$ (or rather the orientation obtained from $B_o^*$, $B_e^*$ by disconnecting the root-edges at the root vertex) are also bipolar orientations, as a consequence of duality for bipolar orientations~\cite{de1995bipolar}. 

\fig{width=\linewidth}{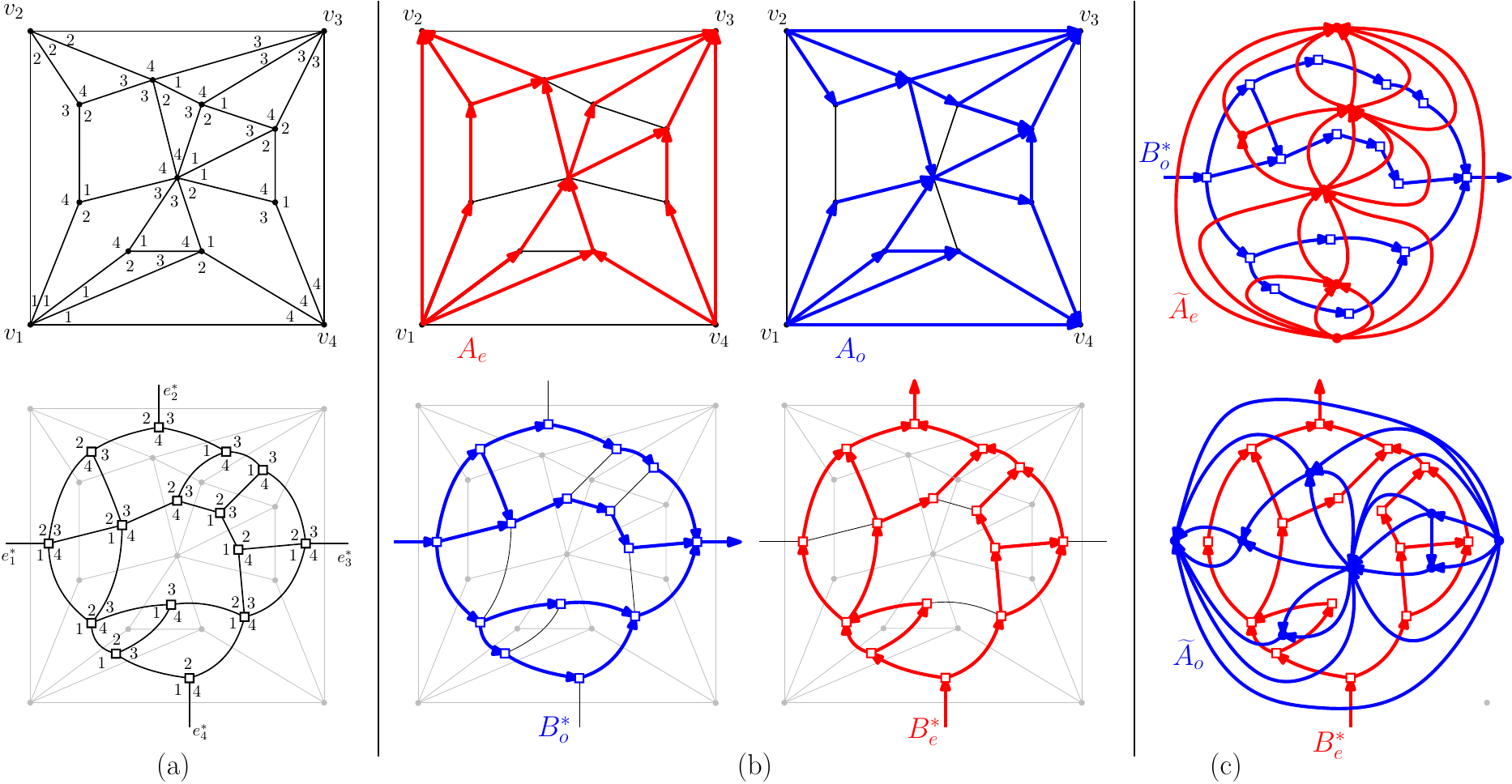}{(a) A 4-GS labeling $\mathcal{L}$ (top) and its dual labeling $\cL^*$ (bottom). (b) The acyclic orientations $A_e$ and $A_o$ (top) and the dual bipolar orientations $B_o^*$ and $B_e^*$. (c) Duality relation between $\wtA_e$ and $B_o^*$ (top) and between $\wtA_o^-$ and $B_e^*$ (bottom).}


\section{Straight-line drawing of \plms with faces-degree at most 4}\label{sec:primal_algo}

In this section we define and analyze straight-line drawing algorithms for adapted 3,4-angulations of the square. We first define a ``face-counting'' algorithm generalizing the algorithms from~\cite{Fu07b,Barriere-Huemer:4-Labelings-quadrangulation}, and then a more general (hence more compact) ``increasing functions'' algorithm.

\subsection{Face-counting drawing algorithm for 3,4-angulations}\label{sec:face-counting-primal}
In this subsection we define a straight-line drawing algorithm for 3,4-angulations of the square based on 4-GS structures.

Let $G$ be an adapted 3,4-angulation of the square. Let $\cL$ be a 4-GS labeling of $G$, and let $W_1$, $W_2$, $W_3$, $W_4$ be the corresponding 4-GS trees (which are rooted at $v_2$, $v_3$, $v_4$, $v_1$ respectively). 
Recall from Section~\ref{sec:4GS2} that the \emph{odd and even bipolar orientations} associated to $\cL$ are defined as $B_o=W_1\cup W_3^{-}$ (which has source $v_2$ and sink $v_4$) and $B_e=W_4\cup W_2^{-}$ (which has source $v_1$ and sink $v_3$). 
For each vertex $v$ of $G$, and for $i\in \four$ we define the path $P_i(v)$ as the directed path from $v$ to the root $v_{i-1}$ in the rooted tree $W_i$. 
We also define the \emph{odd and even paths} for $v$ as 
$$P_o(v)=P_1(v)\cup P_3(v)^-\subset B_o\textrm{ and }P_e(v)=P_4(v)\cup P_2(v)^-\subset B_e,$$
where $P_i(v)^-$ is the directed path of $G$ obtained by reversing $P_i(v)$. Note that $P_o(v)$ is a directed path from $v_2$ to $v_4$ in the bipolar orientations $B_o$, hence this path is simple (does not contain any cycle). Similarly, $P_e(v)\subset B_e$ is a simple directed path from $v_1$ to $v_3$.

We define the \emph{face-counting coordinates} of $v$ as $(\ell_e(v),r_o(v))$, where $r_o(v)$ is the number of inner faces of $B_o$ on the right of the directed path $P_o(v)$, and $\ell_e(v)$ is the number of inner faces of $B_e$ on the left of the directed path $P_e(v)$. The face-counting coordinates are represented in Figure~\ref{fig:face-counting}(a).

\begin{thm}\label{thm:face-counting-primal}
Let $G$ be an adapted 3,4-angulation of the square and let $\cL$ be a 4-GS labeling of~$G$. 
Placing each vertex $v$ of $G$ at its face-counting coordinates $(\ell_e(v),r_o(v))\in \RR^2$, and drawing each edge of $G$ as a segment, yields a planar straight-line drawing of $G$.  
\end{thm}

\fig{width=\linewidth}{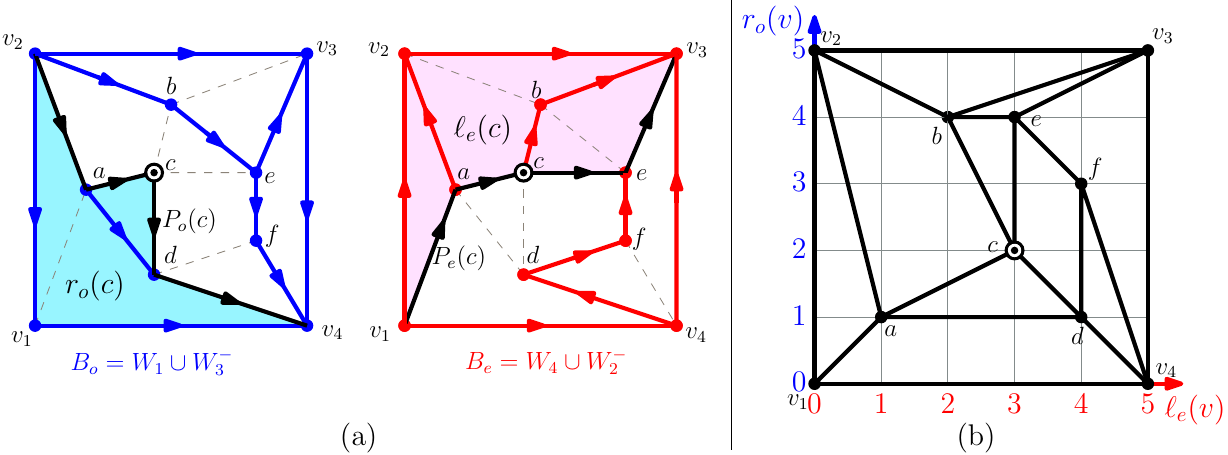}{(a) The paths and regions used to compute the face-counting coordinates of a vertex. For the vertex $c$ we get $\ell_e(c)=3$ and $r_o(c)=2$. (b) The drawing of $G$ obtained by the face-counting algorithm.}

We will prove Theorem~\ref{thm:face-counting-primal} in Section~\ref{sec:proof-planar-primal} as a corollary of a more general result.

We now compare the face counting algorithm of Theorem~\ref{thm:face-counting-primal} to the algorithms from~\cite{Fu07b} and~\cite{Barriere-Huemer:4-Labelings-quadrangulation} described in the introduction. We start with the algorithm from~\cite{Fu07b}. Let $G$ be an irreducible triangulation of the square. 
As explained in Section~\ref{sec:4GS1}, there is a natural bijection between the transversal structures of $G$ and the 4-GS labelings of $G$, which is given by the labeling rule of Figure~\ref{fig:labeling-separating}(left). 
Moreover it is clear from the definitions that for any inner vertex $v$, the odd path $P_o(v)$ of $\cL$ is equal to the blue separating path $\Pb(v)$ of $\cT$, and the even path $P_e(v)$ of $\cL$ is equal to the red separating path $\Pr(v)$ of $\cT$. Note also that the odd bipolar orientation $B_o=W_1\cup W_3^-$ of $\cL$ is contained in the blue map $G_b$ of $\cT$, and the even bipolar orientation $B_e=W_4\cup W_2^-$ is contained in the red map $G_r$. Precisely, $B_o$ is made of the edges of the blue map $G_b$ which are on a blue separating path (for some inner vertex), while $B_e$ is made of the edges of the red map $G_r$ which are on a red separating path. Therefore, when applying the face-counting algorithm of Section~\ref{sec:face-counting-primal} to an irreducible triangulation, one gets the $x$-coordinate (resp. $y$-coordinate) of each inner vertex $v$ by counting the faces at the left of the red separating path $\Pr(v)$ (resp. right of the blue separating path $\Pb(v)$) in the submap $B_e$ of $G_r$ (resp. submap $B_o$ of $G_b$). Hence our algorithm gives an optimization of the algorithm based on transversal structures presented in the introduction. This optimization is precisely the one given in~\cite[Section 5.2]{Fu07b} (compact version of the algorithm where every grid-line is occupied by at least one vertex). 

Next we discuss the face-counting algorithm from~\cite{Barriere-Huemer:4-Labelings-quadrangulation} based on separating decompositions. Let $G$ be a simple quadrangulation. As explained in Section~\ref{sec:4GS1}, there is a natural bijection between the separating decompositions of $G$ and the even 4-GS labelings of $G$, which is given by the labeling rule of Figure~\ref{fig:labeling-separating}(right). Let $\cS$ be a separating decomposition, and let $\cL$ be the associated 4-GS labeling.
It is clear from the definitions that the blue and orange trees $T_b$ and $T_o$ of $\cS$ are equal to the trees $W_1$ and $W_3$ of $\cL$. Hence every edge of $G$ belongs to the odd bipolar orientation $B_o=W_1\cup W_3^-$. Similarly, every edge of $G$ belongs to the even bipolar orientation $B_e=W_4\cup W_2^-$. 
Note also that, for any inner vertex $v$ and any $i\in\four$, the path $P_i(v)$ defined from the separating decomposition $\cS$ is equal to the path from $v$ to $v_i$ in the tree $W_i$. Therefore, the separating paths $P_o(v)$ and $P_e(v)$ defined from $\cS$ are equal to the odd and even paths defined from $\cL$, and the placement of vertices by the face-counting algorithm in~\cite{Barriere-Huemer:4-Labelings-quadrangulation} based on $\cS$ coincides with the placement from the face-counting algorithm based on~$\cL$. 

In conclusion, the straight-line drawing algorithm from Theorem~\ref{thm:face-counting-primal} interpolates between the compact version of the algorithm from~\cite{Fu07b}, and the face-counting algorithm from~\cite{Barriere-Huemer:4-Labelings-quadrangulation}.


\begin{remark} 
It is shown in \cite{Barriere-Huemer:4-Labelings-quadrangulation} that the face-counting algorithm can be supplemented by a compaction procedure that significantly reduces the grid-size. It is an open question to determine how to best extend this compaction procedure to the more general setting of the present article. 
\end{remark}

\subsection{Increasing function drawing algorithm for 3,4-angulations}\label{sec:increasing-primal}
In this section we present a drawing algorithm for adapted 3,4-angulations of the square which generalizes the one presented in the previous subsection.

Let $G$ be an adapted 3,4-angulation of the square, and let $\cL$ be a 4-GS labeling of~$G$. Recall from Section~\ref{sec:4GS2} the definition of the \emph{odd and even acyclic orientations} $A_o$ and $A_e$ of~$G$. By Lemma~\ref{lem:primal-bipolar} the oriented \plms $\wtA_o=A_o/\ovAo$ and $\wtA_e=A_e/\ovAe$ are bipolar orientations. We call them the \emph{odd and even bipolar contractions} of $G$. Recall that the set of edges $\ovAo$ and $\ovAe$ are forests of $G$ (i.e. acyclic subgraphs). Hence the faces of $\wtA_o$ (resp. $\wtA_e$) are in 1-to-1 correspondence with 
the faces of $G$. We will now triangulate some of the faces of $\wtA_o$ and $\wtA_e$ of degree 4, in order to construct the ``augmented'' bipolar orientations $\aA_o$ and $\aA_e$. This is represented in Figure~\ref{fig:diagonal-edges}.

\fig{width=\linewidth}{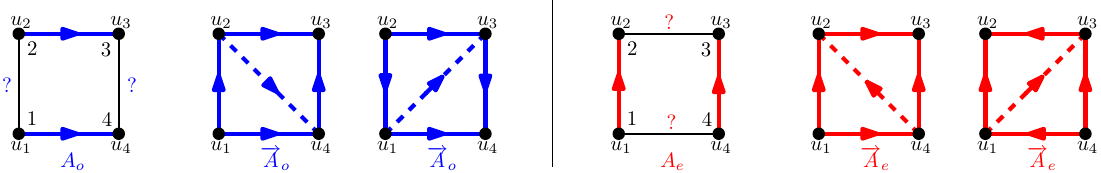}{Adding diagonal arcs to some of the quadrangular faces of $A_o$ (left) and $A_e$ (right) in order to construct the augmented bipolar contractions $\protect\aA_o$ and $\protect\aA_e$.}


Let $f$ be a face of $G$ of degree 4, and let $u_1,u_2,u_3,u_4$ be the vertices of $G$ incident to the corners of $f$ labeled $1,2,3,4$. By definition of $A_o$ (see, for example, Figure~\ref{fig:labels-around-edge}(b)), the two arcs $(u_2,u_3)$ and $(u_1,u_4)$ are in $A_o$; The situation is illustrated in  Figure~\ref{fig:diagonal-edges}(left).
\begin{compactitem}
\item If the arcs $(u_1,u_2)$ and $(u_4,u_3)$ are in $A_o$, then we call $(u_2,u_4)$ an \emph{odd diagonal arc} of~$G$.
\item If the arcs $(u_2,u_1)$ and $(u_3,u_4)$ are in $A_o$, then we call $(u_1,u_3)$ an \emph{odd diagonal arc} of~$G$.
\end{compactitem}
The \emph{augmented odd bipolar contraction} is the oriented map $\aA_o$ obtained from $A_o$ by adding all the odd diagonal arcs and then contracting the edges in $\ovAo$. Note that $\aA_o$ is obtained from $\wtA_o$ by adding diagonal arcs in some of its faces of degree 4 as represented in Figure~\ref{fig:diagonal-edges}(left).


\fig{width=\linewidth}{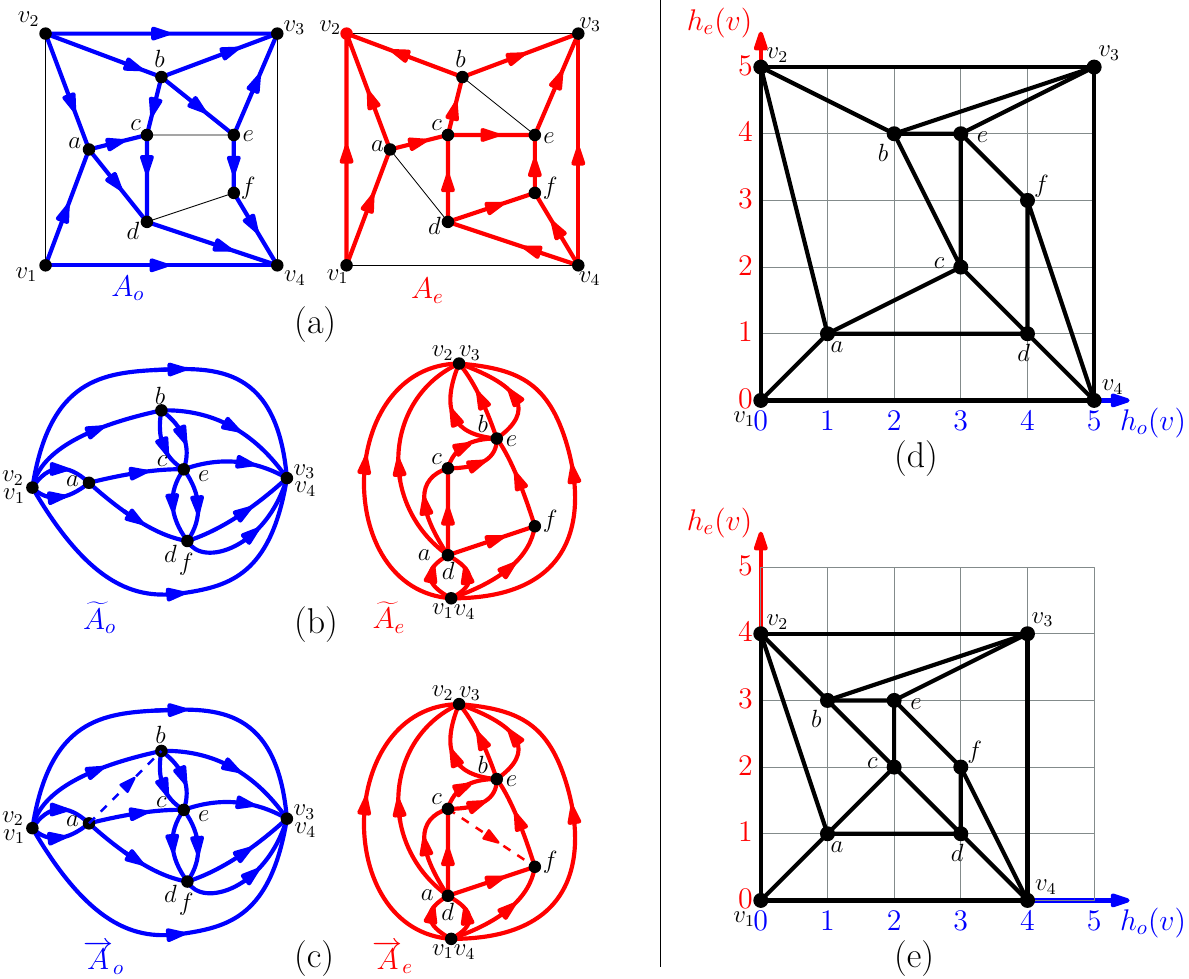}{(a) The acyclic orientations. (b) The bipolar contractions. (c) The augmented bipolar contractions. (d) The drawing obtained by choosing some increasing function $h_o$ and $h_e$. (e) The drawing obtained by choosing quasi-increasing functions (as defined in Section~\ref{sec:optimization}) $h_o$ and $h_e$.}

Similarly, 
the two arcs $(u_1,u_2)$ and $(u_4,u_3)$ around the face $f$ are in $A_e$; see Figure~\ref{fig:diagonal-edges}(right).
\begin{compactitem}
\item If the arcs $(u_1,u_4)$ and $(u_2,u_3)$ are in $A_e$, then we call $(u_4,u_2)$ an \emph{even diagonal arc} of~$G$.
\item If the arcs $(u_4,u_1)$ and $(u_3,u_2)$ are in $A_e$, then we call $(u_1,u_3)$ an \emph{even diagonal arc} of~$G$.
\end{compactitem}
The \emph{augmented even bipolar contraction} is the oriented map $\aA_e$ obtained from $A_e$ by adding all the even diagonal arcs and contracting the edges in $\ovAe$. Note that $\aA_e$ is obtained from $\wtA_e$ by adding diagonal arcs in some of its faces of degree 4 as represented in Figure~\ref{fig:diagonal-edges}(right). 
As explained in Remark~\ref{rk:vecAobipolar} below, the oriented maps $\aA_o$ and $\aA_e$ are indeed bipolar orientations. 


Let $V$ be the vertex set of $G$. 
An \emph{$\aA_o$-increasing} function for $G$ is a function $h:V\to \RR$ that satisfies the following 2 conditions: \begin{compactitem}
\item $h(u)=h(v)$ for every edge $\{u,v\}$ in $\ovAo$,
\item $h(u)<h(v)$ for every arc $(u,v)$ in $A_o$ and every odd diagonal arc $(u,v)$
\end{compactitem}
Such functions can be identified with functions defined on the vertices of $\aA_o$ which are strictly increasing along each arc of $\aA_o$. Similarly, an \emph{$\aA_e$-increasing} function for $G$ is a function $h:V\to \RR$ such that $h(u)=h(v)$ for every edge $\{u,v\}$ in $\ovAe$, and $h(u)<h(v)$ for every arc $(u,v)$ in $A_e$ and every even diagonal arc $(u,v)$.

\begin{thm}\label{thm:increasing-primal}
Let $G$ be an adapted 3,4-angulation of the square, let $\cL$ be a 4-GS labeling, and let $\aA_o$ and $\aA_e$ be the corresponding augmented bipolar contractions of $G$. 
Let $h_o$ be an $\aA_o$-increasing function, and let $h_e$ be an $\aA_e$-increasing function.
Placing each vertex $v$ of $G$ at coordinates $(h_o(v),h_e(v))\in \RR^2$, and drawing each edge of $G$ as a segment, yields a planar straight-line drawing of $G$. 
\end{thm}

\begin{example} Theorem~\ref{thm:increasing-primal} is illustrated in Figure~\ref{fig:increasing-functions}(d). For the example in that figure, an $\aA_o$-increasing function is any function $h_o$ defined on $V=\{v_1,v_2,v_3,v_4,a,b,c,d,e,f\}$ such that 
$$h_o(v_1)=h_o(v_2)<h_o(a)<h_o(b)<h_o(c)=h_o(e)<h_o(d)=h_o(f)<h_o(v_3)=h_o(v_4),$$
and an $\aA_e$-increasing function is any function $h_e:V\to \RR$ such that
 $$h_e(v_1)=h_e(v_4)<h_e(a)=h_e(d)<h_e(c)<h_e(f)<h_e(b)=h_e(e)<h_e(v_2)=h_e(v_3).$$
Such functions $h_o,h_e$ are used in Figure~\ref{fig:increasing-functions}(d).
 \end{example}

The proof of Theorem~\ref{thm:increasing-primal} will be given in Section~\ref{sec:proof-planar-primal}. 
We will also give a slightly stronger result in Section~\ref{sec:optimization}, which allows for more compact grid drawings such as the one shown in Figure~\ref{fig:increasing-functions}(e).

As we now explain, Theorem~\ref{thm:increasing-primal} implies Theorem~\ref{thm:face-counting-primal} about the face-counting algorithm. Indeed we will prove the following result. 
\begin{lemma}\label{lem:face-cout-increasing}
Let $G$ be a 3,4-angulation of the square endowed with a 4-GS labeling.
Let $\aA_o$ and $\aA_e$ be the augmented bipolar contractions of $G$. The face-counting coordinates $(\ell_e(v),r_o(v))$ defined in Section~\ref{sec:face-counting-primal} are such that the function $\ell_e$ is $\aA_o$-increasing and the function $r_o$ is $\aA_e$-increasing.
\end{lemma}

\begin{remark}\label{rk:vecAobipolar}
Lemma~\ref{lem:face-cout-increasing} shows that there exist some $\aA_o$ and $\aA_e$-increasing functions. This shows that $\aA_o$ and $\aA_e$ are acyclic orientations. It is also clear that these orientations have a unique source and unique sink (since this is the case for $\wtA_o$ and $\wtA_e$), hence $\aA_o$ and $\aA_e$ are bipolar orientations. This could also have been shown by the same arguments as in the proof of Lemma~\ref{lem:primal-bipolar}, upon observing that the contours of the new faces of degree $3$ (created by the insertion of diagonals) are acyclic.
\end{remark}

\begin{proof}[Proof of Lemma~\ref{lem:face-cout-increasing}]
Let us prove that the face-counting function $\ell_e$ is $\aA_o$-increasing. 
 There are three facts to prove:
\begin{compactenum}[(a)]
\item if $\{u,v\}$ is an edge in $\ovAo$, then $\ell_e(u)=\ell_e(v)$,
\item if $(u,v)$ is an arc in $A_o$, then $\ell_e(u)<\ell_e(v)$,
\item if $(u,v)$ is a diagonal arc of $\aA_o$, then $\ell_e(u)<\ell_e(v)$.
\end{compactenum}
Recall that $\ell_e(v)$ is defined as the number of faces of the bipolar orientation $B_e$ on the left of the path $P_e(v)=P_4(v)\cup P_2(v)^-$.

We start by proving (a). The statement is clear for outer edges, so we now consider an inner edge $e=\{u,v\}\in \ovAo$. Recall that the corners incident to $e$ are labeled as indicated in Figure~\ref{fig:labels-around-edge}(b), so we have $(u,v)\in W_4$ and $(v,u)\in W_2$. In this case the paths $P_e(u)$ and $P_e(v)$ are equal, hence $\ell_e(u)=\ell_e(v)$.

Next, we prove (b). The statement is clear for outer arcs, so we now consider an inner arc $a=(u,v)$ in $A_o$. The condition satisfied by the labels of the corners incident to the arc $a$ is indicated in Figure~\ref{fig:labels-around-edge}(b). Let us first suppose that $\linit(a)=1$. In this case $\lterm(a)\in\{3,4\}$. We contend that the path $P_e(u)$ stays (weakly) to the left of $P_e(v)$. The situation is represented in Figure~\ref{fig:face-count-increasing}(a). Consider the path $P=P_2(u)^-\cup \{a\} \cup P_4(v)$. The label of the corners around $a$ implies that the first arc of $P_4(u)$ is strictly to the left of $P$ as indicated in Figure~\ref{fig:face-count-increasing}(a). Moreover, $P_4(u)$ cannot cross $P$. Indeed, $P_4(u)$ and $P_2(u)$ have no vertex in common except $u$ (otherwise there would be a directed cycle in the bipolar orientation $B_e$); and if $P_4(u)$ and $P_4(v)$ have a common vertex then they merge after that vertex. Hence $P_4(u)$ stays on the left of $P$. Thus $P_e(u)$ stays to the left of $P$.
Similarly, $P_2(v)$ starts weakly to the right of $P$, hence it stays weakly to the right of $P$. Thus $P_e(v)$ stays to the right of $P$. In conclusion $P_e(u)$ is to the left of $P$ which is to the left of $P_e(v)$. This gives the desired inequality $\ell_e(u)\leq \ell_e(v)$ for the number of faces on the left of $P_e(u)$ and $P_e(v)$. The inequality is actually strict since $P_e(u)\neq P$. This completes the proof for the case $\linit(a)=1$.

A similar proof can be given in the cases $\rinit(a)=2$, or $\lterm(a)=3$, or $\rterm(a)=4$ which are represented in Figure~\ref{fig:face-count-increasing}(b). Indeed in all these cases the label condition implies that $a$ is on the right of $P_e(u)$ and on the left of $P_e(v)$, which is enough to prove that $\ell_e(u)< \ell_e(v)$ by reasoning as above.

Lastly we prove (c). There are two types of diagonal arcs, as represented in Figure~\ref{fig:face-count-increasing}(c). In both cases, the diagonal $a=(u,v)$ is to the right of $P_e(u)$ and to the left of $P_e(v)$. Hence, as before we can see that the path $P=P_2(u)^-\cup \{a\} \cup P_4(v)$ is at the right of $P_e(u)$ and at the left of $P_e(v)$, which implies $\ell_e(u)< \ell_e(v)$ as desired. This completes the proof that the function $\ell_e$ is $\aA_o$-increasing. 

The proof that the function $r_o$ is $\aA_e$-increasing is almost identical and left to the reader. 
\end{proof}

\fig{width=.9\linewidth}{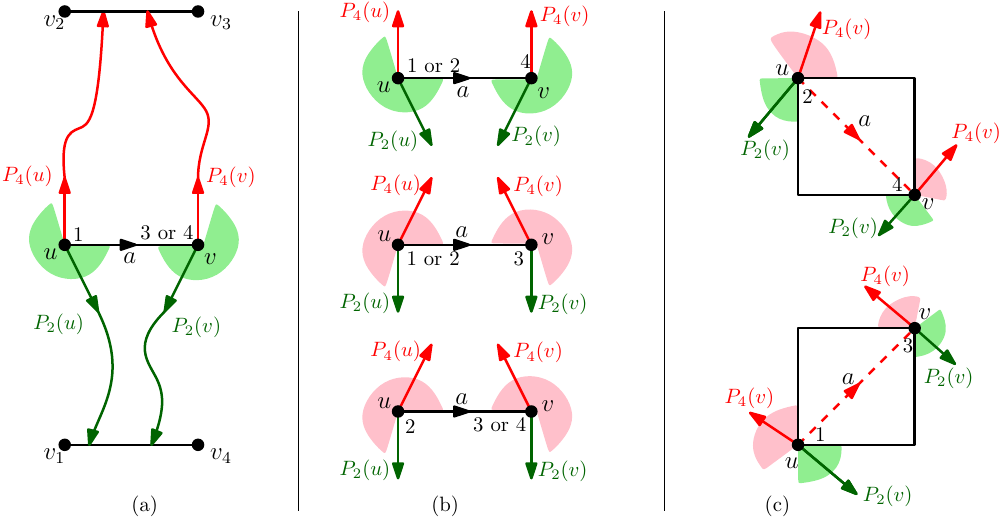}{Proof of Lemma~\ref{lem:face-cout-increasing}. (a) An arc $a\!=\!(u,v)$ such that $\linit(a)\!=\!1$. The clockwise order of $a$, $P_2(u)$ and $P_4(u)$ around $u$ is indicated, with the range for position of the first arc of $P_2(u)$ indicated by a green sector. Similarly the clockwise order of $a$, $P_2(v)$ and $P_4(v)$ around $v$ is indicated. (b) Local clockwise order of the relevant paths around an arc $a=(u,v)$ satisfying either $\lterm(a)=4$ (top) $\rinit(a)$ (middle) or $\rterm(a)=3$ (bottom). (c) Local clockwise order of the relevant paths around a diagonal arc $a=(u,v)$.} 


\subsection{Proof of planarity}\label{sec:proof-planar-primal}
In this subsection we prove Theorem~\ref{thm:increasing-primal}. We will rely on a general planarity criterion that we now state. It is well-known and appears in different forms in the literature, for instance as Lemma 4.4 in~\cite{CDV:these}, with references. We provide a short proof for completeness.

Let $G$ be a 2-connected \plm (for instance, an adapted 3,4-angulation of the square). Note that the contour of every face of $G$ is a simple cycle. 
A \emph{piecewise-linear drawing} of a graph $G$ is a drawing of $G$ in the plane where every edge is represented as a (non-empty) self-avoiding concatenation of segments.
The drawing is called \emph{valid} for a face $f$ of $G$ if the curve corresponding to the clockwise contour\footnote{\emph{Clockwise} means that the contour is traversed with $f$ on its right if $f$ is an inner face, and with $f$ on its left if $f$ is the outer face.} of $f$ is a polygonal simple curve in the plane whose interior is on its right. 

\begin{lemma}\label{lem:planar-criteria}
Let $G$ be a 2-connected \plmm. 
Consider a piecewise-linear drawing of the graph underlying $G$.
If the drawing is valid for every face of $G$ (including the outer face), then the drawing is planar (and corresponds to the \plm $G$). 
\end{lemma}

\begin{proof}
Consider a piecewise-linear drawing of a 2-connected \plm $G$, which is valid for every face. 
Let $F$ be the set of faces of $G$. 
For a face $f\in F$, let $P_f$ be the polygon representing $f$ in the drawing, and let $\Sigma_f$ be the set of points corresponding to the vertices incident to $f$. If $f$ is an inner face  (resp. is the outer face), we define $N_f:\mathbb{R}^2\backslash \Sigma_f\to\{0,\frac{1}{2},1\}$ to be the indicator function for $f$, assigning value $1$ to the points in the interior (resp. exterior) of $P_f$, value $1/2$ to the points on the boundary of $P_f$ (excluding $\Sigma_f$),  and value $0$ to the points in the exterior (resp. interior) of $P_f$. Let $\Sigma$ be the set of points in the plane that are images of vertices of $G$, and let $N$ be the function on $\mathbb{R}^2\backslash\Sigma$ defined as $N:=\sum_{f\in F} N_f$. 
As the drawing is valid for every face, the function $N$ is locally constant, hence it is constant. Since it is equal to $1$ for a point in the drawing of the outer face far enough so that it 
avoids the polygons representing inner faces, it is constant equal to $1$. 
Moreover, at each vertex $v$, the fact that the drawing is valid for all faces ensures that the set of faces in clockwise order around $v$ forms a fan that wraps $k$ times (i.e., of total angle $2k\pi$), for a positive integer $k$. Since $N=1$ in the neighborhood of~$v$, 
we must have $k=1$. Similarly, no two vertices can be at the same point, nor a vertex intersect an edge, nor two edges intersect each other (any of these situations would imply that $N\geq 2$ in a neighborhood). This proves the planarity of the drawing.
\end{proof}

In Section~\ref{sec:dual_algo} we will also need the following complementary lemma:

\begin{lemma}\label{lem:outer_right_visible}
Let $G$ be a 2-connected \plmm, endowed with a piecewise-linear drawing. Let $C=(C(t))_{t\in[0,1]}$ be the directed closed curve given by the drawing of the outer face contour of $G$ traversed clockwise.   
A point $p=C(t)$ is called \emph{outer-right-visible} if there is a ray starting from $p$ to the right side of $C$ (locally at $p$) and reaching to infinity without meeting $C$ again.
If the drawing is valid for every \emph{inner} face of~$G$, then there is no outer-right-visible point.   
\end{lemma}
\begin{proof}
Assume the drawing is valid at every inner face, and consider a ray $\gamma$ starting on the right side of $C$ at a point $p=C(t)$ and reaching to infinity.  
Let $p'$ be the farthest point on $\gamma$ that belongs to the drawing of an edge or vertex. 
Clearly, after leaving $p'$ the ray $\gamma$ cannot visit any point on an edge or vertex, nor in an area corresponding to the  drawing of an inner face. 
Hence, $p'$ cannot be the image of a point on an inner edge, nor the image of an inner vertex $v$ (for this we need the validity for inner faces, which ensures that the inner faces incident to $v$ form a fan that 
wraps $k$ times around $v$ in the drawing, for some positive integer $k$). Hence $p'\in C$. 
It only remains to show that $p'\neq p$.  
If $p$ belongs to the drawing of an edge of $C$, then this edge is on the contour of an inner face $f$ whose drawing is valid, hence $\gamma$ starts within the interior of $f$ and will meet the contour of $f$ again, thus $p'\neq p$. Suppose now that  $p$ corresponds to the drawing of a vertex $v$ of $C$. In this case, the drawing of the inner faces incident to $v$ (which are all valid) form a fan that covers the right side of $C$ at $p$ (and possibly wraps around $p$ several more times). Hence $\gamma$ starts in the interior of an inner face, or in the drawing of an inner edge. In every case we conclude that $p'\neq p$ which completes the proof.  
\end{proof}

\begin{remark}
The condition that the drawing is valid at every inner face is not enough to avoid self-intersections of the outer contour, as illustrated in Figure~\ref{fig:PL_drawing}, where it can be seen  that the drawing has no outer-right-visible point.
\end{remark}

\fig{width=3cm}{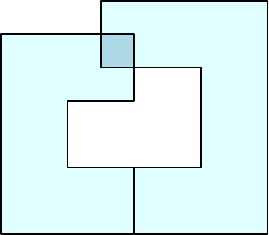}{A piecewise-linear drawing of a \plm which is valid for its two inner faces but not planar: its outer contour intersects itself.}

We now use Lemma~\ref{lem:planar-criteria} to prove Theorem~\ref{thm:increasing-primal}. Let $G$ be an adapted 3,4-angulation of the square, and let $V$ be its vertex set. We call a function $\al:V\to \RR^2$ a \emph{vertex-placement} for $G$. A vertex placement $\al$ induces a \emph{straight-line drawing} of $G$ (in which each edge is represented by a segment). We say that $\al$ is \emph{valid} for a face $f$, if the straight-line drawing of $f$ obtained from this vertex placement is valid. 
Let $\cL$ be a 4-GS labeling of $G$, let $h_o:V\to \RR$ be an $\aA_o$-increasing function, and let $h_e:V\to \RR$ be an $\aA_e$-increasing function. In order to prove Theorem~\ref{thm:increasing-primal} it suffices to prove that the vertex placement $\al:v\mapsto (h_o(v),h_e(v))$ is valid for every face of $G$. 

It is clear that $\alpha$ is valid for the outer face of $G$. Indeed $\al(v_1)=(\min(h_o),\min(h_e))$, $\al(v_2)=(\min(h_o),\max(h_e))$, $\al(v_3)=(\max(h_o),\max(h_e))$, and $\al(v_4)=(\max(h_o),\min(h_e))$.

Next we show that $\alpha$ is valid for every inner face of degree 4. Such a face $f$ has its corners labeled 1, 2, 3 and 4 in clockwise order in $\cL$. Let $u_1$, $u_2$, $u_3$ and $u_4$ be the vertices incident to the corners labeled 1, 2, 3 and 4 respectively. We claim that 
\begin{equation}\label{eq:separation-square}
\max(h_o(u_1),h_o(u_2))<\min(h_o(u_3),h_o(u_4))\textrm{ and } \max(h_e(u_1),h_e(u_4))<\min(h_e(u_2),h_e(u_3)).
\end{equation}
Let us justify the inequality for $h_o$. The situation is represented in Figure~\ref{fig:planarity}(a).
Since $h_o$ is $\aA_o$-increasing, and the arcs $(u_1,u_4)$ and $(u_2,u_3)$ are in $A_o$, we get $h_o(u_1)<h_o(u_4)$ and $h_o(u_2)<h_o(u_3)$. It remains to check the inequalities $h_o(u_1)<h_o(u_3)$ and $h_o(u_2)<h_o(u_4)$. If the edge $\{u_1,u_2\}$ belongs to $\ovAo$, then $h_o(u_1)=h_o(u_2)$ hence these inequalities hold. Similarly if $\{u_3,u_4\}\in\ovAo$, then the inequalities hold. So we need only consider the situation where all the edges incident to $f$ are in the support of $A_o$. 
If either of the arcs $(u_1,u_2)$ or $(u_4,u_3)$ is in $A_o$ then we get $h_o(u_1)<h_o(u_3)$ from these arcs, while if neither of these arcs are in $A_o$ then the diagonal arc $(u_1,u_3)$ belongs to $\aA_o$ which ensures $h_o(u_1)<h_o(u_3)$; see Figure~\ref{fig:diagonal-edges}. Similarly, one checks $h_o(u_2)<h_o(u_4)$ in all cases. The proof of the inequality concerning $h_e$ is similar.
 
By \eqref{eq:separation-square} there is a vertical line separating the placements of $u_1$, $u_2$ from the placement of $u_3$ and $u_4$, and an horizontal line separating the placements of $u_1$, $u_4$ from the placement of $u_2$ and $u_3$ as indicated in Figure~\ref{fig:planarity}(a). 
This easily implies that the vertex placement $\alpha$ is valid for the face $f$.

It remains to show that the placement $\al$ is valid for every inner face $f$ of $G$ of degree $3$. There are 4 possibilities for the labels of the corners of $f$: it can be missing the label 1, 2, 3 or 4. Suppose for concreteness that the missing label is 4, and let $u_1$, $u_2$ and $u_3$ be the vertices incident to the corners of $f$ labeled 1, 2 and 3, respectively. We claim that 
\begin{equation}\label{eq:separation-triangle}
h_o(u_1)\leq h_o(u_2)<h_o(u_3)\textrm{ and }h_e(u_1)< h_e(u_3)\leq h_e(u_2).
\end{equation}
Let us justify the inequalities for $h_o$. The situation is represented in Figure~\ref{fig:planarity}(b). First, the arc $(u_2,u_3)$ is in $A_o$ (since $\rterm(u_2,u_3)=3$), hence $h_o(u_2)<h_o(u_3)$. Second, Condition (L3) of 4-GS labelings implies that the label $\lterm(u_1,u_2)$ is not equal to 2, hence this label must be ether 3 or 4. Therefore, there are only 3 possibilities for the labels at the left of the arc $(u_1,u_2)$: 
\begin{compactitem}
\item either $\linit(u_1,u_2)=4$ and $\lterm(u_1,u_2)=3$, in which case $\{u_1,u_2\}\in \ovAo$ and $h_o(u_1)=h_o(u_2)$;
\item or $\linit(u_1,u_2)=1$ and $\lterm(u_1,u_2)=3$, in which case $(u_1,u_2)\in A_o$ and $h_o(u_1)<h_o(u_2)$;
\item or $\linit(u_1,u_2)=1$ and $\lterm(u_1,u_2)=4$, in which case $(u_1,u_2)\in A_o$ and $h_o(u_1)<h_o(u_2)$.
\end{compactitem}
Hence $h_o(u_1)\leq h_o(u_2)$ in all cases. The inequalities for $h_e$ is justified in a similar manner after observing that Condition (L3) of 4-GS labelings implies $\lterm(u_2,u_3)\in \{1,4\}$. By \eqref{eq:separation-triangle}, the placement of vertices are as represented at the bottom of Figure~\ref{fig:planarity}(b), which clearly implies that the vertex placement $\al$ is valid for $f$.

\fig{width=\linewidth}{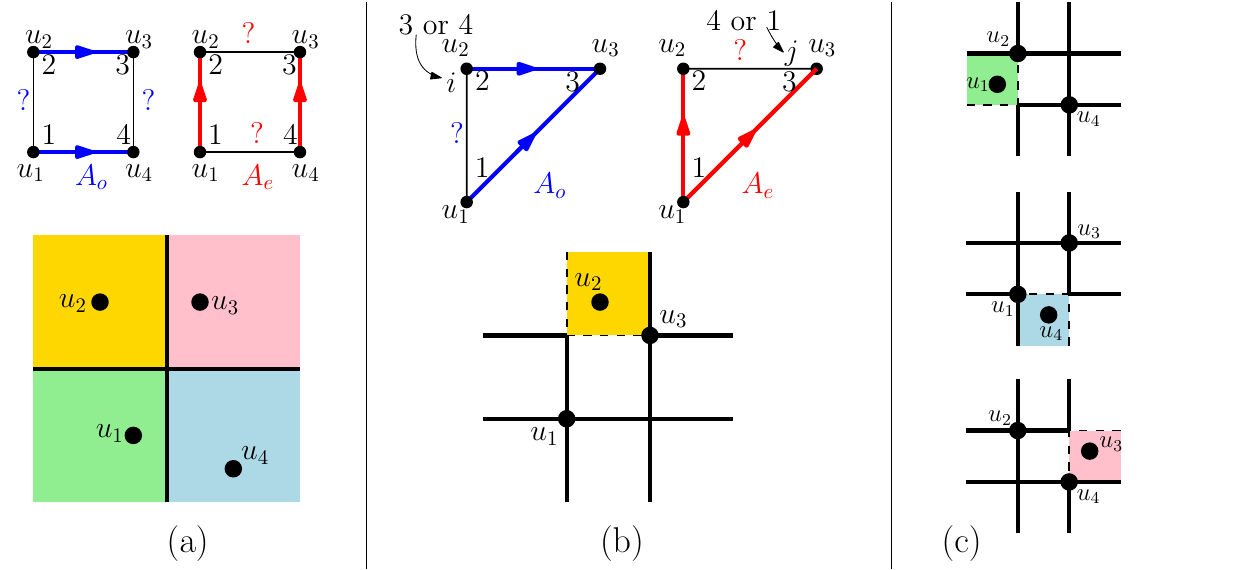}{Proof of planarity. (a) Top: orientations $A_o$ and $A_e$ for a face $f$ of degree 4. 
Bottom: Vertex placement for the face $f$ of degree 4. 
(b) Top: orientations $A_o$ and $A_e$ for a face $f$ of degree 3 having corner labels 1, 2 and~3. 
Vertex placement for the face $f$. The shaded region corresponds to the placement of vertex $u_2$ (the dotted lines indicates the weak inequalities).
(c) Vertex placement for a face of degree 3 with missing corner label 3 (top), 2 (middle), or~1 (bottom).}

A similar reasoning applies to the faces of degree 3 with other labeling situations as represented in Figure~\ref{fig:planarity}(c): compared to the case treated above the geometric constraints for the case of a face of degree 3 with missing label 3 (resp. 2, 1) is obtained by a rotation of $90$ degrees (resp. 180 degrees, 270 degrees). This completes the proof that the vertex placement $\al$ is valid for all inner faces. Hence, by Lemma~\ref{lem:planar-criteria}, $\al$ induces a planar drawing. This proves Theorem~\ref{thm:increasing-primal}.

\begin{remark}\label{rk:true-bipolar}
Any drawing of a \plm $G$ using the increasing function algorithm of Theorem~\ref{thm:increasing-primal} completely determines the 4-GS labeling used to obtain the drawing. Indeed, from the above discussion on the different types of faces of $G$ and the corresponding vertex placements, it is not hard to see that the labels of the corners inside each inner face can be determined from the placement of the vertices incident to this face.
\end{remark}


\subsection{Time complexity of the increasing function drawing algorithm}\label{sec:time-complexity}
In this subsection we discuss the time complexity of the increasing function drawing algorithm presented in Section~\ref{sec:increasing-primal}.
We say that an algorithm for 3,4-angulations can be performed in \emph{linear time} if it can be done in a number of arithmetic operations which is linear in the number of vertices. Note that the numbers of edges and faces of a 3,4-angulation are linear in the number of vertices by the Euler relation combined with the face-edge incidence relation. 

Let $G=(V,E)$ be an adapted 3,4-angulation of the square. We aim to draw $G$ with vertices placed at integer coordinates. 
Using the increasing function drawing algorithm for this task amounts to: 
\begin{compactitem}
\item[(i)] computing a 4-GS labeling $\cL$, and the corresponding augmented bipolar contractions $\aA_o$, $\aA_e$, and 
\item[(ii)] computing some integer-valued $\aA_o$-increasing and $\aA_e$-increasing functions $h_o$ and $h_e$. 
\end{compactitem}
For step (ii) we aim for a small grid size, so we want the increasing function $h_o$, $h_e$ to be \emph{tight}. Precisely, an $\aA_o$-increasing function $h_o$ is called \emph{tight} if it takes non-negative integer values and minimizes $\max(h_o)$ among such functions. Tightness for $h_e$ is defined similarly.

Step (i) can be completed in linear time. Indeed, a linear time algorithm for finding a 4-GS labeling $\cL$ was given in~\cite{OB-EF-SL:Grand-Schnyder}. Moreover it is easy to see that computing the orientations $\aA_o$, $\aA_e$ from $\cL$ can be done in linear time.

We now show that Step (ii) can also be completed in linear time. Given an acyclic directed graph $D$, we call \emph{source-level} of a vertex $v$ of $D$ the maximal length of the directed paths ending at $v$ (so that the sources of $D$ are at source-level 0). As we now explain, the source-level of every vertex in a digraph $D=(V_D,E_D)$ can be computed in time linear in $|V_D|+|E_D|$. Recall that the \emph{topological ordering algorithm} on $D$ returns a total order $\cO$ of the vertex set $V_D$ such that every arc of $D$ is oriented from the smaller vertex to the greater vertex for the order $\cO$. It is known that the topological ordering algorithm is linear in $|V_D|+|E_D|$. Once the total order $\cO$ is fixed, one can compute the source-level of every vertex, by treating vertices in the order given by the total order $\cO$ (from smallest to greatest):
for every vertex $v$ the source-level of $v$ is 0 if $v$ is a source and 
$$\textrm{source-level}(v)=\max\{\textrm{source-level}(u)+1 \mid (u,v)\in E_D\}$$
otherwise. This computation is also clearly linear in $|V_D|+|E_D|$.


Consider now the function $h_o:V\to \NN$ defined by setting $h_o(v)$ to be the source-level of the vertex corresponding to $v$ in $\aA_o$ for all $v\in V$. It is easy to see that $h_o(v)$ is $\aA_o$-increasing (recall that $\aA_o$-increasing functions can be identified with functions on the vertex set of $\aA_o$ which are increasing along the arcs of $\aA_o$). Moreover, the increasing function $h_o$ is tight. Indeed, any $\aA_o$-increasing function $f:V\to \NN$ satisfies 
$$\forall v\in V,~h_o(v)\leq f(v).$$ 
This property is easily shown by induction on the source-level of $v$ (or rather the source-level of the vertex corresponding to $v$ in $\aA_o$). Indeed we can identify $h_o$ and $f$ with some functions $h_o'$ and $f'$ on the vertex set $V'$ of $\aA_o$, and we want to show $h_o'(v)\leq f'(v)$ for all $v\in V'$. For the source $v_0$ of $\aA_o$, one has $h_o'(v_0)=0\leq f'(v_0)$. And for any vertex $v\in V'$ at source-level $k>0$, there is an arc $(u,v)$ of $\aA_o$ such that $u$ is at source-level $k-1$, hence the induction hypothesis gives $h_o'(v)=h_o'(u)+1\leq f'(u)+1\leq f'(v)$, which completes the proof.
Hence $h_o$ is a tight $\aA_o$-increasing function which can be computed in linear time.

Similarly the function $h_e(v):V\to \NN$ defined by setting $h_e(v)$ to be the source-level of the vertex corresponding to $v$ in $\aA_e$ is a tight $\aA_e$-increasing function which can be computed in linear time. We summarize our findings as follows:

\begin{prop}
The increasing function algorithm for drawing an adapted 3,4-angulation of the square can be completed in linear time in the number of vertices. 
\end{prop}

\subsection{Bounds on the grid size}\label{sec:grid-size}
In this subsection we discuss the size of the grid needed for our straight-line grid-drawing algorithms. 

Let $G$ be an adapted 3,4-angulation of the square, and let $\cL$ be a 4-GS labeling. Let $h_o,h_e$ be the tight increasing functions computed by the increasing-function algorithm. By definition the grid required for the vertex placement has size $m_o\times m_e$, where $m_o=\max(h_o)$ and $m_e=\max(h_e)$. We want to give some upper bounds on $m_o$ and $m_e$.

Since the face-counting algorithm presented in Section~\ref{sec:face-counting-primal} produces some increasing functions $\ell_e$ and $r_o$, we get the following upper bounds for the dimension of the grid: $m_o\leq \max(\ell_e)$ and $m_e\leq \max(r_o)$. By definition, $\max(r_o)$ is the number of inner faces in the bipolar orientation $B_o$. Using the Euler relation, one gets $\max(r_o)=\ee_o-\vv+1$, where $\vv$ is the number of inner vertices of $G$, and $\ee_o$ is the number of inner edges in $B_o=W_1\cup W_3^{-}$. Next, we observe that $\ee_o=2\vv-\dd_o$, where $\dd_o$ is the number of inner edges of $G$ carrying both an arc in $W_1$ and an arc in $W_3$. Thus $\max(r_o)=\vv+1-\dd_o$. Similarly $\max(\ell_e)=\vv+1-\dd_e$ where $\dd_e$ is the number of inner edges of $G$ carrying both an arc in $W_2$ and an arc in $W_4$. We summarize:


\begin{prop} \label{prop:grid-bound}
Let $G$ be an adapted 3,4-angulation of the square, and let $\cL$ be a 4-GS labeling.
The increasing function algorithm provides a grid drawing of $G$ on a grid of size $m_o\times m_e$, with $m_o\leq \vv+1-\dd_e$ and $m_e\leq \vv+1-\dd_o$, where $\dd_o$ (resp. $\dd_e$) is the number of inner edges of $G$ carrying both an arc in $W_1$ and an arc in $W_3$ (resp. an arc in $W_2$ and an arc in $W_4$). 
\end{prop}

In the case where $G$ is a triangulation of the square, $\dd_e+\dd_o$ is the number $\ee_2$ of inner edges having 2 colors in the 4-GS woods (because each arc has at most one color for any 4-GS structure on any triangulation of the square). Let $\ee_0,\ee_1,\ee_2$ be the number of inner edges having $0,1,2$ colors respectively, and let $\ee=\ee_0+\ee_1+\ee_2$ be the total number of inner edges. The Euler relation (combined with the incidence relation between faces and edges) gives $\ee=3\vv+1$, and the fact that each inner vertex is incident to one outgoing arc of each color gives
$\ee_1+2\ee_2=4\vv$. This implies $\ee_2=\vv-1+\ee_o$, and $m_o+m_e=2\vv+2-\ee_2=\vv+3-\ee_0\leq \vv+3$. Thus, for a triangulation of the square, the half-perimeter of the grid drawing is strictly less than the number of vertices.

\subsection{Reduction of the grid size}\label{sec:optimization}
In this subsection we discuss some optimizations of the increasing function drawing algorithm aimed at reducing the grid size. These improvements will be obtained by relaxing some requirements related to diagonal arcs for the increasing functions $h_o,h_e$.

In order to state which requirements can be relaxed, we first need to classify the diagonal arcs of $\aA_o$ and $\aA_e$. We refer to  Figure~\ref{fig:diagonal-edges} for the definition of the diagonal arcs.
For a diagonal arc $a$ of $A_o$ or $A_e$, we denote by $f(a)$ the inner face of $G$ containing $a$. Each of the 4 edges of $G$ incident to the face $f(a)$ can be either \emph{missing} if they are in $\ovAo\cup \ovAe$, or \emph{converging} if they are oriented the same way in $A_o$ and $A_e$, or \emph{diverging} if they are oriented in opposite directions in $A_o$ and $A_e$.
We let $\miss(a)$, $\conv(a)$ and $\dive(a)$ be the number of missing, converging and diverging edges incident to $f(a)$. Note that $\miss(a)\in \{0,1,2\}$ and $\miss(a)+\conv(a)+\dive(a)=4$.

We call all the diagonal arcs of $\aA_o$ and $\aA_e$ \emph{weak}, except for the diagonal arcs $a$ of $\aA_e$ such that $\conv(a)=4$ or $\dive(a)=4$ which are called \emph{strict}. We call a diagonal arc $a$ \emph{mandatory} unless $|\conv(a)-\dive(a)|=\miss(a)$. It is easy to check that a diagonal arc $a$ is non-mandatory if and only if it satisfies one of the following conditions:
\begin{compactitem}
\item[(i)] $\miss(a)=0$ and $\conv(a)=2$, 
\item[(ii)] $\miss(a)=1$ and $\conv(a)\in \{1,2\}$,
\item[(iii)] $\miss(a)=2$.
\end{compactitem}
Indeed, when $\miss(a)=0$ the condition $|\conv(a)-\dive(a)|=0$ is equivalent to $\conv(a)=2$, (since $\conv(a)+\dive(a)=4$), while when $\miss(a)=2$ one can check that $\conv(a)\in \{0,2\}$ (see Figure~\ref{fig:diagonal-edges} and observe that the missing edges are either both in $\ovAo$ or both in $\ovAe$) hence  $|\conv(a)-\dive(a)|=2$ (since  $\conv(a)+\dive(a)=2$).
 
Note that a strict diagonal arc is always mandatory.

\begin{example} 
For the diagonal arc $a\in \aA_o$ on the left of Figure~\ref{fig:increasing-functions}(c) one has $\miss(a)=0$, $\conv(a)=1$, and $\dive(a)=3$, hence this arc is weak and mandatory.
For the  diagonal arc $b\in \aA_e$ on the right of Figure~\ref{fig:increasing-functions}(c) one has $\miss(b)=2$, $\conv(a)=0$, and $\dive(a)=2$, hence this arc is weak and non-mandatory.
\end{example}

We say that a function $g:V\to\RR$ is \emph{quasi $\aA_o$-increasing} (resp. \emph{quasi $\aA_e$-increasing}) if it satisfies the following 3 conditions:
\begin{compactitem}
\item $g(u)=g(v)$ for every edge $\{u,v\}$ in $\ovAo$ (resp. $\ovAe$),
\item $g(u)\leq g(v)$ for every weak mandatory diagonal arc $(u,v)$ in $\aA_o$ (resp. $\aA_e$),
\item $g(u)<g(v)$ for every arc $(u,v)$ in $A_o$ (resp. $A_e$) and every strict mandatory  diagonal arc $(u,v)$ in $\aA_o$ (resp. $\aA_e$).
\end{compactitem}
Clearly any $\aA_o$-increasing function is quasi $\aA_o$-increasing (and similarly for $\aA_e$). We can now state a slightly stronger version of Theorem~\ref{thm:increasing-primal}.

\begin{thm}\label{thm:increasing-primal-optimized}
Let $G$ be an adapted 3,4-angulation of the square, let $\cL$ be a 4-GS labeling, and let $\aA_o$ and $\aA_e$ be the corresponding augmented bipolar contractions of $G$. 
Let $h_o$ be a quasi $\aA_o$-increasing function and $h_e$ be a quasi $\aA_e$-increasing function.
Placing each vertex $v$ of $G$ at coordinates $(h_o(v),h_e(v))\in \RR^2$ gives a planar straight-line drawing of $G$.    
\end{thm}

Theorem~\ref{thm:increasing-primal-optimized} is illustrated in Figure~\ref{fig:increasing-functions}(e). In that figure the vertex placement is given by some quasi increasing functions $h_o$ and $h_e$. Note that in that example, the function $h_o$ is not $\aA_o$-increasing since $h_o(a)=h_o(b)$ and the function $h_e$ is not $\aA_e$-increasing since $h_e(c)=h_o(f)$.

\begin{proof}[Proof of Theorem~\ref{thm:increasing-primal-optimized}]
Our proof of Theorem~\ref{thm:increasing-primal-optimized} relies on the same line of arguments as for Theorem~\ref{thm:increasing-primal}. Compared to Theorem~\ref{thm:increasing-primal}, some requirements associated to the diagonal arcs of $\aA_o$ and $\aA_e$ have been relaxed, and we need to check that the proof remains valid. In the proof of Theorem~\ref{thm:increasing-primal}, the requirements associated to a diagonal arc $a$ for a placement function $\al=(h_o,h_e)$ are used solely to prove that the placement $\al$ is valid for the face $f(a)$ of $G$. Hence, we only need to prove the following two claims. 
\begin{compactitem}
\item[(i)] For any non-mandatory diagonal arc $a=(u,v)$, the placement $\al$ is valid for $f(a)$ even if the requirement on $h_o$ or $h_e$ associated to $a$ is dropped.
\item[(ii)] For any weak diagonal arc $a=(u,v)$, the placement $\al$ is valid for $f(a)$ even if the strict inequality requirement on $h_o$ or $h_e$ associated to $a$ is replaced by a weak inequality.
\end{compactitem}

Let us start by proving Claim (i). Let $a$ be a non-mandatory arc of $\aA_o$, and let $u_1,u_2,u_3,u_4$ be the vertices of $G$ incident to the corners of $f(a)$ labeled 1,2,3,4 respectively (see again Figure \ref{fig:diagonal-edges}). Recall that $(u_2, u_3)$ and $(u_1, u_4)$ are in $A_o$. There are two cases to check: either $a=(u_1,u_3)$ or $a=(u_2,u_4)$. Let us start by assuming that $a=(u_1,u_3)$ as in Figure~\ref{fig:planarity-optimization}(b). In this situation we know
$$\{(u_2,u_1),(u_1,u_4),(u_2,u_3),(u_3,u_4)\}\subseteq A_o~\textrm{ and }~\{(u_1,u_2),(u_4,u_3)\}\subseteq A_e,$$
and since $a$ is non-mandatory there are only four possible configurations for the orientation of $f(a)$ in $\aA_e$, which are represented in Figure~\ref{fig:planarity-optimization}:
\begin{compactitem}
\item $(u_1,u_4)\in A_e$ and $(u_2,u_3)\in A_e$ (this gives $\miss(a)=0$ and $\conv(a)=2$),
\item or $(u_1,u_4)\in A_e$ and $\{u_2,u_3\}\in \ovAe$ (this gives $\miss(a)=1$ and $\conv(a)=1$),
\item or $\{u_1,u_4\}\in \ovAe$ and $(u_2,u_3)\in \aA_e$ (this gives $\miss(a)=1$ and $\conv(a)=1$),
\item or $\{u_1,u_4\}\in \ovAe$ and $\{u_2,u_3\}\in \ovAe$ (this give $\miss(a)=2$ and $\conv(a)=0)$.
\end{compactitem}
In all four configurations, the quasi increasing functions $h_o$ and $h_e$ satisfy 
\begin{equation}\label{eq:placement-non-mandatory}
h_o(u_2)< h_o(u_1),h_o(u_3)< h_o(u_4)~\textrm{ and }~h_e(u_1)\leq h_e(u_2),h_e(u_4)\leq h_e(u_3).
\end{equation}
These inequalities (together with the observation $h_e(u_1)< h_e(u_3)$) are enough to ensure that the placement $\al=(h_o,h_e)$ is valid for the face $f(a)$ as claimed; see Figure~\ref{fig:planarity-optimization}(c). 

\fig{width=\linewidth}{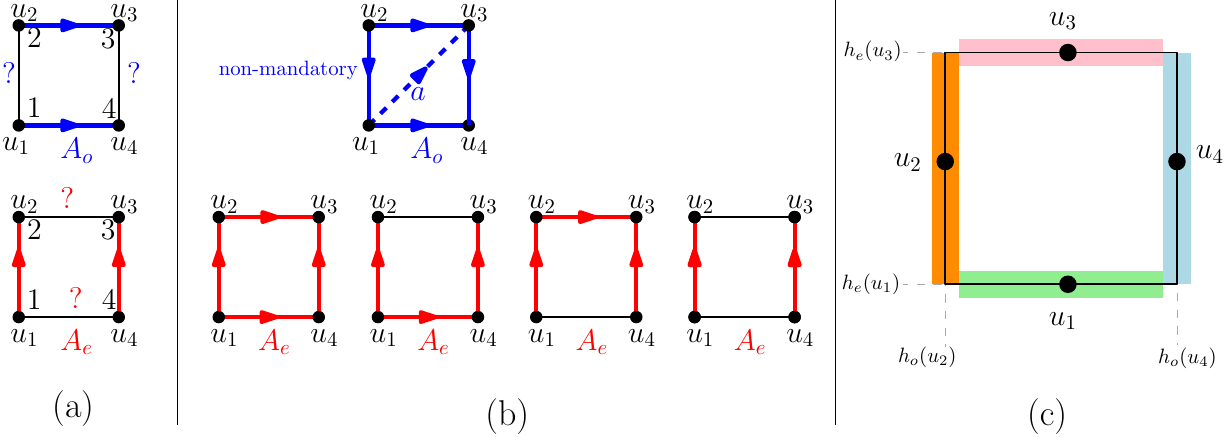}{Proof of Claim (i) for a non-mandatory arc $a$ of $\protect\aA_o$. (a) The face $f(a)$ and the arcs which are necessarily present in $A_o$ and $A_e$. (b) The possible configurations for $A_o$ and $A_e$ assuming that $a=(u_1,u_3)$ is a non-mandatory arc $a$ of $\protect\aA_o$. (c) Vertex placement $\al$ for the face $f(a)$ given by the inequalities in \eqref{eq:placement-non-mandatory}: the shaded regions indicate where each vertex could be placed by $\al$. In every case the placement $\al$ is valid for $f(a)$.}

Similarly, in the situation $a=(u_2,u_4)$ we have 
$$\{(u_1,u_2),(u_2,u_3),(u_1,u_4),(u_4,u_3)\}\subseteq A_o~\textrm{ and }~\{(u_1,u_2),(u_4,u_3)\}\subseteq A_e,$$
and either $(u_1,u_2)\in A_e$ or $\{u_1,u_2\}\in \ovAe$ and either $(u_4,u_3)\in A_e$ or $\{u_4,u_3\}\in \ovAe$. This gives 
$$h_o(u_1)< h_o(u_2),h_o(u_4)< h_o(u_3)~\textrm{ and }~h_e(u_4)\leq h_e(u_1),h_e(u_3)\leq h_e(u_2),$$
which (together with the observation $h_e(u_4)< h_e(u_2)$) is enough to ensure that the placement $\al$ is valid for the face $f(a)$.
The situation where $a$ is a non-mandatory arc of $\aA_e$ is similar.

It remains to prove Claim (ii). Let $a$ be weak diagonal arc of $\aA_o$ or $\aA_e$ and let $u_1,u_2,u_3,u_4$ be the vertices of $G$ incident to the corners of $f(a)$ labeled 1,2,3,4 respectively. Since we have already proved Claim (i) we can assume that any diagonal arc in the face $f(a)$ is mandatory, hence the quasi increasing functions $h_o$ and $h_e$ satisfy 
\begin{equation*}
h_o(u_1),h_o(u_2)\leq h_o(u_3),h_o(u_4)~\textrm{ and }~h_e(u_1),h_e(u_4)\leq h_e(u_2),h_e(u_3).
\end{equation*}
In fact, there are 3 possible situations for $h_o$ (since $h_o(u_1)<h_o(u_4)$ and $h_o(u_2)<h_o(u_3)$):
\begin{compactitem} 
\item[(a$_o$)] either $h_o(u_1),h_o(u_2)< h_o(u_3),h_o(u_4)$,
\item[(b$_o$)] or $h_o(u_2)<h_o(u_1)=h_o(u_3)<h_o(u_4)$ (which can happen if $a=(u_1,u_3)$ is a weak diagonal arc of $\aA_o$), 
\item[(c$_o$)] or $h_o(u_1)<h_o(u_2)=h_o(u_4)<h_o(u_3)$ (which can happen if $a=(u_2,u_4)$ is a weak diagonal arc of $\aA_o$).
\end{compactitem}
Similarly, there are 3 possible situations for $h_e$:
\begin{compactitem} 
\item[(a$_e$)] either $h_e(u_1),h_e(u_4)< h_e(u_2),h_e(u_3)$,
\item[(b$_e$)] or $h_e(u_4)<h_e(u_1)=h_e(u_3)<h_e(u_2)$ (which can happen if $a=(u_1,u_3)$ is a weak diagonal arc of $\aA_e$), 
\item[(c$_e$)] or $h_e(u_1)<h_e(u_2)=h_e(u_4)<h_e(u_3)$ (which can happen if $a=(u_4,u_2)$ is a weak diagonal arc of $\aA_e$).
\end{compactitem}
It is easy to see that if (a$_o$) or (a$_e$) holds, then the placement $\al$ is valid for $f(a)$ as illustrated in Figure~\ref{fig:planarity-optimization2}(b).
If both (b$_o$) and (c$_e$) hold, then the placement $\al$ is valid for $f(a)$ as illustrated in Figure~\ref{fig:planarity-optimization2}(c).
Similarly, if both (b$_e$) and (c$_o$) hold, then the placement $\al$ is valid for $f(a)$.

\fig{width=\linewidth}{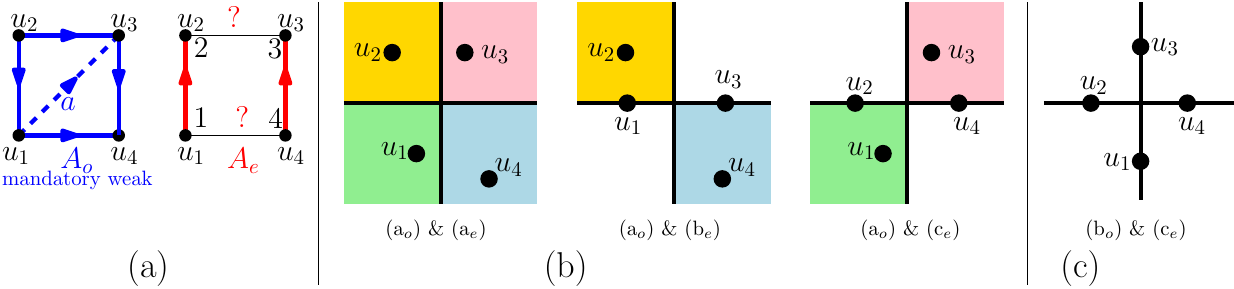}{Proof of Claim (ii) for a weak diagonal arc $a=(u_1,u_3)$ of $\protect\aA_o$. (a) The face $f(a)$ and the arcs which are necessarily present in $A_o$ and $A_e$. (b) Vertex placement $\al$ for the face $f(a)$ assuming (a$_o$) holds. (c) Vertex placement for $f(a)$ assuming both (b$_o$) and (c$_e$) hold.}

The only remaining situations are when both (b$_o$) and (b$_e$) hold, or both (c$_o$) and (c$_e$) hold. However these situations are incompatible with our convention for weak diagonal arcs. For instance, in order for both (b$_o$) and (b$_e$) to hold, the arc $(u_1,u_3)$ needs to be a diagonal arc of both $\aA_o$ and $\aA_e$; but this implies $\conv(a)=4$ so that $(u_1,u_3)$ is a strict diagonal arc of $\aA_e$, which contradicts (b$_e$). This completes the proof of Claim (ii), and Theorem~\ref{thm:increasing-primal-optimized}.
\end{proof}

Let us call \emph{optimized} increasing function drawing algorithm, the grid drawing algorithm obtained from Theorem~\ref{thm:increasing-primal-optimized}. Just as the increasing function algorithm discussed in Section~\ref{sec:time-complexity}, this optimized version can be performed in linear time. Indeed, starting for an adapted 3,4-angulation of the square $G$, the computation of the augmented contractions $\aA_o$ and $\aA_e$ can be done in linear time as explained above. Removing the non-mandatory diagonal arcs, and identifying the weak diagonal arcs can clearly be done in linear time. It then remains to compute ``tight'' quasi increasing functions $h_o$ and $h_e$.

Given an acyclic directed graph $D=(V_D,E_D)$ with arcs marked as either \emph{weak} or \emph{strict}, we call \emph{strict source-level} of a vertex $v\in V_D$, the maximal number of strict arcs on directed paths ending at $v$. The strict source-level of every vertex can be computed in time linear in $|V_D|+|E_D|$. Indeed, one can run the \emph{topological ordering algorithm} on $D$ to get a total order $\cO$ on the vertex set $V_D$ (such that every arc of $D$ is oriented from the smaller vertex to the greater vertex for the order $\cO$). Then, one can compute the strict source-level of every vertex, by treating vertices in the order given by the total order $\cO$ (from smallest to greatest):
for every vertex $v\in V_D$ the strict source-level of $v$ is 0 if $v$ is a source, and 
$$\textrm{strict-level}(v)=\max\big(\{\textrm{strict-level}(u)\mid(u,v)\in E_D \textrm{ weak}\}\cup\{\textrm{strict-level}(u)+1\mid (u,v)\in E_D\textrm{ strict}\}\big)$$
otherwise. 

It is easy to see that the strict source-levels for the digraphs $\aA_o$ and $\aA_e$ (after removal of the non-mandatory diagonal arcs) give a tight quasi $\aA_o$-increasing function $h_o$ and a tight quasi $\aA_e$-increasing function $h_e$. Hence we get the following result.


\begin{prop}
The optimized increasing function algorithm for drawing an adapted 3,4-angulation of the square can be completed in linear time. 
\end{prop}

\section{Orthogonal drawing of \plms with vertex-degree at most 4}\label{sec:dual_algo}



In this section we turn our attention to the dual picture: we will define orthogonal drawing algorithms for dual-adapted rooted 3,4-maps, which are the maps dual to the adapted 3,4-angulations of the square.

Recall that a \emph{planar orthogonal drawing} of a graph is a planar drawing where each edge is represented as a sequence of horizontal and vertical line segments. For convenience (and in view of displaying drawings with at most one bend per edge), we adopt the convention of not drawing the root-vertex $v_\infty$, as in~\cite{OB-EF:Schnyder,FelsnerKV14}. 
\begin{definition}
Let $G^*$ be a rooted 3,4-map, let $\vinfty$ be its root vertex, and let $e_1^*,\ldots,e_4^*$ be the root edges.
A \emph{suspended planar orthogonal drawing} for $G^*$, or \emph{SPO drawing} for short, is a planar orthogonal drawing of the map obtained from $G^*$ by disconnecting the edges $e_1^*,\ldots,e_4^*$ from $\vinfty$, in which the edges $e_1^*,\ldots,e_4^*$ are drawn as dangling arrows in the outer face whose directions are \emph{west}, \emph{north}, \emph{east} and \emph{south} respectively\footnote{This type of orthogonal drawing is called \emph{anchored} in~\cite{FelsnerKV14} (which considers non-necessarily planar drawings).}.
\end{definition}
For instance, the drawing of Figure \ref{fig:face_counting_dual_big}(b)
is an SPO drawing.


The algorithms we present in this section are based on \emph{dual 4-GS structures}, which were introduced in~\cite{OB-EF-SL:Grand-Schnyder}. 
We first define a ``face-counting'' SPO drawing algorithm which generalizes Bernardi and Fusy's algorithm~\cite{OB-EF:Schnyder}, and then a more general (and more compact) ``increasing functions'' SPO drawing algorithm, which also  generalizes He's algorithm~\cite{He93:reg-edge-labeling}.


\subsection{Face-counting drawing algorithm for 3,4-maps}\label{sec:face-counting-dual}
In this subsection we define our first SPO drawing algorithm for rooted 3,4-maps. In this algorithm, the coordinates of each vertex are obtained by counting the number of faces in certain regions of the odd and even dual bipolar orientations. We shall refer to this algorithm as the \emph{face-counting SPO drawing algorithm}.

Let $G^*$ be a dual-adapted rooted 3,4-map. 
Let $\mathcal{L}^*$ be a dual 4-GS labeling of $G^*$ satisfying the additional requirement that no edge is \emph{fully-colored}, that is, belongs to all four trees of the dual 4-GS wood $\cW^* = \{W_1^*,W_2^*,W_3^*,W_4^*\}$ associated with $\cL^*$. We will explain the necessity of this additional requirement in Remark~\ref{rem:uncolored}. For now let us briefly justify the existence of such a structure. 
\begin{lemma}
Any dual-adapted rooted 3,4-map admits a dual 4-GS labeling without any fully-colored edge. Moreover such a structure can be found in linear time.
\end{lemma}
\begin{proof} Via duality, it suffices to prove the following equivalent property: any adapted 3,4-angulation of the square $G$ admits a 4-GS labeling $\cL$ without any uncolored edge (edge not belonging to any of the trees of the associated 4-GS wood). 
It is shown in~\cite{OB-EF-SL:Grand-Schnyder}, that the set of 4-GS woods of $G$ admits the structure of a distributive lattice. 
This lattice structure is easier to define using a different incarnation of 4-GS structures called \emph{4-GS angular orientations}. 
We refer the reader to~\cite{OB-EF-SL:Grand-Schnyder} for the definition of 4-GS angular orientations and the correspondence with 4-GS woods.
The relevant information here is that if $\cL$ has an uncolored edge $e$, then the associated angular orientation $\cA$ has a counterclockwise cycle (made of the 4 ``star edges'' surrounding $e$). Moreover if $\cL$ is the maximal element of the distributive lattice, then $\cA$ has no counterclockwise cycle, hence $\cL$ has no uncolored edge. Lastly, it is shown in~\cite{OB-EF-SL:Grand-Schnyder}, that the maximal element $\cL$ can be constructed in linear time.
\end{proof}

Recall from Section~\ref{sec:dual4GS} that 
 $$B_o^* = W_3^* \cup (W_1^*)^{-}~\textrm{ and }~B_e^* = W_2^* \cup (W_4^*)^{-}$$ 
are bipolar orientations.
For each non-root vertex $v$, we denote by $P_i^*(v)$ the directed path from $v$ to the root vertex $v_\infty$ in the tree $W_i^*$, and define the \emph{odd and even dual paths} for $v$ as 
$$P_o^*(v) = P_3^*(v) \cup P_1^*(v)^{-} \subseteq B_o^* \text{ and } P_e^*(v) = P_2^*(v) \cup P_4^*(v)^{-} \subseteq B_e^*.$$
Note that $P_o^*(v)$ is a simple directed path whose support contains the root-edges $e_3^*$ and $e_1^*$, and $P_e^*(v)$ is a simple directed path whose support contains the root-edges $e_2^*$ and $e_4^*$.

We define the \emph{face-counting coordinates} of $v$ as $(\ell_e^*(v), r_o^*(v))$, where $\ell_e^*(v)$ is the number of inner faces of $B_e^*$ on the left of $P_e^*(v)$, and $r_o^*(v)$ is the number of inner faces of $B_o^*$ on the right of $P_o^*(v)$. The face-counting coordinates are represented in Figure~\ref{fig:face_counting_dual_big}(a).

\fig{width=\linewidth}{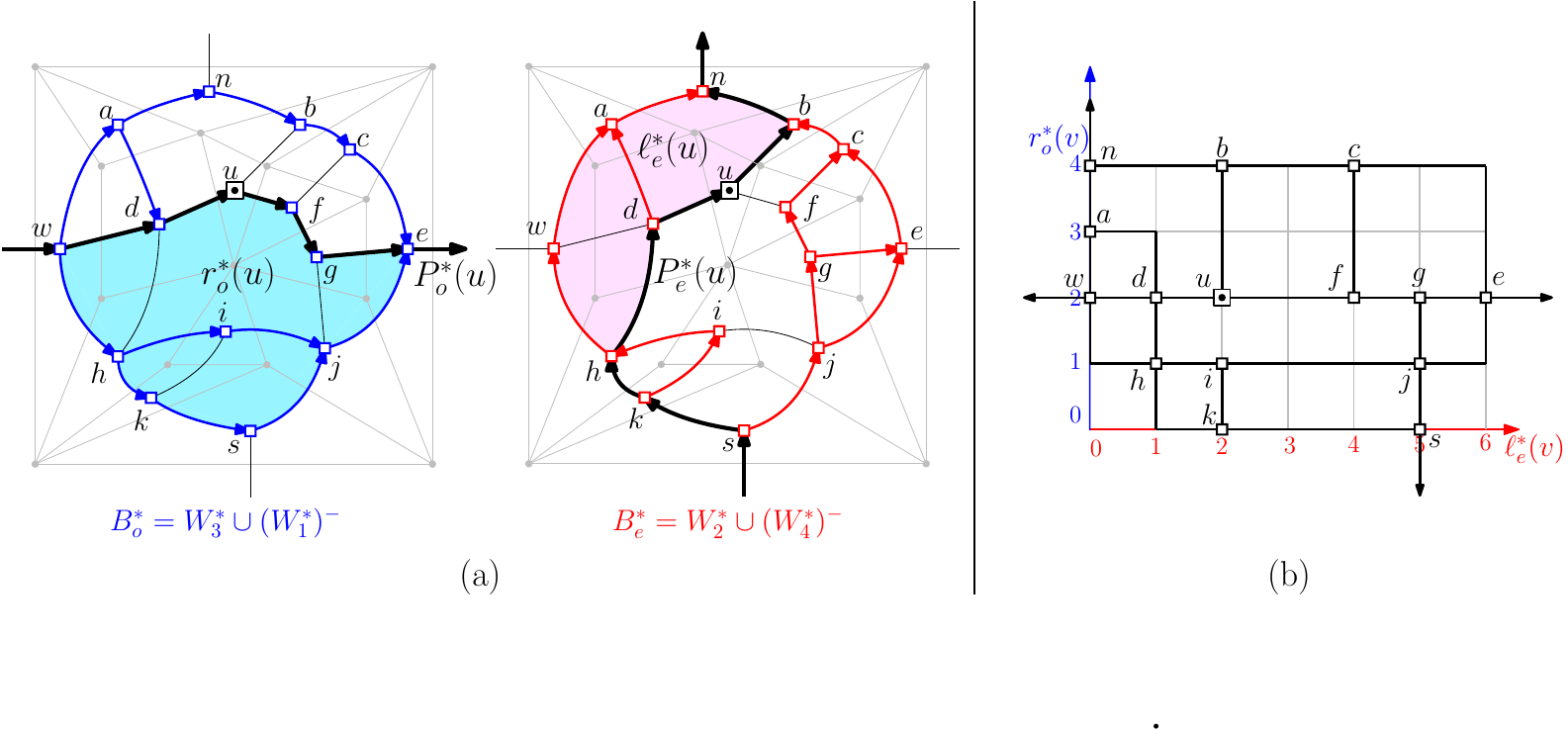}{(a) The paths and regions considered to compute the face-counting coordinates of the vertex $u$. Here $\ell_e^*(u) = 2$ and $r_o^*(u) = 2$. (b) The face-counting drawing of $G^*$.}
We now need to specify how to draw the edges of $G^*$. 
For a point $A$ in $\br^2$, we denote by 
$\gamma_1(A)$ (resp. $\gamma_2(A)$, $\gamma_3(A)$, $\gamma_4(A)$)
the half-line starting from $A$ going in the negative $x$-direction (resp. positive $y$-direction, positive $x$-direction, negative $y$-direction). See Figure~\ref{fig:bent_edge2}(a).

\begin{definition} \label{def:rule-draw-deges}
Let $G^*$ be a dual-adapted rooted 3,4-map endowed with a dual 4-GS labeling~$\cL^*$. We consider the coloring of the arcs of $G^*$ given by the 4-GS wood $(W_1^*,W_2^*,W_3^*,W_4^*)$, and recall that every non-root arc has at least one color in $\four$. 

Given a vertex placement $\alpha:V'\to \br^2$, where $V'$ is the set of non-root vertices of $G^*$, the \emph{$\cL^*$-rule} for drawing each non-root edge $e=\{u,v\}$ is as follows:
\begin{compactitem}
\item If $e$ has exactly two colors, and these colors $i,j$ are of different parity, then the edge $e$ is drawn as a bent-edge. Precisely, if the arc $(u,v)$ has color $i$ and the arc $(v,u)$ has color $j$, then $e$ is drawn as the union of line segments $[\al(u),X]\cap [X,\al(v)]$, where $X$ is the intersection of the rays $\gamma_i(\al(u))$ and $\gamma_j(\al(v))$ (assuming the intersection $X$ is non empty); see Figure~\ref{fig:bent_edge2}(b).
\item In all the other situation, $e$ is drawn as a straight edge.
\end{compactitem}
\end{definition}

\fig{width=\linewidth}{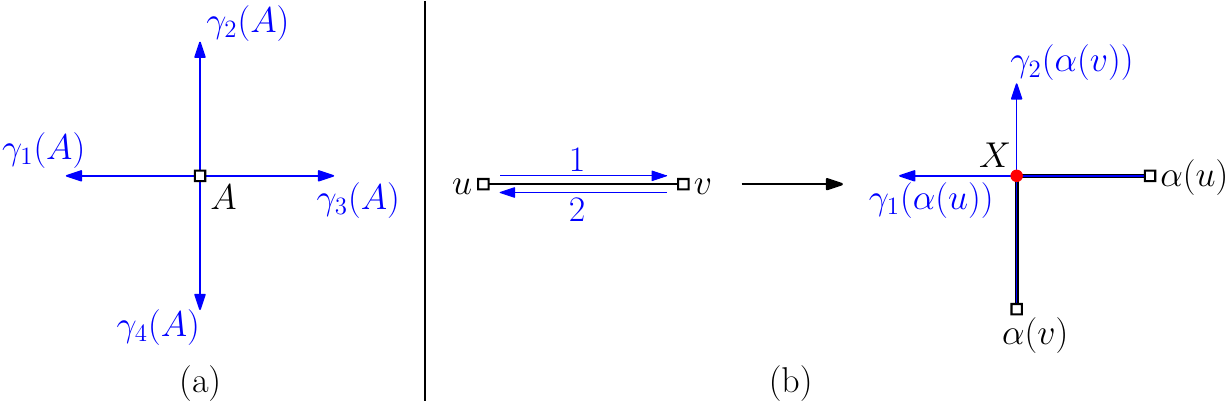}{(a) The four rays $\gamma_1(A),\gamma_2(A),\gamma_3(A),\gamma_4(A)$ emanating from a point~$A$. (b) The bent-edge drawn according to the $\cL^*$-rule for an edge $e=\{u,v\}$ having arc $(u,v)$ colored 1 and arc $(v,u)$ colored 2.}

The \emph{face-counting drawing} of $G^*$ (based on $\cL^*$) is defined as follows. 
\begin{compactitem}
\item Each non-root vertex $v$ is placed at the face-counting coordinates $(\ell_e^*(v), r_o^*(v))$.
\item Each non-root edge $e$ is drawn according to the $\cL^*$-rule (Definition~\ref{def:rule-draw-deges}).
\item The root edges $e_1^*,e_2^*,e_3^*,e_4^*$ are drawn as four dangling arrows whose directions are west, north, east, and south, respectively.
\end{compactitem}

An example is represented in Figure~\ref{fig:face_counting_dual_big}(b). Of course, one needs to justify that it is possible to draw the edges according to the $\cL^*$-rule (and that the straight edges are all either horizontal or vertical). In fact, we will prove the following stronger result.

\begin{thm}\label{thm:face_counting_dual}
 Let $G^*$ be a dual-adapted rooted 3,4-map and let $\mathcal{L}^*$ be a dual 4-GS labeling of $G^*$ such that no edge is fully-colored. 
 The face-counting drawing of $G^*$ based on $\cL^*$ is an SPO drawing.
\end{thm}

We will prove Theorem~\ref{thm:face_counting_dual} after presenting a more general algorithm in the next subsection. 

\begin{remark}\label{rem:uncolored}
Now we explain why the requirement that the dual 4-GS structure does not have fully-colored edges is necessary: if one uses a dual 4-GS labeling $\cL^*$ having a fully colored edge $e=\{u,v\}$, then the face-counting coordinates of $u$ and $v$ are equal, hence the face-counting drawing is not planar.
\end{remark}


To close this subsection, we compare the face-counting SPO drawing algorithm of Theorem~\ref{thm:face_counting_dual} to the algorithm 
from~\cite{OB-EF:Schnyder} to draw the dual $G^*$ of a simple quadrangulation $G$. 
Recall from the introduction that this algorithm relies on a corner labeling $\cL^*$ of $G^*$ obtained from a separating decomposition $\cS$ of $G$. 
As explained in Section~\ref{sec:4GS1}, the separating decomposition $\cS$ corresponds to an even 4-GS labeling $\cL$ of $G$, and it clear from the definition that the corner labeling $\cL^*$ is the dual of $\cL$. 
For each non-root vertex $v$, the path $P_i'(v)$ used by the algorithm in~\cite{OB-EF:Schnyder} is equal to the path $P_i^*(v)$ in the tree $W_i^*$, and the odd and even separating paths $P_o'(v)$ and $P_e'(v)$ are equal to the odd and even dual paths $P_o^*(v)$ and $P_e^*(v)$. 
Therefore, the placement of vertices by the algorithm in~\cite{OB-EF:Schnyder} coincides with the placement from the face-counting algorithm based on the dual 4-GS labeling. 
Furthermore, the rule for drawing edges in~\cite{OB-EF:Schnyder} also coincides with the rule described before. (A consequence of the evenness of $\cL$ property is that every non-root edge of $G^*$ has exactly two colors of different parity hence is drawn as a bent edge.) 
In conclusion, the algorithm of~\cite{OB-EF:Schnyder} is the same as the face-counting algorithm from Theorem~\ref{thm:face_counting_dual} restricted to the dual of simple quadrangulations endowed with even dual 4-GS labelings.



\subsection{Increasing function drawing algorithm for 3,4-maps}\label{sec:increasing_dual}

In this subsection we present another SPO drawing algorithm for dual-adapted rooted 3,4-maps. We shall call this algorithm the \emph{increasing function SPO drawing algorithm}. 
This algorithm generalizes the face-counting algorithm presented in the previous subsection, hence can provide some improvement on the grid size. 
It relies on some acyclic orientations that we now present.

Let $G^*$ be a dual-adapted rooted 3,4-map, and let $\cL^*$ be a dual 4-GS labeling (whose associated dual 4-GS wood may have fully-colored edges). The \emph{odd} and \emph{even dual acyclic orientations} $A_o^*$ and $A_e^*$ are defined as:
\begin{align*}
 A_o^* = &\{a \text{ non-root arc} \mid 4 \in [\text{right-init}(a),\text{right-term}(a)[ \text{ or } 2 \in [\text{left-term}(a), \text{left-init}(a)[\\
 &\text{ or } (\text{right-init}(a),\text{right-term}(a),\text{left-term}(a), \text{left-init}(a)) = (1,1,3,3)\}\\
 A_e^* = &\{a \text{ non-root arc} \mid 3 \in [\text{right-init}(a),\text{right-term}(a)[ \text{ or } 1 \in [\text{left-term}(a), \text{left-init}(a)[\\
 &\text{ or } (\text{right-init}(a),\text{right-term}(a),\text{left-term}(a), \text{left-init}(a)) = (4,4,2,2)\}
\end{align*}
The definition of $A_o^*$ and $A_e^*$ is represented in Figure~\ref{fig:a_dual_def} (top row), and a concrete example is represented in Figure~\ref{fig:bipolars_dual_big}(b). 
One can check that $A_o^*$ is a subset of the  acyclic orientation $B_o^*$, and $A_e^*$ is a subset of the acyclic orientation $B_e^*$, hence they are acyclic.

Let $\ov{A_o^*}$ be the set of edges not in the support of $A_o^*$, and let $\ovAes$ be the set of edges not in the support of $A_e^*$. From the definition of $A_o^*$ and $A_e^*$ and the conditions for dual 4-GS labelings, it is not hard to see that \begin{align*}
 \ovAos = &\{ e \text{ non-root edge} \mid \text{one side of $e$ has labels $\{1,2\}$ or $\{3,4\}$} \\ 
 &\text{or the four labels around $e$ are } \{2,2,4,4\}\},\\
 \ovAes = &\{ e \text{ non-root edge} \mid \text{one side of $e$ has labels $\{2,3\}$ or $\{4,1\}$} \\
 &\text{or the four labels around $e$ are } \{1,1,3,3\}\}.
\end{align*}

The definition of $\ovAos$ and $\ovAes$ is also represented in Figure~\ref{fig:a_dual_def} (middle and bottom rows). 
We claim that $\ovAos$ (resp. $\ovAes$) do not contain cycles. 
To see this, first note that by construction $\ovAos$ (resp. $\ovAes$) is a subset of the support of $B_e^*$ (resp. $B_o^*$). 
Moreover, it is easy to check that a cycle in $\ovAos$ (resp. $\ovAes$) always supports a directed cycle in $B_e^*$ (resp. $B_o^*$), which is impossible since we have established in Section~\ref{sec:dual4GS} that $B_e^*$ and $B_o^*$ are bipolar orientations.

\fig{width=\linewidth}{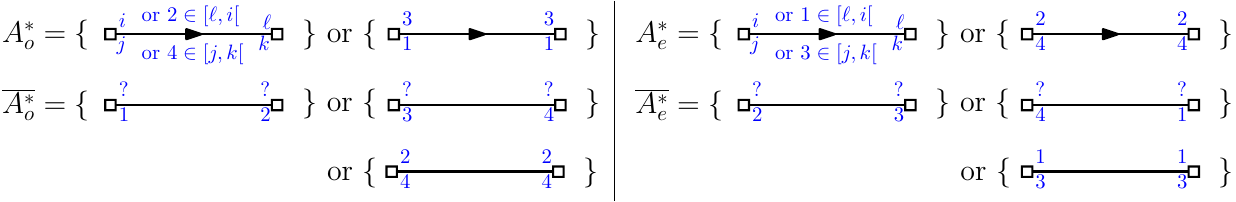}{The non-root arcs in the dual acyclic orientations $A_o^*$ and $A_e^*$, and the non-root edges not in their support.}

We now define the dual augmented bipolar contractions in two steps. The \emph{odd} (resp. \emph{even}) \emph{dual augmented bipolar contraction} $\aA_o^*$ (resp. $\aA_e^*$) is obtained from $A_o^*$ (resp. $A_e^*$) by the following procedure: 
\begin{compactitem}
 \item[(i)] For each non-root face of $G^*$, add an \emph{odd} (resp. \emph{even}) \emph{diagonal arc} from the last corner labeled 4 (resp. 3) clockwise around $f$ to the last corner labeled 2 (resp. 1) clockwise around $f$, if both corners exist. 
 \item[(ii)] Contract all the edges in $\ovAos$ (resp. $\ovAes$).
\end{compactitem}
This process is represented in Figure~\ref{fig:increasing_dual_big}(a)-(b). 
Note that Step (ii) does not depend on Step (i) and can be executed independently. In fact, it is convenient to introduce some separate notations for the structures obtained by performing Step (ii) only: we define $\wt{A}_o^* := A_o^*/\ovAos$ and $\wt{A}_e^* := A_e^*/\ovAes$, and call them the \emph{odd} and \emph{even bipolar contractions}. 
We will prove later that $\wt{A}_o^*$, $\wt{A}_e^*$, $\aA_o^*$ and $\aA_e^*$ are bipolar orientations.

An $\aA_o^*$ (resp. $\aA_e^*$)-\emph{increasing} function for $G^*$ is a function $h^*$ defined on the non-root vertices of $G^*$ that satisfies the following 2 conditions: \begin{compactitem}
\item $h^*(u) = h^*(v)$ for every edge $\{u, v\}$ in $\ovAos$ (resp. $\ovAes$),
\item $h^*(u)< h^*(v)$ for every arc $(u,v)$ in $A_o^*$ (resp. in $A_e^*$) and every odd (resp. even) diagonal arc $(u,v)$. 
\end{compactitem}
Such functions can be identified with functions defined on vertices of $\aA_o^*$ (resp. $\aA_e^*$) that are strictly increasing along arcs. 


\begin{thm}\label{thm:increasing_dual}
 Let $G^*$ be a dual-adapted rooted 3,4-map, let $\mathcal{L^*}$ be a dual 4-GS labeling, and let $\aA_o^*$ and $\aA_e^*$ be the corresponding augmented dual bipolar contractions of $G^*$. Let $h_o^*$ be an $\aA_o^*$-increasing function, and let $h_e^*$ be an $\aA_e^*$-increasing function. Placing each non-root vertex $v$ at coordinates $(h_o^*(v), h_e^*(v)) \in \br^2$, and drawing each non-root edge of $G^*$ using the $\cL^*$-rule (Definition~\ref{def:rule-draw-deges}), yields an SPO drawing of $G^*$.
\end{thm}

\fig{width=\linewidth}{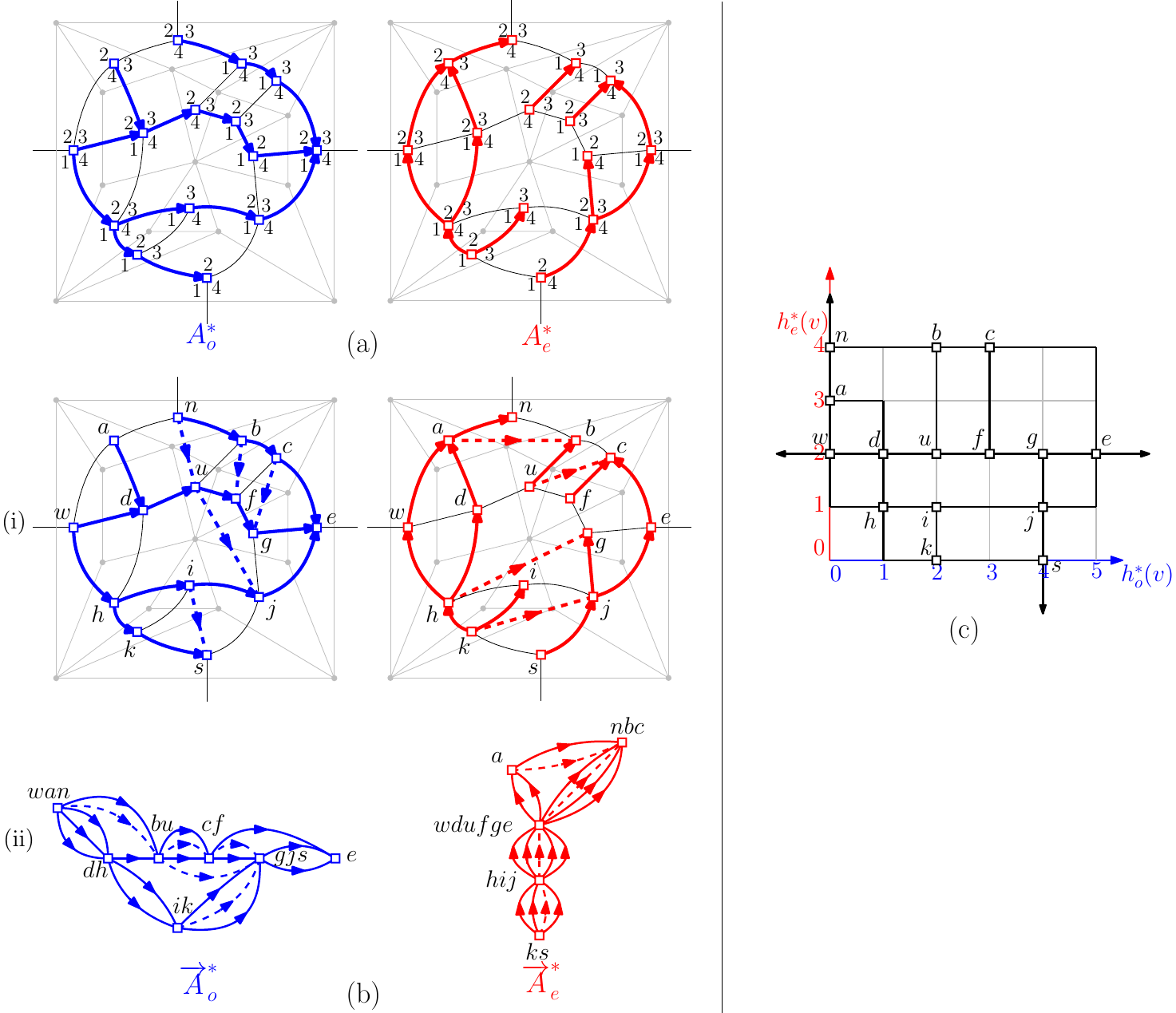}{(a) The dual acyclic orientations. (b) The augmented dual bipolar contractions. (c) The drawing obtained by choosing some tight increasing functions $h_o^*$ and $h_e^*$.}

An SPO drawing based on increasing functions is represented in Figure~\ref{fig:increasing_dual_big}. One can compare this drawing to the one in 
Figure~\ref{fig:face_counting_dual_big} obtained by the face-counting algorithm on the same map, and observe that the increasing-function drawing is slightly more compact.

Theorem~\ref{thm:increasing_dual} about the increasing function algorithm actually implies Theorem~\ref{thm:face_counting_dual} about the face-counting algorithm. This is because face-counting coordinates are increasing functions as stated below.

\begin{lem}\label{lem:dual_face_cnt_increasing}
Let $G^*$ be a dual-adapted rooted 3,4-map, and let $\mathcal{L^*}$ be a dual 4-GS labeling such that no edge is fully-colored. The face-counting coordinates $(\ell_e^*(v), r_o^*(v))$ defined in the previous subsection are such that the function $\ell_e^*$ is $\aA_o^*$-increasing and the function $r_o^*$ is $\aA_e^*$-increasing.
\end{lem}

\begin{proof} 
 We prove the statement for the function $\ell_e^*$ (the situation for $r_o^*$ is similar). We need to show three facts: 
 \begin{compactitem}
 \item[(a)] If $\{u,v\} \in \ovAos$, then $\ell_e^*(u) = \ell_e^*(v)$.
 \item[(b)] If $(u,v) \in A_o^*$, then $\ell_e^*(u) < \ell_e^*(v)$.
 \item[(c)] If $(u,v)$ is a diagonal arc in $\aA_o^*$, then $\ell_e^*(u) < \ell_e^*(v)$.
 \end{compactitem}

We prove (a) first. 
 One can easily check that any edge in $\ovAos$ bears a pair of opposite arcs of colors 2 and 4.
 Hence $P_e^*(u) = P_e^*(v)$ and $\ell_e^*(u) = \ell_e^*(v)$. 

 For (b), suppose $(u,v) \in A_o^*$. Given that $\mathcal{L^*}$ has no edge with incident corners labeled $\{1,1,3,3\}$, it is easy to see that the arc $(u,v)$ does not have color 4, and the opposite arc $(v,u)$ does not have color 2.
 Hence $P_e^*(u) \neq P_e^*(v)$. 
 Moreover, $P_e^*(u)$ is weakly on the left of $P_e^*(v)$. The situation is represented in Figure~\ref{fig:dual_edge}(a).
 In particular, $P_e^*(u)$ and $P_e^*(v)$ cannot cross each other: any common vertex of $P_e^*(u)$ and $P_e^*(v)$ is either a common vertex of $P_4^*(u)$ and $P_4^*(v)$ or a common vertex of $P_2^*(u)$ and $P_2^*(v)$. 
 But $P_2^*(u)$ and $P_2^*(v)$ (and similarly, $P_4^*(u)$ and $P_4^*(v)$) merge once they have  a common vertex. 
  This implies $\ell_e^*(u) < \ell_e^*(v)$.

\fig{width=.6\linewidth}{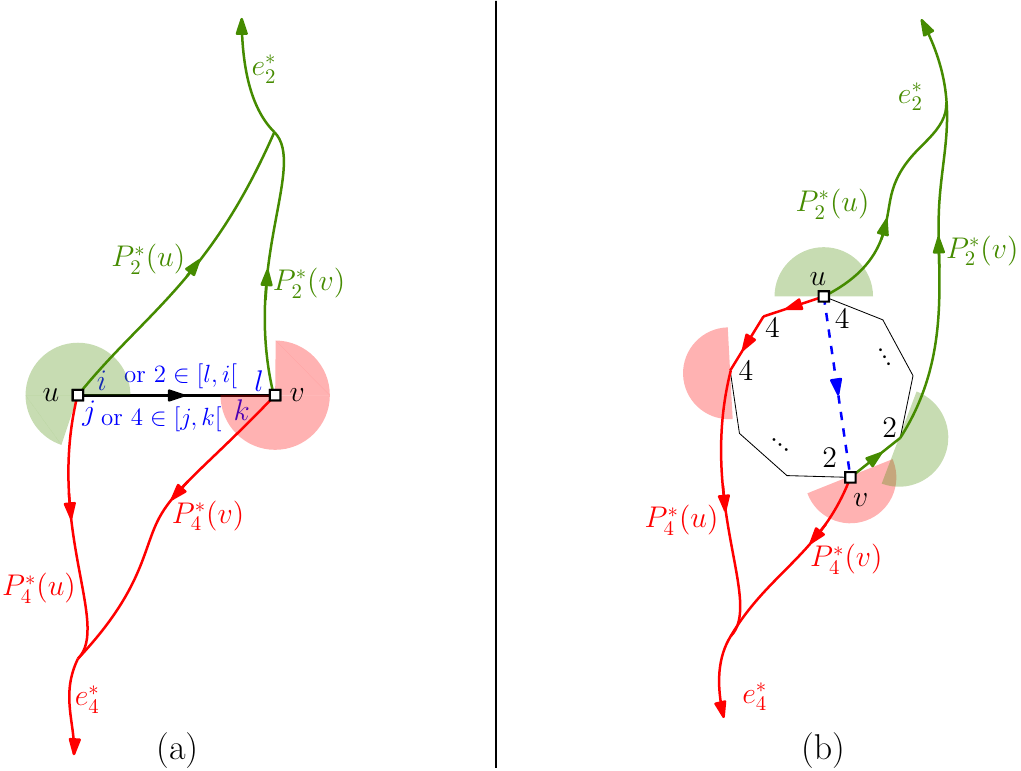}{(a) An arc $a=(u,v) \in A_o^*$. The range for position of the first arc of $P_2^*(u)$ is indicated by the green sector, and the range for position of the first arc of $P_4^*(v)$ is indicated by the red sector. (b) A non-root face, on which $u$ is incident to the last corner with label 4, and $v$ is incident to the last corner with label 2. The range for position of the first arc of $P_2^*(u)$ is indicated by the green sector, and the range for position of the first arc of $P_4^*(v)$ is indicated by the red sector.}

 For (c), recall that in $\aA_o^*$, a diagonal arc inside a non-root face is from the last corner with label 4 to the last corner with label 2, if both corners exist. Let us refer to the two vertices incident to these two corners as $u$ and $v$. We claim that $P_e^*(u) \neq P_e^*(v)$ and $P_e^*(u)$ is weakly on the left of $P_e^*(v)$. Indeed, $P_2^*(u)$ leaves the face at $u$, while $P_4^*(u)$ follows the boundary of the face until the first corner with label 4; $P_4^*(v)$ leaves the face at $v$, while $P_2^*(v)$ follows the boundary of the face until the first corner with label 2. This situation is presented in Figure~\ref{fig:dual_edge}(b). 
 This implies $\ell_e^*(u) < \ell_e^*(v)$, which completes the proof that $\ell_e^*$ is $\aA_o^*$-increasing.
\end{proof}

Next we show that the orientations $\aA_o^*$ and $\aA_e^*$ are acyclic (hence admit increasing functions).



\begin{lem}\label{lem:dual_bipolar}
The maps $\wt{A}_o^*$, $\wt{A}_e^*$, $\aA_o^*$ and $\aA_e^*$ are bipolar orientations.
\end{lem}

\begin{proof} We will prove the property for $\wt{A}_o^*$ and $\aA_o^*$ (the proof for $\wt{A}_e^*$ and $\aA_e^*$ is similar). 
In order to prove that $\aA_o^*$ is acyclic, we first show that the face-counting function $\ell_e$ is ``almost'' $\aA_o^*$-increasing. Precisely, by reasoning as in the proof of Lemma~\ref{lem:dual_face_cnt_increasing}, the following properties are easy to show:
 \begin{compactitem}
 \item[(a')] If $\{u,v\} \in \ovAos \cup \{\text{edges with labels \{1,1,3,3\}}\}$, then $\ell_e^*(u) = \ell_e^*(v)$.
 \item[(b')] If $(u,v) \in A_o^* \backslash \{\text{edges with labels \{1,1,3,3\}}\}$, then $\ell_e^*(u) < \ell_e^*(v)$.
 \item[(c)] If $(u,v)$ is a diagonal arc in $\aA_o^*$, then $\ell_e^*(u) < \ell_e^*(v)$.
 \end{compactitem}
Supposing for contradiction that $\aA_o^*$ has a directed cycle $C$, the above properties imply that $\ell_e^*$ is constant on $C$, and $C$ consists solely of fully-colored edges with labels $\{1,1,3,3\}$. This implies that $C$ is contained in the tree $W_3^*$ (and in fact every tree of the dual 4-GS wood), which is a contradiction. Hence $\aA_o^*$ is acyclic. 
Since $\aA_o^*$ is obtained from $\wt{A}_o^*$ by adding some edges and identifying some vertices, $\wt{A}_o^*$ must also be acyclic.

It remains to show that $\wt{A}_o^*$ and $\aA_o^*$ have a unique source and a unique sink. For this it suffices to show that for any non-root vertex $v$ of $G^*$, the vertex $\wt{v}$ of $\wt{A}_o^*$ corresponding to $v$ is incident to an ingoing and an outgoing arc in $\wt{A}_o^*$ (hence also in $\aA_o^*$). Observe that for every non-root vertex $v$ of $G^*$, the directed path $P_3^*(v)\cup P_1^*(v)^-$ is contained in $A_o^*\cup \ovAos$, which implies that the vertex $\wt{v}$ of $\wt{A}_o^*$ is on a directed path (hence is incident to an ingoing and an outgoing arc in $\wt{A}_o^*$). This completes the proof.
\end{proof}

We conclude this subsection by discussing the relationship between $\wt{A}_o^*$ and $\wt{A}_e^*$ and the primal bipolar orientations $B_o$ and $B_e$.

\fig{width=\linewidth}{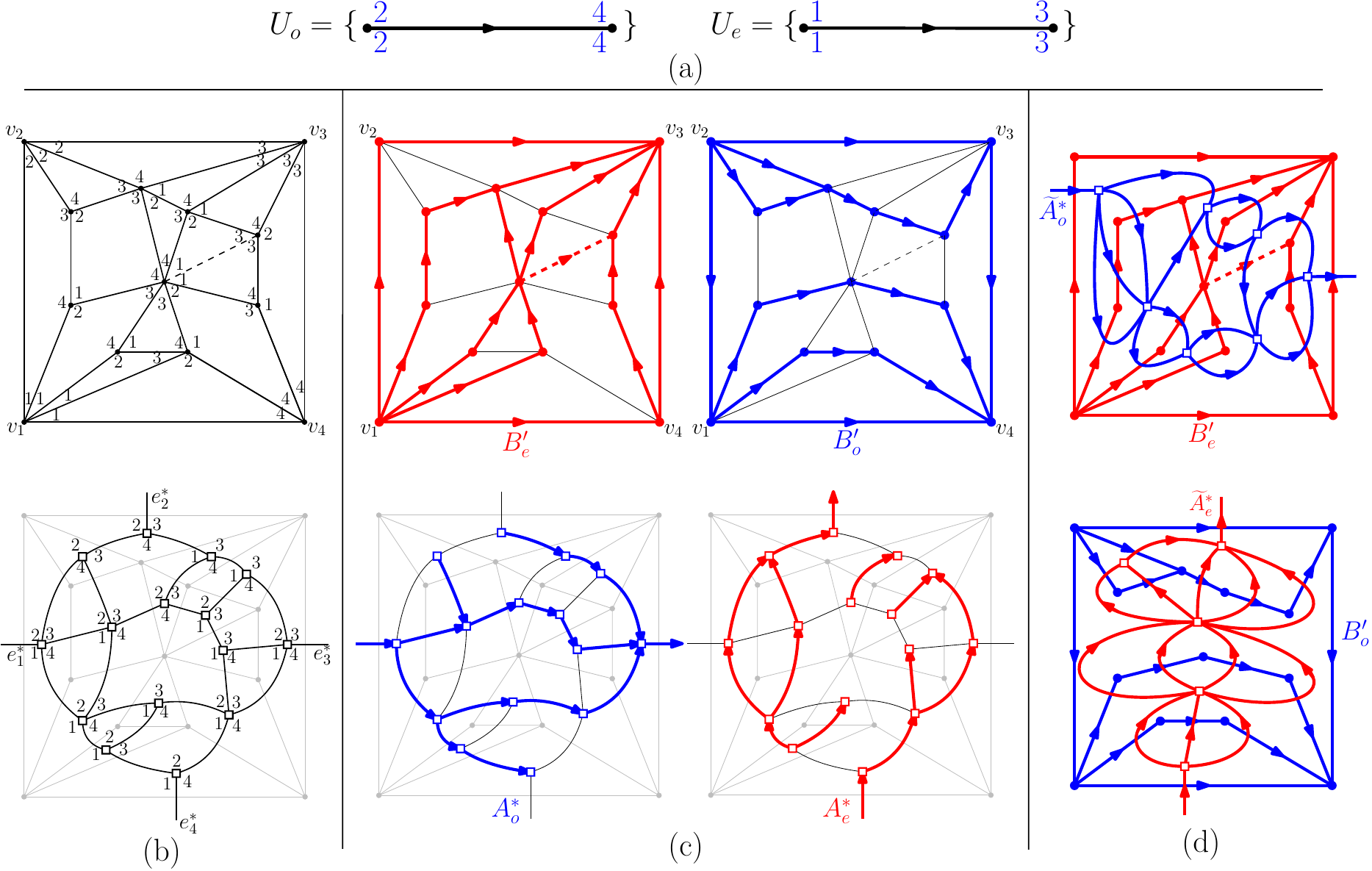}{(a) The inner arcs in $U_o$ and $U_e$. (b) A 4-GS labeling $\mathcal{L}$ (top), where the uncolored edge is presented as a dashed line, and its dual labeling $\mathcal{L}^*$ (bottom). (c) The augmented bipolar orientations $B_e'$ and $B_o'$ (top) and the dual acyclic orientations $A_o^*$ and $A_e^*$ (bottom). (c) Duality relations between $\wt{A}_o^*$ and $B_e'$ (top) and between $(\wt{A}_e^*)^-$ and $B_o'$ (bottom).}

\begin{remark}\label{rem:trans_bipolar}
Given $G$ and $\cL$ as above, we define the  \emph{odd} and \emph{even augmented bipolar orientations} as the following oriented submaps of $G$: 
$$B_o' = B_o \cup U_o\textrm{ and }B_e' = B_e \cup U_e$$
where 
\begin{align*}
 U_o &= \{a \text{ inner arc} \mid (\text{left-init}(a), \text{right-init}(a),\text{right-term}(a),\text{left-term}(a)) = (2,2,4,4) \}\\
 U_e &= \{a \text{ inner arc} \mid (\text{left-init}(a), \text{right-init}(a),\text{right-term}(a),\text{left-term}(a)) = (1,1,3,3)\}
\end{align*}
The definition of $U_o$ and $U_e$ is illustrated in Figure~\ref{fig:dual_primal_pair}(a). Note that when the primal map $G$ is an irreducible triangulation, the augmented bipolar orientations $B_o'$ and $B_e'$ associated to a 4-GS labeling $\cL$ of $G$ are exactly the blue and red submaps $G_b$ and $G_r$ of the transversal structure corresponding to $\cL$.

It easily follows from the definition that the bipolar orientation $\wt{A}_o^*$ is the \emph{dual} of $B_e'$, while the  bipolar orientation $\wt{A}_e^*$ is the \emph{dual} of $(B_o')^-$. This is represented in Figure~\ref{fig:dual_primal_pair}. This can be compared with the duality between $\wt{A}_o$ (resp. $\wt{A}_e$) and $(B_e^*)^-$ (resp. $B_o^*$) represented in Figure~\ref{fig:dual-orientations}. These two duality relations are summarized in Figure \ref{fig:two-dual-bipolars}. 
\end{remark}
\fig{width=\linewidth}{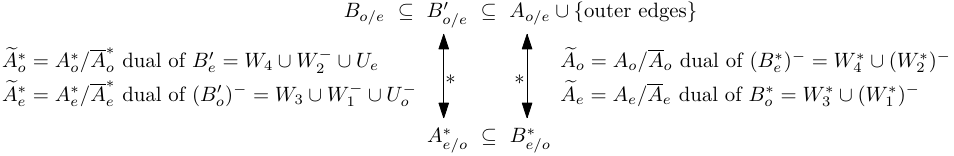}{Duality relations for the bipolar orientations underlying the increasing function algorithms. One uses the  orientations $\wt{A}_o$ and $\wt{A}_e$ (which are the dual of $(B_e^*)^-$ and $B_o^*$) for drawing 3,4-angulations, and the orientations  $\wt{A}^*_o$ and $\wt{A}^*_e$ (which are the dual of $B_e'$ and $(B_o')^-$) for drawing 3,4-valent maps.}

\subsection{Proof of planarity}\label{sec:planarity_ortho}
In this subsection we prove Theorem~\ref{thm:increasing_dual}.

Let $G^*$ be a dual-adapted rooted 3,4-map, and let $\mathcal{L}^*$ be a dual 4-GS labeling on $G^*$. 
Let $G'$ be the submap $G^*\backslash \{v_\infty\}$, and let $V'$ be the set of non-root vertices of $G$. 
Let $h_o^*: V' \to \br$ be an $\aA_o^*$-increasing function, and let $h_e^*: V' \to \br$ be an $\aA_e^*$-increasing function.

We want to show that the vertex placement $\al: v \mapsto (h_o^*(v), h_e^*(v))$ induces an SPO drawing of~$G^*$. 
The proof again relies on the planarity criterion established in Section~\ref{sec:proof-planar-primal} (see Lemmas~\ref{lem:planar-criteria} and~\ref{lem:outer_right_visible}). 
Recall from Section~\ref{sec:proof-planar-primal} that a drawing of $G^*$ is called \emph{valid} for a face $f$ of $G^*$ if the curve corresponding to the clockwise contour of $f$ is a simple curve in the plane whose interior is on its right.

\begin{lemma} \label{lem:innerfacesdual-valid}
The vertex placement $\al$ is such that each non-root edge can be drawn according to the $\cL^*$-rule (Definition~\ref{def:rule-draw-deges}). It produces an orthogonal drawing of $G'$ which is valid for every inner face of $G'$.
\end{lemma}

\begin{proof}
Let $f$ be an inner face of $G'$. 
The labels of the corners inside $f$ in clockwise order form an interval of 1, an interval of 2, an interval of 3 and an interval of 4 (where at most two of these intervals can be empty). 
We claim that the face $f$ will be embedded as represented in Figure~\ref{fig:dual_face_valid} (the shape of the embedding depends on which labels appears inside $f$).

\fig{width=\linewidth}{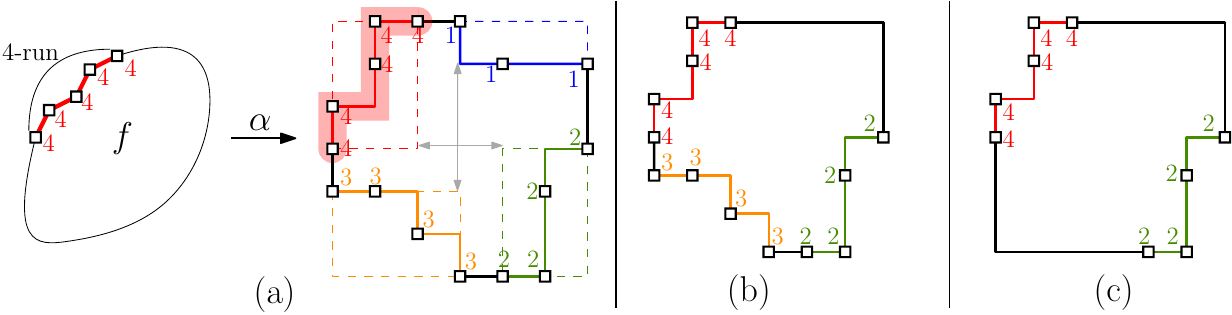}{(a) The embedding of an inner face $f$ with no missing labels. The different runs are indicated by different colors. (b) The embedding of an inner face $f$ missing the label 1. (c) The embedding of an inner face $f$ missing labels 1 and 3.}

For $i\in\four$, we call \emph{$i$-run} the sequence of arcs having $f$ on their right, with the two labels on their right equal to $i$.
We consider the drawing of the arcs in the $i$-run. Let us consider $i=4$ for concretness. Let $a=(u,v)$ be an arc in the 4-run of $f$. It is easy to see that the pair of labels $(j,k)=(\linit(a),\lterm(a))$ on the left of $a$ can be either $(3,1)$, $(3,2)$, $(2,1)$ or $(2,2)$. These situations are represented in Figure~\ref{fig:dual_edge2}(a).
\begin{compactitem}
\item If $(j,k)=(3,1)$, then $a$ has color 3 while $-a$ has color 4, so $a\in A_o^* \cap A_e^*$. In this case the $\cL^*$-rule for drawing $a$ can be applied and $a$ is drawn as a bent edge going ``right-up".
\item If $(j,k)=(3,2)$, then $a$ has color 3 while $-a$ has colors 1 and 4, so $a\in A_o^* \cap \ovAes$. In this case the $\cL^*$-rule for drawing $a$ can be applied and $a$ is drawn as an horizontal edge going ``right".
\item If $(j,k)=(2,1)$ (resp. $(2,2)$), then $a$ has colors 2 and 3 while $-a$ has color 4 (resp. colors 1 and 4), so $a\in \ov{A_o^*} \cap A_e^*$. In this case the $\cL^*$-rule for drawing $a$ can be applied and $a$ is drawn as a vertical edge going ``up".
\end{compactitem}
In conclusion, the 4-run of $f$ goes (weakly) right-up. Symmetrically, the 1-run goes down-right, the 2-run goes left-down, and the 3-run goes up-left. 
\fig{width=\linewidth}{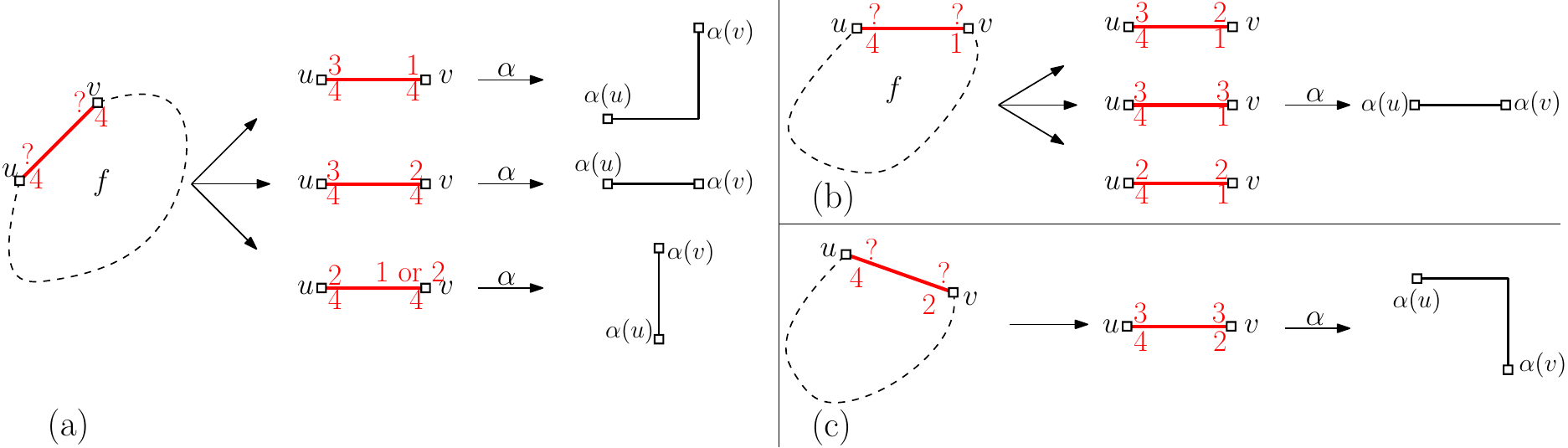}{(a) Embedding of edges with labels $\{4,4\}$ inside $f$. (b) Embedding of edges with labels $\{4,1\}$ inside $f$. (c) Embedding of edges with labels $\{4,2\}$ inside $f$.}

Next we consider the edges between runs. 
Consider an arc $a$ having $f$ on its right with labels $(\rinit(a),\rterm(a))=(4,1)$. Then, the pair of labels $(j,k)=(\linit(a),\lterm(a))$ on the left of $a$ can be either $(3,2)$, $(3,3)$ or $(2,2)$. In every case the $\cL^*$-rule for drawing $a$ can be applied and $a$ is drawn as an horizontal edge going ``right". This is represented in Figure~\ref{fig:dual_edge2}(b).
Symmetrically, 
\begin{compactitem}
\item if $a$ is an arc having $f$ on its right with labels $(\rinit(a),\rterm(a))=(1,2)$, then it is drawn down,
\item if $a$ is an arc having $f$ on its right with labels $(\rinit(a),\rterm(a))=(2,3)$, then it is drawn left,
\item if $a$ is an arc having $f$ on its right with labels $(\rinit(a),\rterm(a))=(3,4)$, then it is drawn up.
\end{compactitem}
Finally, an arc $a$ having $f$ on its right with labels $(\rinit(a),\rterm(a))=(4,2)$ will be drawn as a bent edge going ``right then down". This is represented in Figure~\ref{fig:dual_edge2}(c). Symmetric results olds for an arc with label $(\rinit(a),\rterm(a))=(1,3)$ or $(2,4)$ or $(3,1)$. 

To summarize, every non-root edged of $G^*$ can be drawn following the $\cL^*$-rule which produces an orthogonal drawing of $G'$. Moreover, for each inner face $f$, the drawing of the clockwise contour of $f$ is as follows: 
the 1-run drawn going down-right, a down segment, the 2-run drawn left-down, a left segment, the 3-run drawn up-left, an up segment, the 4-run going right-up, a right segment. 

It is clear from the properties stated above that consecutive runs of $f$ (that is, a $i$-run and a $(i+1)$-run) cannot intersect one another. Hence the drawing is valid for the face $f$ unless opposite runs (that is, a $i$-run and a $(i+2)$-run) intersect each other. 
Now, the function $h_o^*$ is increasing along odd diagonal arcs inside $f$ (which goes from the end of the 4-run to the end of the 2-run), hence the 4-run is strictly to the left of the 2-run; and the function $h_e^*$ is increasing along even diagonal arcs inside $f$ (which goes from the end of the 3-run to the end of the 1-run), hence the 3-run is strictly below the 1-run. In conclusion, the drawing of $f$ is as illustrated in Figure~\ref{fig:dual_face_valid}, and it is valid.
\end{proof}

In order to finish the proof of Theorem~\ref{thm:increasing_dual}, we need to show that the orthogonal drawing of $G'$ is valid for the outer face.

Consider the vertices on the outer boundary of $G'$ that are adjacent to the root vertex $v_\infty$ in $G^*$. Let us denote the vertex incident to the root-edge $e_1^*$ (resp. $e_2^*$, $e_3^*$, $e_4^*$) by $v_1^*$ (resp. $v_2^*$, $v_3^*$, $v_4^*$).  The vertices $v_1^*,v_2^*,v_3^*,v_4^*$ divide the clockwise outer contour of $G'$ into four parts. We refer to the part from $v_i^*$ to $v_{i+1}^*$ as the \emph{outer $(i+1)$-run}, as all the labels on the outer face of this part are $(i+1)$. Some of the outer runs can be empty since we could have $v_i^*=v_{i+1}^*$; however we cannot have $v_i^*=v_{i+2}^*$ unless $G'$ consists of a single vertex.

It is not hard to see that, for any pair of increasing functions $h_o^*$ and $h_e^*$, the four runs of the outer contour satisfy some specific inequalities about $h_o^*$ and $h_e^*$, which translate into the following orientation properties:
\begin{compactitem}
 \item The outer 2-run is directed up-right.
 \item The outer 3-run is directed right-down.
 \item The outer 4-run is directed down-left.
 \item The outer 1-run is directed left-up.
\end{compactitem}

In other words, the shapes of the four outer runs are similar to the shapes of the four runs of a generic inner face (with index of the run offset by 2). See Figure~\ref{fig:cut}(a) (left) for an example.
This alone however does not guarantee that the drawing of the outer face of $G'$ is valid. 
The difference is that in the case of inner faces, the condition on diagonal arcs ensure that opposite runs do not intersect. 
For the outer face however, not only are such diagonal arc requirements non-existent, but also the outer contour may have cut vertices in the first place.

Fortunately, the shapes of the four runs imply that all cut vertices, if they exist, must belong to one of the two pairs of opposite outer runs, say the outer 3-run and the outer 1-run. The cut vertices divide $G'$ into a chain of 2-connected subgraphs, that we call \emph{2-blocks} of $G$ (consecutive 2-blocks intersect at a cut vertex). It is clear from the orientation property of the outer-face of $G'$ that the drawing of $G'$ is planar if and only if the drawing of each 2-block is planar. An example is illustrated in Figure~\ref{fig:cut}(a) (right).
%

\fig{width=\linewidth}{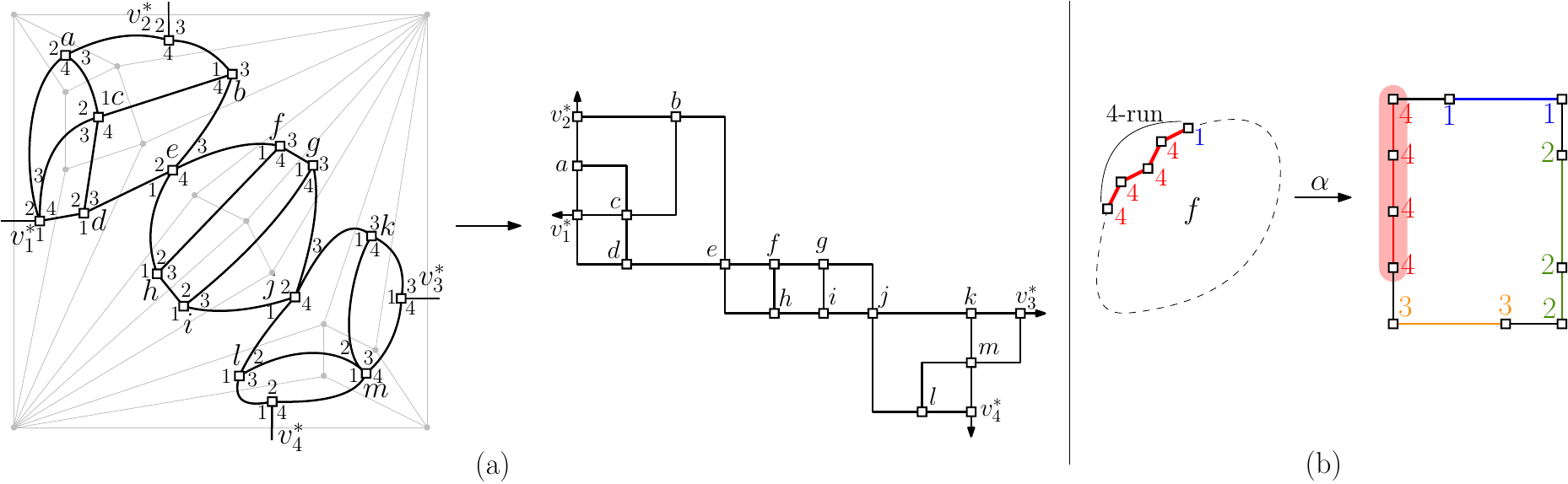}{(a) Left: A dual 4-GS labeling $\mathcal{L}^*$, with two cut vertices~$e$ and~$j$. They divide $G'$ into three 2-blocks. Right: The drawing obtained by choosing some increasing functions. (b) The embedding of a face by He's algorithm.}

Let $H$ be a 2-block of $G$. Recall that the drawing of each inner face of $H$ is valid by Lemma~\ref{lem:innerfacesdual-valid}. Hence, if we can show that the drawing of the outer face of $H$ is valid, then we can invoke Lemma~\ref{lem:planar-criteria} to conclude that the drawing of $H$ is planar. Let us suppose for contradiction that the drawing of the outer face of $H$ is not valid. Given the orientation property of the runs, this implies that the boundary of $H$ must have some outer-right-visible point (recall the definition from Lemma \ref{lem:outer_right_visible} in Section~\ref{sec:proof-planar-primal}). This is impossible according to Lemma~\ref{lem:outer_right_visible}. Therefore the outer face of $H$ is valid, and by Lemma~\ref{lem:planar-criteria}, the drawing of $H$ is planar.

Since the drawing of every 2-block of $G$ is planar, the drawing of the whole graph $G'$ is planar. 
It is also clear that the placements of $v_i^*$'s are extremal in their respective directions, and hence the four dangling arrows can be added to the outer face, completing a valid SPO drawing of $G^*$. 
This completes the proof of Theorem~\ref{thm:increasing_dual}.

\subsection{Comparison to He's algorithm}\label{sec:comparison-He}

We now compare the increasing function algorithm of Theorem~\ref{thm:increasing_dual} to He's algorithm from~\cite{He93:reg-edge-labeling}. He's algorithm is based on a transversal structure on an irreducible triangulation $G$ and a pair of increasing functions on the faces of the red and blue submaps $G_r$ and $G_b$. 
As explained in Section~\ref{sec:4GS1}, each transversal structure also corresponds to a 4-GS labeling, and Remark~\ref{rem:trans_bipolar} points out that $G_r$ and $G_b$ are equal to the augmented bipolar orientations $B_e'$ and $B_o'$. 
Given the duality relations between the dual bipolar contractions $\wt{A}_o^*, \wt{A}_e^*$ and $B_e', B_o'$, it is easy to see that the increasing functions on the faces of $G_r$ and $G_b$, which are used in He's algorithm, naturally corresponds to a pair of increasing functions on the vertices of $\wt{A}_o^*$ and $\wt{A}_e^*$.
We claim that such functions are also increasing functions on $\aA_o^*$ and $\aA_e^*$, because they automatically increase along the diagonal arcs.
To see this, we can investigate the shape of faces in the drawing.
First, note that an edge connecting two vertices of degree 3 cannot have exactly two colors with different parity, and hence every edge of $G^*$ is drawn with no bend. 
Next, inside each inner face of $G^*$, all four labels must appear, and with an analysis as in Section~\ref{sec:planarity_ortho}, one can observe that the shape of the four runs around an inner face of $G^*$ degenerates to something simple. Namely, 
\begin{compactitem}
 \item the 1-run is drawn as a sequence of edges directed horizontally to the right, 
 \item the 2-run is drawn as a sequence of edges directed vertically downward, 
 \item the 3-run is drawn as a sequence of edges directed horizontally to the left, 
 \item the 4-run is drawn as a sequence of edges directed vertically upward.
\end{compactitem}
Consequently, each inner face $f$ of $G^*$ is drawn as a rectangle, with the four corners of that rectangle corresponding to the last vertices incident to a corner of label 1,2,3,4 respectively in clockwise order around $f$. See Figure~\ref{fig:cut}(b) for an illustration. Hence the functions used in He's algorithm already ensures strict increments along diagonal arcs.
In short, the algorithm of~\cite{He93:reg-edge-labeling} is the same as the increasing function algorithm from Theorem~\ref{thm:increasing_dual} restricted to 
the case where all non-rooted vertices have degree $3$. 


\subsection{Properties of the SPO drawings}\label{sec:prop-edges-dual}
In this subsection we discuss some additional properties of the SPO drawings obtained by the increasing function algorithm. In previous subsections, we have investigated the shapes of the embedded faces, and we also remarked that an edge connecting two degree 3 vertices must be straight. Here we investigate some further orientation properties of embedded edges. These properties will provide some intuitions for the bend-minimization method discussed in Section~\ref{sec:bend}. 

There is a particular type of edges we are interested in. We say that a non-root arc $a = (u,v)$ is \emph{special} if we have 
$$(\text{right-init}(a),\text{right-term}(a),\text{left-term}(a), \text{left-init}(a)) = (k,k,k+2,k+3)$$ 
for some $k\in\four$. If one arc of an edge $e$ is special, then we also call $e$ \emph{special}. See Figure~\ref{fig:special}(a) for an example.

\fig{width=0.7\linewidth}{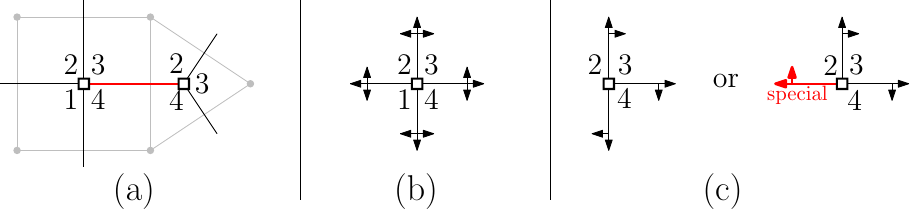}{(a) A special edge with labels $\{4,4,2,3\}$. (b) Orientation properties around vertices of degree 4. (c) Orientation properties around vertices of degree 3 with incident labels $\{2,3,4\}$ depending on whether it is incident to a special edge or not (the red arc is part of a special edge).}

Special edges satisfy some unusual properties. 
For one thing, from Condition (L3$^*$) of the definition of dual 4-GS labelings, one can easily see that if an arc $a = (u,v)$ is special, then $u$ must be a vertex of degree 4 and $v$ of degree 3. 
Furthermore, using the terminology of $i$-runs introduced in Subsection \ref{sec:planarity_ortho}, one can easily verify that each run around a face contains at most one special edge, and if there is one then it has to be the last edge on that run (in clockwise order).


We now state the orientation properties of edges in our drawings.


\begin{prop}\label{prop:orientation}
 An SPO drawing obtained by the increasing function algorithm satisfies the following properties: 
\begin{compactitem}
 \item An arc whose initial vertex has degree 3 never turns left: it is either drawn as a straight segment, or it turns right.
 \item An arc with a unique color $k$ is drawn so that it starts by going in direction $k$.
 \item An arc with  colors $k$ and $k+1$ is drawn so that it starts by going in direction $k$, unless it is part of a special edge, in which case it starts by going in direction $k+1$.
\end{compactitem}
\end{prop}

Observe that the edge properties listed in Proposition \ref{prop:orientation} imply that the situation around vertices is as illustrated in Figure \ref{fig:special}(b-c).

\begin{proof}
 Let us prove the first property.
 First, by $\cL^*$-rule for drawing edges, an arc with more than one color will be drawn as a straight segment. Now suppose an arc $a = (u,v)$ has a unique color $k$, and $u$ has degree 3. For it to turn left, the opposite arc $-a$ needs to have a unique color $k+1$, but this violates Condition (L3$^*$) of dual 4-GS labelings. 

 The later two properties can be verified using a similar analysis as in previous subsections.
\end{proof}



\subsection{On allowing vertices of degree 2}\label{sec:digon}
In this subsection we explain how to draw graphs having some vertices of degree 2.

For an arbitrary graph with a vertex $v$ of degree 2, we call \emph{erasing} $v$ the operation of deleting $v$ and replacing the 2 incident edges $\{v,u_1\}$ and $\{v,u_2\}$ by a single edge $\{u_1,u_2\}$. 
See figure~\ref{fig:digon}(a). 

\fig{width=0.9\linewidth}{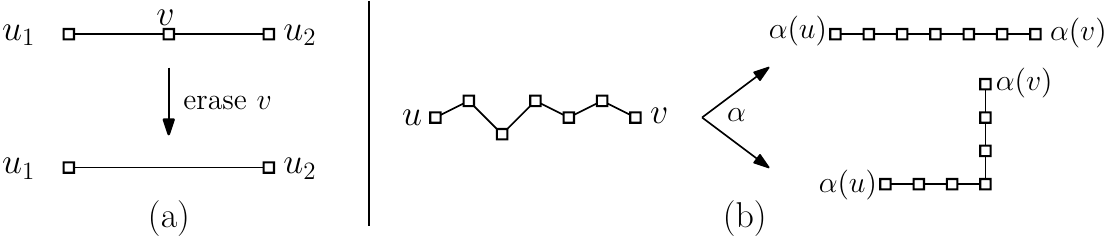}{(a) Erasing a vertex of degree 2. (b) Two possible embeddings of a chain with five vertices of degree 2.}

Now, let $G^\circ$ be a \plm whose vertex degrees are in $\{2,3,4\}$, with a designated root vertex $v_\infty$. For simplicity let us suppose $v_\infty$ is not incident to any vertex of degree 2. Let $G^*$ be the \plm obtained by erasing all the vertices of degree 2. Then every non-root edge of $G^*$ corresponds to a path of edges in $G^\circ$. 

Assuming that $G^*$ is a dual-adapted rooted 3,4-map, we claim that we can easily obtain a grid SPO drawing for $G^\circ$. This takes two steps: 
\begin{compactitem}
 \item[(1)] Compute a dual 4-GS labeling on $G^*$, and run the increasing function algorithm with an additional requirement: if a non-root edge $e$ of $G^*$ corresponds to a path in $G^\circ$ with $N_e$ edges, then $e$ needs to be drawn with total length at least $N_e$ (either as a straight segment or as a bent-edge).
 \item[(2)] Insert the vertices of degree 2 on the edges of $G^*$ they belong to. Due to the additional requirement in Step (1), they all can be placed at grid points.
\end{compactitem}

The additional requirement in Step (1) is easy to achieve: for each edge $e$ of $G^*$, record the length $N_e$ of the corresponding path in $G^\circ$. Then, when computing the increasing functions, make sure that the increments along $e$ in $\aA_o^*$ and $\aA_e^*$ sum to at least $N_e$. This ensures that the total length of $e$ in the drawing is at least $N_e$. See Figure~\ref{fig:digon}(b) for an illustration.




\subsection{Time complexity}
In this subsection we discuss the time complexity and the grid size that can be achieved by the increasing function orthogonal drawing algorithm. 



Let $G^*$ be a dual-adapted rooted 3,4-map. To perform the increasing function orthogonal drawing algorithm one needs to: 
\begin{compactitem}
  \item[(i)] compute a dual 4-GS labeling $\mathcal{L^*}$, and the corresponding orientations $\aA_o^*, \aA_e^*$, and
  \item[(ii)] compute some integer-valued $\aA_o^*$-increasing and $\aA_e^*$-increasing functions $h_o^*$ and $h_e^*$ that are \emph{tight} (recall the definition from Section~\ref{sec:time-complexity}).
\end{compactitem}

Step (i) can be completed in linear time since it was proved in~\cite{OB-EF-SL:Grand-Schnyder} that 4-GS labeling could be completed in linear time. Step (ii) can be completed using a similar procedure as in the primal case, which is elaborated in Section~\ref{sec:time-complexity}. We summarize these claims into the following statement. 

\begin{prop}
  The increasing function orthogonal drawing algorithm for drawing a dual-adapted rooted 3,4-map can be completed in linear time in the number of vertices.
\end{prop}

We mention that the time complexity remains linear when allowing for vertices of degree 2 using the method discussed in Section~\ref{sec:digon}.



\subsection{Bound on the grid size}

In this subsection we discuss the size of the grid needed for the orthogonal grid-drawing algorithms.

Let $G^*$ be a dual-adapted rooted 3,4-map, and let $\mathcal{L}^*$ be a dual 4-GS labeling. Denote by $m_o^* \times m_e^*$ the grid size resulting from the tight increasing functions associated to $\mathcal{L}^*$. Note that $m_o^*$ is the maximal length of directed paths in $\aA_o^*$, and $m_e^*$ is the maximal length of directed paths in $\aA_e^*$. 
In general, it seems difficult to estimate $m_o^*$ and $m_e^*$.
However, an obvious bound is $m_o^* \leq \vv(\aA_o^*)-1$, where $\vv(\aA_o^*)$ is the number of vertices of $\aA_o^*$, and similarly $m_e^* \leq \vv(\aA_e^*)-1$, where $\vv(\aA_e^*)$ is the number of vertices of $\aA_e^*$.

Now, recall from Section~\ref{sec:increasing_dual} the duality relations between $\wt{A}_o^*$ and $B_e'=B_e \cup U_e$, and between $(\wt{A}_e^*)^-$ and $B_o'=B_o \cup U_o$ (see Remark~\ref{rem:trans_bipolar}). 
This implies $m_o^* \leq \ff(B_e')-1$ and $m_e^* \leq \ff(B_o')-1$, where $\ff(B_e')$ and $\ff(B_o')$ are the number of inner faces of $B_e'$ and $B_o'$. Since $B_e'=B_e \cup U_e$ and $B_o'=B_o \cup U_o$ we get 
$$\ff(B_e')=\ff(B_e)+|U_e|\textrm{ and } \ff(B_o')=\ff(B_e)+|U_o|,$$
where $\ff(B_e)$ is the number of inner faces of $B_e$, and $\ff(B_o)$ is the number of inner faces of $B_o$. 

It was shown in Section~\ref{sec:grid-size} that $\ff(B_e)=\mathbf{v}+1-\mathbf{d}_e$, and $\ff(B_o)=\mathbf{v}+1-\mathbf{d}_o$, where $\mathbf{v}$ is the number of inner vertices of the primal map $G$, and $\mathbf{d}_e$ (resp. $\mathbf{d}_o$) is the number of inner edges of $G$ with colors $\{2,4\}$ (resp. colors $\{1,3\}$).  
Combining these results and expressing them in terms of parameters of $G^*$ gives the following bound for the grid size.


\begin{prop}\label{prop:grid_size_dual}
  Let $G^*$ be a dual-adapted rooted 3,4-map, and let $\mathcal{L}^*$ be a dual 4-GS labeling. 
  The increasing function orthogonal drawing algorithm produces a grid drawing of $G^*$ on a grid of size $m_o^* \times m_e^*$, 
  with 
  $$m_o^* \leq \ff^*-\mathbf{d}_o^*+\mathbf{u}_o^*,\textrm{ and }m_e^* \leq \ff^*-\mathbf{d}_e^*+\mathbf{u}_e^*,$$
  where $\ff^*$ (resp. $\mathbf{d}_o^*,\mathbf{d}_e^*,\mathbf{u}_o^*,\mathbf{u}_e^*$) is the number of non-root faces (resp. edges of color $\{1,3\}$, edges of color $\{2,4\}$, fully-colored edges with labels $\{1,1,3,3\}$, fully-colored edges with labels $\{2,2,4,4\}$) of $G^*$.
\end{prop}

In the case where $G$ is a quadrangulation and the 4-GS labeling $\cL$ is even, there are no edges of $G^*$ with exactly two colors of the same parity nor any fully-colored edges, hence $m_o^*$ and $m_e^*$ are both bounded by $\ff^*$. 
This bound matches the grid size of the Bernardi and Fusy's algorithm \cite{OB-EF:Schnyder} (before the compaction step).

In the case where $G$ is a triangulation, $\mathbf{d}_o^*+\mathbf{d}_e^*$ is equal to the number $e_2$ of edges of $G$ bearing exactly 2 colors, while $\mathbf{u}_o^*+\mathbf{u}_e^*$ is equal to the number $e_0$ of edges of $G$ bearing no colors, and the calculation in Section~\ref{sec:grid-size} gives $\mathbf{d}_o^*+\mathbf{d}_e^*= \mathbf{f}^*-1+\mathbf{u}_o^*+\mathbf{u}_e^*$. In this case, the bound on the half-perimeter of the drawing is simply $m_o^*+m_e^*\leq \mathbf{f}^*+1$, which matches the bound given by He's algorithm~\cite{He93:reg-edge-labeling}, a bound also achieved by the algorithm in~\cite{RahmanNN98} with a different approach.


\begin{remark}\label{rk:biedl_bound}
Our drawings have the property that every grid line is occupied by at least one vertex. In that case, letting $\mathbf{v}^*$ be the number of non-root vertices, and $s_h$ (resp. $s_v$) be the number of straight horizontal (resp. vertical) edges, 
An easy argument (see~\cite{biedl1996optimal}) then ensures that the grid-width is at most $\mathbf{v}^*-1-s_v$ and the grid-height is at most $\mathbf{v}^*-1-s_h$.  This matches the bounds for $m_o^*$ and $m_e^*$ from Proposition~\ref{prop:grid_size_dual}, since $\mathbf{v}^*-s_v=\vv(\aA_o^*)$ (contracted edges in $\aA_o^*$ are those counted by $s_v$), and similarly $\mathbf{v}^*-s_h=\vv(\aA_e^*)$. 	
\end{remark}

\begin{remark}
There is another well-known approach, due to Tamassia~\cite[Theorem 1]{tamassia1987embedding}, to obtain a grid realization of a given (combinatorial) plane  orthogonal drawing. It consists in refining the orthogonal drawing into a (combinatorial) rectangular drawing by  suitable insertions of vertices and edges, and then applying He's algorithm~\cite{He93:reg-edge-labeling} to obtain a grid realization of the rectangular drawing, which yields a grid realization of the original orthogonal drawing. Since every grid-line is occupied by at least one vertex, the bounds of Proposition~\ref{prop:grid_size_dual} also hold with this strategy.   
Our approach shows that one can directly apply an increasing function strategy in order to obtain a grid realization, without the preliminary refinement step.    	
\end{remark}

We now state some further bounds on the number of straight edges and grid size in terms of parameters of the planar map solely.

\begin{proposition}\label{prop:nb-bend-ortho}
Let $G^*$ be a dual-adapted rooted 3,4-map, let $\mathbf{e}^*$ be its number of non-root edges, and let $\mathbf{v}_3^*$ (resp. $\mathbf{v}_4^*$)  be its number of vertices of degree $3$ (resp. non-root vertices of degree $4$). 

For any dual 4-GS labeling $\cL^*$, the increasing function orthogonal drawing algorithm produces a grid drawing of $G^*$ whose number of straight edges is at least $\mathbf{v}_3^*/2$ and number of bent edges is at most $4\mathbf{v}_4^*$. Consequently, the semi-perimeter of the grid is at most $\mathrm{min}(\mathbf{e}^*, \mathbf{e}^*+2\mathbf{v}_4^*-\mathbf{v}_3^*+2)$. 

Moreover, upon choosing $\cL^*$ appropriately (linear time), the number of bent edges is at most $3\mathbf{v}_4^*$ and the semi-perimeter of the grid is at most $\mathrm{min}(\mathbf{e}^*, \mathbf{e}^*+\mathbf{v}_4^*-\mathbf{v}_3^*+2)$. 
\end{proposition}

\begin{proof}
The second and third items in Proposition~\ref{prop:orientation} imply  that, for every vertex of degree~$3$, the arc with two colors is on a straight edge. Hence the number of straight edges is at least $\mathbf{v}_3^*/2$. On the other hand, the first item in Proposition~\ref{prop:orientation} ensures that every bent edge turns left from a non-root vertex of degree $4$. Hence, the number of bent edges is at most $4\mathbf{v}_4^*$. 

Then it follows from Biedl's bound (see Remark~\ref{rk:biedl_bound})  that the semi-perimeter is at most 
$$\mathrm{min}(2\mathbf{v}^*-2-\mathbf{v}_3^*/2,2\mathbf{v}^*-2-\mathbf{e}^*+4\mathbf{v}_4^*)=\mathrm{min}(\mathbf{e}^*,\mathbf{e}^*+2\mathbf{v}_4^*-\mathbf{v}_3^*+2),$$
where $\mathbf{v}^* = \mathbf{v}_4^* + \mathbf{v}_3^*$, and we used $\mathbf{e}^*=2\mathbf{v}_4^*+\tfrac{3}{2}\mathbf{v}_3^*-2$.

Finally, we observe that for any vertex of degree $4$ with contribution $4$ to the number of bent edges, one can reduce this contribution to 0 by decreasing by $1$ (mod $4$) the labels at the 4 corners around that vertex. These operations can be performed simultaneously at all such  vertices, and yield a new  dual 4-GS labeling for which each vertex of degree $4$ has contribution at most $3$ to the number of bends. Thus, the number of bent edges after these operations is at most $3\mathbf{v}_4^*$, and the grid semi-perimeter is at most  $\mathrm{min}(\mathbf{e}^*, \mathbf{e}^*+\mathbf{v}_4^*-\mathbf{v}_3^*+2)$. 
\end{proof}

\begin{remark}
The bound $\mathbf{e}^*$ on the semi-perimeter is also satisfied by the algorithm of Biedl and Kant~\cite{biedl1998better}, which applies to the (larger) family of 2-connected planar maps with vertex degrees at most 4.
 However, edges can have up to 2 bends with that algorithm. 
\end{remark}

\begin{remark}
The bound $3\mathbf{v}_4^*$ for the number of bends of the drawing of $G^*$ (in Proposition \ref{prop:nb-bend-ortho}) shows that our algorithm produce a drawing at distance at most $3\mathbf{v}_4^*$ from the optimal bend-number over all plane orthogonal drawings of $G^*$. Having distance  at most $3\mathbf{v}_4^*$ to the optimal bend-number also appears in the orthogonal drawing algorithm given in~\cite{BhatiaLK18} (building on~\cite{RahmanNN99}) for 3-connected planar maps with vertex degrees at most 4
(a family which intersects ours but without inclusion of one into the other), where however the optimal bend number can itself be larger than $3\mathbf{v}_4^*$. 
\end{remark}

\subsection{Reduction of the grid size}\label{sec:reduce-grid-dual}
In this subsection we discuss two methods to reduce the grid size of the increasing function orthogonal drawing algorithm. Similarly as for the primal case, these improvements will be obtained by relaxing some constrains for the increasing functions $h_o^*,h_e^*$.

\begin{remark} Throughout this section we work with a fixed dual 4-GS labeling. We note however that in \cite{He93:reg-edge-labeling}, the author discussed some methods to reduce the grid size by changing the transversal structure. We leave as an open question whether these methods can be generalized.
\end{remark}


\bigskip


\ni \textbf{Relaxing the equalities for edges in $B_o^*\setminus A_o^*$ and $B_e^*\setminus A_e^*$:} 
The first optimization we present is about relaxing some equalities for the increasing functions $h_o^*,h_e^*$.
Recall that the increasing functions $h_o^*$ (resp. $h_e^*$) are defined on all non-root vertices of $G^*$, and are constant along edges in $\ovAos$ (resp. $\ovAes$). We will now relax this condition for some of the edges in $B_o^*\cap \ovAos$ (resp. $B_e^*\cap \ovAes$).
Recall that $A_o^*\subseteq B_o^*$ and $A_e^*\subseteq B_e^*$. 
A function $g^*$ defined on the set of non-root vertices of $G^*$ is called \emph{$B_o^*$-increasing} (resp. \emph{$B_e^*$-increasing}) if it satisfies the following conditions: 
\begin{compactitem}
  \item[(i)] $g^*(u) = g^*(v)$ for every non-root edge $\{u,v\}$ not in $B_o^*$ (resp. not in $B_e^*$),
  \item[(ii)] $g^*(u) \leq g^*(v)$ for every arc $(u,v)$ in $B_o^* \backslash A_o^*$ (resp. $B_e^* \backslash A_e^*$),
  \item[(iii)] $g^*(u) < g^*(v)$ for every arc $(u,v)$ in $A_o^*$ and (resp. $A_e^*$) and every odd (resp. even) diagonal arc $(u,v)$.
\end{compactitem}

The only difference between the definition of $B_o^*$-increasing functions and $\aA_o^*$-increasing functions is that for $\aA_o^*$-increasing functions the requirement $g^*(u) \leq g^*(v)$ in Condition (ii) would be replaced by $g^*(u) = g^*(v)$.
Hence, any $\aA_o^*$-increasing function is $B_o^*$-increasing. Similarly, any $\aA_e^*$-increasing function is $B_e^*$-increasing. 
Hence replacing $\aA_o^*$ and $\aA_e^*$-increasing functions by $B_o^*$ and $B_e^*$-increasing functions in the increasing function algorithm could potentially reduce the grid size, as long as we can show this does not break planarity. We will indeed prove the following result.

\begin{prop}\label{prop:method1-planar}
Let $G^*$ be a dual-adapted rooted 3,4-map and let $\mathcal{L^*}$ be a dual 4-GS labeling. Let $g_o^*$ be a $B_o^*$-increasing function, and let $g_e^*$ be a $B_e^*$-increasing function. Placing each non-root vertex $v$ at coordinates $(g_o^*(v), g_e^*(v)) \in \br^2$ (and drawing  each non-root edges following the rules specified below) yields an SPO drawing of $G^*$.
\end{prop}

Proposition~\ref{prop:method1-planar} is illustrated in Figure~\ref{fig:method1_draw}. In that example, using $B_o^*$ and $B_e^*$-increasing functions (instead of $\aA_o^*$ and $\aA_e^*$-increasing functions) reduces the horizontal dimension of the drawing by 1.

\fig{width=\linewidth}{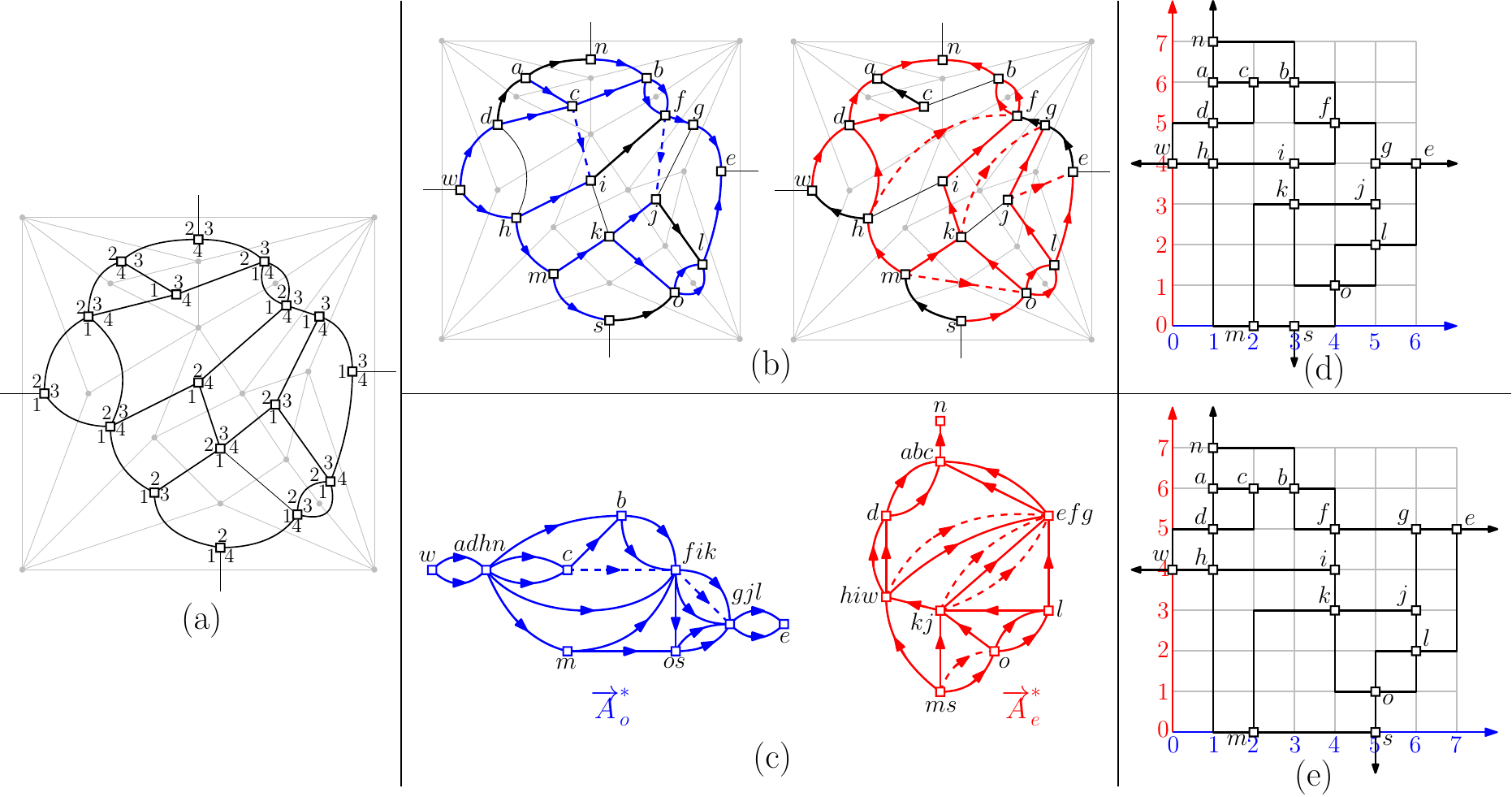}{(a) A rooted 3,4-map with a dual 4-GS labeling. (b) The bipolar orientations $B_o^*$ and $B_e^*$ together with the diagonal arcs. (c) The augmented bipolar contractions $\protect\aA_o^*$ and $\protect\aA_e^*$. (d) The drawing output using the $B_o^*$ and $B_e^*$-increasing functions. (e) The drawing output using the $\protect\aA_o^*$ and $\protect\aA_e^*$-increasing functions.}

\begin{proof}
The proof of Proposition~\ref{prop:method1-planar} is similar to the one of Theorem~\ref{thm:increasing_dual}. We need to check that the drawing is still valid with the relaxed condition on the coordinate functions.

Let $f$ be a non-root face of $G^*$. Recall the definition of the $i$-runs of $f$ from Section~\ref{sec:planarity_ortho}.
Using a similar analysis as in Section~\ref{sec:planarity_ortho}, one can conclude that the 4-run is still drawn weakly right-up; the only difference is that arcs that were previously drawn as a straight segment (going either up or right) may now be drawn as a bent-edge going right-up.
This is illustrated in Figure~\ref{fig:dual_edge3}(a). This figure also gives the rule for drawing the edge in each situation.
By symmetry, the 1-run (resp. 2-run, 3-run) is drawn weakly down-right, (resp. left-down, up-left).\footnote{We mention, that the rule for drawing fully-colored edges (for instance an edge that is on a 4-run of a face, and a 2-run of the opposite face) cannot be made fully symmetric: the bend has to be convex for one of the adjacent face, unlike the other edges on $i$-runs. However, the precise rule for drawing these edges is irrelevant for the planarity argument given below.} 

\fig{width=\linewidth}{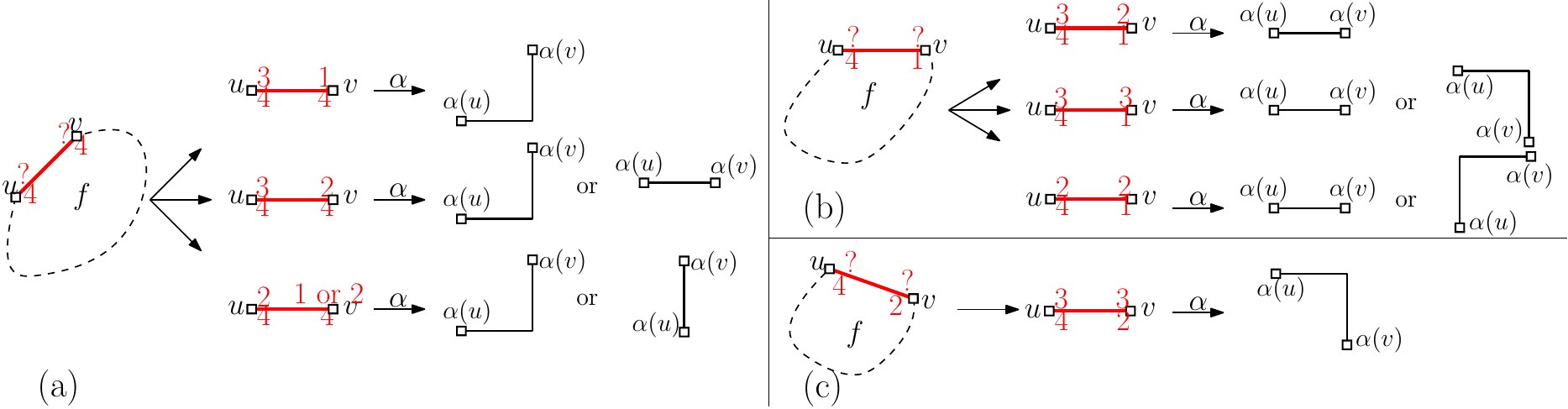}{SPO drawing with $B_o^*,B_e^*$ increasing coordinates. (a) Embedding of edges with labels $\{4,4\}$ inside a face $f$. (b) Embedding of edges with labels $\{4,1\}$ inside $f$. (c) Embedding of edges with labels $\{4,2\}$ inside $f$.}

Moreover an edge connecting the 4-run and the 1-run of $f$ (that is, an edge incident with corner labels $\{4,1\}$ around $f$) is drawn "rightward", that is, either right, or right-up or down-right. This is illustrated in Figure~\ref{fig:dual_edge3}(b). 
This orientation property ensures that the 4-run and 1-run do not touch each other. By symmetry, consecutive runs cannot touch each other. Figure~\ref{fig:dual_edge3} also gives the rule for drawing the edge connecting the 4-run and the 1-run (if it exists), or the edge connecting the 4-run and the 2-run (if it exists). Note that the 
rule for drawing edges ensures that if an edge around $f$ is incident to two distinct labels and has a bend, then this bend is convex with respect to $f$. Therefore these edges cannot intersect with the runs. Lastly, the condition on diagonal arcs prevent opposite runs (e.g. the 4-run and the 2-run) from touching each other. This implies that every inner face has a valid drawing.



We can also investigate the shape of the outer face.
As before one can show that the $(e_1^*,e_2^*)$-run (resp. $(e_2^*,e_3^*)$-run, $(e_3^*,e_4^*)$-run, $(e_4^*,e_1^*)$-run) is directed right-up (resp. down-right, left-down, up-left). This is enough to conclude that the drawing of the outer face is valid by using the same reasoning as before. This implies that the overall drawing is planar. 
\end{proof}


\medskip


\ni \textbf{Relaxing the inequalities on some diagonal arcs:} 
The second optimization we present is about weakening some inequalities for the $B_o^*$ and $B_e^*$-increasing functions $g_o^*,g_e^*$. In the increasing function algorithm of Proposition~\ref{prop:method1-planar}, the increasing functions $g_o^*, g_e^*$ are required to be strictly increasing along the diagonal arcs. We claim that some of these inequalities can be relaxed to weak inequalities.

Let $G^*$ be a rooted $3,4$-map, and let $\cL^*$ be a dual 4-GS labeling. We call \emph{quasi $B_o^*$-increasing} a function $g^*$ on the non-root vertices of $G^*$ which satisfies all the conditions of $B_o^*$-increasing functions except that for odd diagonal edges $(u,v)$ the condition $g^*(u)<g^*(v)$ is replaced by $g^*(u)\leq g^*(v)$. We define \emph{quasi $B_e^*$-increasing} functions similarly.

Note that using quasi $B_o^*$ and $B_e^*$-increasing coordinates (instead of strict $B_o^*$ and $B_e^*$-increasing coordinates) could lead to a non-planar drawing. Indeed, recall that the condition on diagonal arcs was used to ensure the planarity of the drawing of each non-root face (precisely, these conditions ensure that opposite runs around a face do not intersect). But it is clear that in the example represented in Figure~\ref{fig:diag_bad}(a), the weakening of the inequality for the diagonal arc $(u,v)\in\aA_o^*$ could lead to a violation of planarity (the 4-run could intersect the 2-run). 

We now identify which diagonals arcs can be weakened. Consider an SPO drawing of a rooted $3,4$-map $G^*$ obtained by some $B_o^*$ and $B_e^*$-increasing functions $(g_o^*,g_e^*)$. 
Let $f$ be a non-root face of $G^*$ containing a diagonal arc $(u,v)\in \aA_o^*$ (so that the 4-run and the 2-run of $f$ are non-empty). 
Let $C_4$ be the lowest point with the same abscissa as $u$ on the drawing of the 4-run of $f$, and let $C_2$ be the highest point with the same abscissa as $v$ on the drawing of the 2-run of $f$. The points $C_2$ and $C_4$ are indicated in Figure \ref{fig:diag_bad}.
We emphasize that $C_4$ and $C_2$ are geometric features of the embedding of $f$, not features of $f$ itself (in particular, they may correspond to actual vertices on $f$ or just bends of the drawings of edges). 
It is easy to see that as long as $C_4$ is strictly above $C_2$, then the inequality $g_o^*(u) < g_o^*(v)$ induced by the diagonal arc $(u,v)$ can be weakened (that is, replaced by the condition  $g_o^*(u) \leq g_o^*(v)$). This is illustrated in Figure \ref{fig:diag_bad}(b).
In this case, we say that the diagonal arc $(u,v)$ is $\aA_o^*$-\emph{weak}. 

\fig{width=\linewidth}{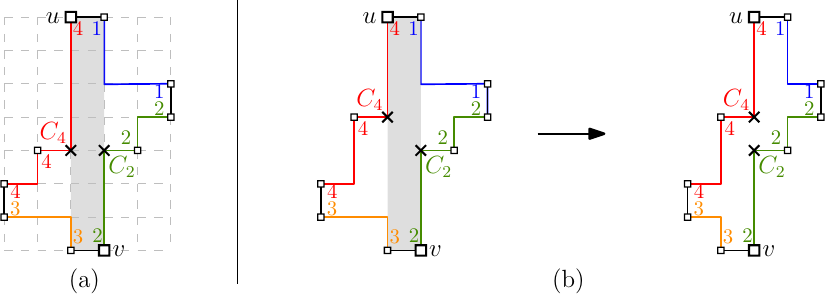}{(a) A face for which the diagonal arc $(u,v) \in \protect\aA_o^*$ cannot be weakened. (b) If $C_4$ is strictly above $C_2$, then $(u,v)$ can be weakened.}

Note that the notion of $\aA_o^*$-weakness for a diagonal arc in $\aA_o^*$ actually depends on the $B_e^*$-increasing function $g_e^*$. Hence this type of optimization needs to be performed one coordinate at a time. 
Given an SPO drawing of $G^*$ where the coordinate functions $(g_o^*,g_e^*)$ are quasi $B_o^*$-increasing and $B_e^*$-increasing respectively, we can define $\aA_e^*$-\emph{weak} diagonals similarly as above: for a diagonal arc $(u,v)\in\aA_e^*$ in a face $f$ of $G^*$ we consider the rightmost point $C_3$ with the same ordinate as $u$ on the drawing of the 3-run of $f$ and the leftmost point $C_1$ with same ordinate as $v$ on the drawing of the 1-run of $f$, and say that the diagonal $(u,v)$ is $\aA_e^*$-\emph{weak} if $C_3$ is strictly to the left of $C_1$. In this case the inequality $g_e^*(u)<g_e^*(v)$ can be weakened without breaking planarity.

This suggests the following optimized SPO drawing algorithm:
\begin{compactenum}
  \item Compute a (strict) $B_e^*$-increasing function $\wt{g_e^*}$, and find all the $\aA_o^*$-weak diagonal arcs.
  \item Compute a quasi $B_o^*$-increasing function $g_o^*$, which is weakly increasing along $\aA_o^*$-weak diagonal arcs (and strictly increasing along the other diagonal arcs).
  \item Compute a quasi $B_e^*$-increasing function $g_e^*$, which is weakly increasing along $\aA_e^*$-weak diagonal arcs (and strictly increasing along the other diagonal arcs).
  \item Draw the graph $G^*$ using the coordinate function $(g_o^*,g_e^*)$.
\end{compactenum}
The above discussion shows that the result of this algorithm is planar (since the validity of the drawing of the inner faces is preserved, and the validity of the drawing of the outer face does not depend on diagonal arcs).

A drawing obtained using this optimized algorithm is represented in Figure~\ref{fig:method2_draw}. In this example, the weakening of diagonal edges inequality leads to a reduction of the vertical dimension by 1.

\begin{remark}
This optimization strategy is analogous to the reduction operations used in the 4-valent case~\cite{OB-EF:Schnyder}. 		
\end{remark} 

\fig{width=\linewidth}{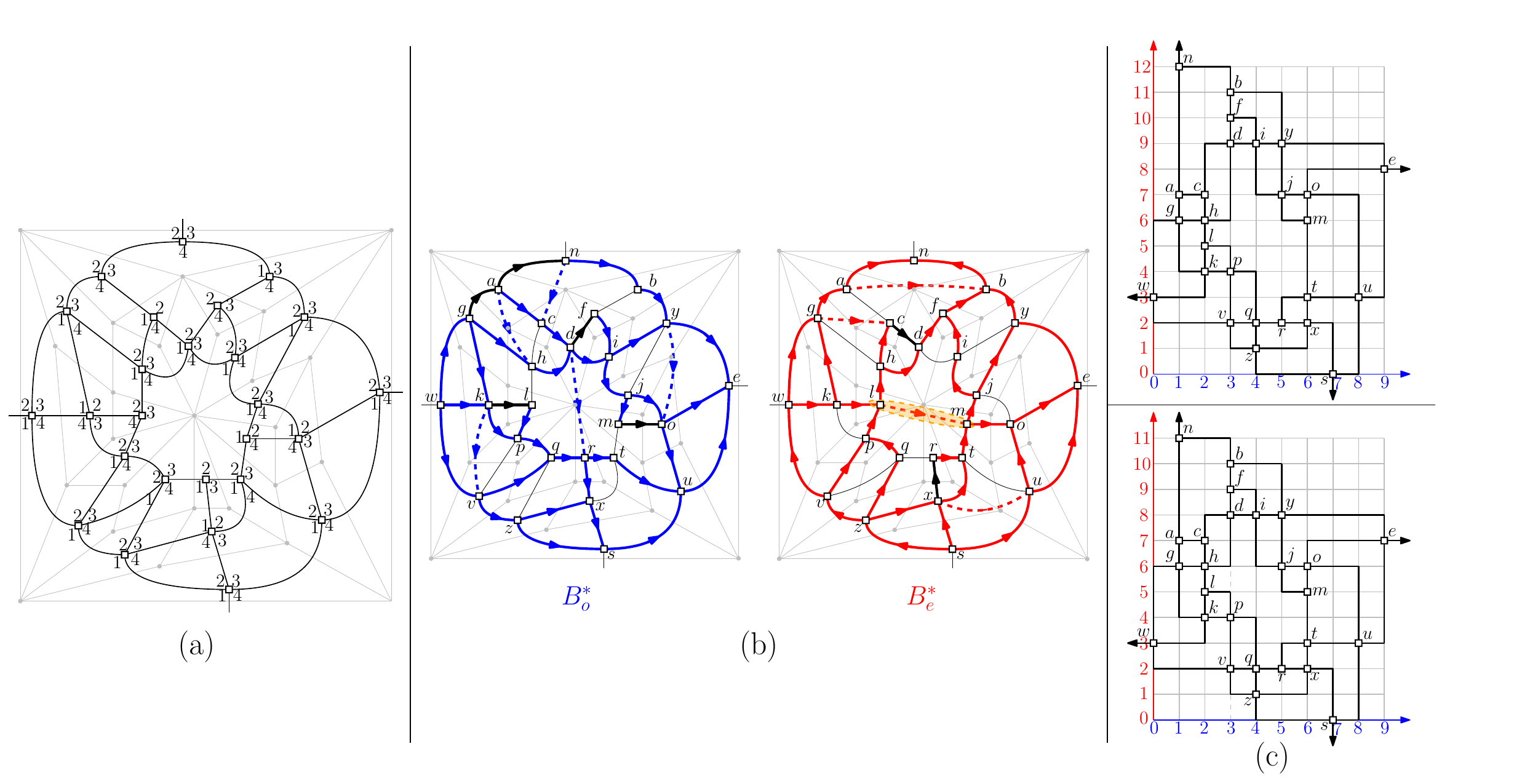}{(a) A dual-adapted rooted 3,4-map with a 4-GS labeling. (b) The bipolar orientations $B_o^*$ and $B_e^*$ plus the diagonal arcs. Those edges in $B_o^* \cap \ov{A}_o^* (B_e^* \cap \ov{A}_e^*)$ are colored as black. The odd diagonal arc $(l,m)$ is weak. (c) Top: The drawing output using strict $B_o^*$ and $B_e^*$-increasing functions. Bottom: The drawing output using the weak $B_o^*$ and $B_e^*$-increasing functions obtained using the algorithm outlined above.}

\section{Additional remarks and optimization}\label{sec:conclusion}
 

\subsection{Variants on the face counting algorithms} \label{sec:variant-face-count}
The face-counting algorithms given in Sections~\ref{sec:face-counting-primal} and~\ref{sec:face-counting-dual} admit some variants that we now describe.

Consider the primal picture first. 
Let $G$ be an adapted 3,4-angulation of the square and let $\cL$ be a 4-GS labeling of~$G$. In the face-counting algorithm of Section~\ref{sec:face-counting-primal}, each vertex $v$ is placed at the point $(\ell_e(v),r_o(v))$, where 
$r_o(v)$ is the number of inner faces of $B_o$ at the right of $P_o(v)$ and $\ell_e(v)$ is the number of inner faces of $B_e$ at the left of $P_e(v)$.  Given some positive weights associated to each inner face of $B_o$ the $B_e$, one can consider the vertex placement $(\ell'_e(v),r'_o(v))$ where 
$r_o'(v)$ is the sum of the weights of the inner faces of $B_o$ at the right of $P_o(v)$ and $\ell_e'(v)$ is the sum of the weights of the inner faces of $B_e$ at the left of $P_e(v)$. We claim that this vertex placement lead to a planar drawing of $G$ for any choice of positive weight. This generalizes Theorem~\ref{thm:face-counting-primal} which corresponds to a weight of 1 for every face.
Indeed, one can prove that $\ell'_e(v)$ is an $\aA_o$-increasing function and $r'_o(v)$ is an $\aA_e$-increasing function exactly as before; hence planarity is again a consequence of Theorem~\ref{thm:increasing-primal}. Similarly one could define some other placement $(\ell''_e(v),r''_o(v))$ using vertex-counting instead, as they would still be a special case of the increasing function algorithm. 

In the dual picture one can also consider variants of the face-counting algorithm by assigning weights to the faces of $B_o^*$ and $B_e^*$; or one could define vertex-counting variants as all of these variants fall under the ``increasing function'' umbrella.

\begin{remark} 
The face-counting algorithm of Section~\ref{sec:face-counting-primal} is somewhat reminiscent of Schnyder's famous drawing algorithm for triangulations~\cite{Schnyder:wood2}. The relation is even closer for a variant of the algorithm that we now present. 
Let $G$ be an adapted 3,4-angulation of the square with $\ff$ inner faces, and let $\cL$ be a 4-GS labeling of~$G$. 
Consider the vertex placement given by $\al'(v)=(\ell_e'(v),r_o'(v))$, where $r_o'(v)$ is the number of inner faces of $G$ at the right of $P_o(v)$ and $\ell_e'(v)$ is the number of inner faces of $G$ at the left of $P_e(v)$. The analogue of Theorem~\ref{thm:face-counting-primal} holds for this vertex placement. Indeed it corresponds to one of the face-weighted variants described above (where the weight of a face of $B_o$ or $B_e$ is the number of faces of $G$ that it contains). Now, for a vertex $v$ of $G$, we can consider the paths $P_1(v),\ldots,P_4(v)$ in $W_1,\ldots,W_4$ between $v$ and $v_4,v_1,v_2,v_3$ respectively. These paths are non-crossing as proved in~\cite{OB-EF-SL:Grand-Schnyder}, and denoting by $r_i(v)$ the number of inner faces between $P_{i-1}(v)$ and $P_{i+1}(v)$ one gets $r_o'(v)=r_2(v)+r_3(v)$ and $\ell_e'(v)=r_3(v)+r_4(v)$. Hence, in barycentric coordinates, one has
$$\al'(v)=\lambda_1(v)\, v_1+\lambda_2(v)\, v_2+\lambda_3(v)\, v_3+\lambda_4(v)\, v_4,$$
where $v_1=(0,0)$, $v_2=(0,\ff)$, $v_3=(\ff,\ff)$, $v_4=(\ff,0)$, and $\lambda_i(v)=\frac{r_i(v)}{\ff}$ for all $i$. This barycentric vertex placement given by the number of inner faces between consecutive Schnyder wood paths is a direct analogue of Schnyder's algorithm~\cite{Schnyder:wood2}.
\end{remark}

\subsection{Triangulating under the 4-connectivity constraint}\label{sec:triangulating}
Given $G$, an adapted 3,4-angulation of the square, it is always possible to add a diagonal in each quadrangular inner face in such a way that the resulting triangulation of the square is adapted. Indeed, it is easy to see that for every quadrangular inner face of $G$, at least one of the two diagonals can be added without creating a separating 3-cycle. As we now explain, if $G$ is endowed with a 4-GS labeling then this decision can be made very easily, and moreover the 4-GS labeling can be maintained all along, in time $O(1)$ for each edge addition, to end with an irreducible triangulation of the square endowed with a 4-GS labeling. 

Let $f$ be a quadrangular inner face of $G$. Let $c$ be a corner incident to $f$, let $e$ be the edge visited after $c$ in clockwise order around $f$, and let $c'$ be the corner incident to $e$ opposite to $c$ (i.e., at the other extremity and on the other side). The corner $c$ is called \emph{special} if the label-jump from $c$ to $c'$ is $1$. 
A corner is called \emph{diagonal} if it is special and preceded by a non-special corner in clockwise order around $f$. 
With $c_1,c_2,c_3,c_4$ the four corners around $f$ (the index indicating the label), we choose a pair $c_i,c_{i+2}$ of opposite corners such that the other corners $c_{i+1}$ and $c_{i-1}$ are not diagonal; such a pair exists since $f$ can not have two consecutive diagonal corners. 

\begin{figure}
\begin{center}
\includegraphics[width=\linewidth]{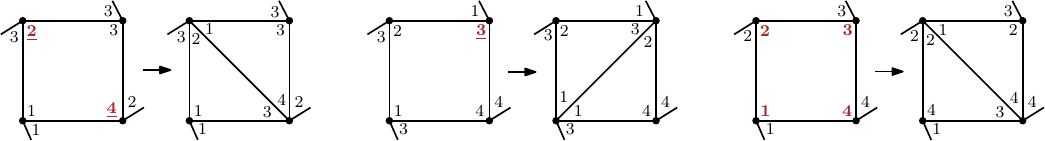}
\end{center}
\caption{Some cases for the configuration of special/non-special labels in a quadrangular inner face, where special labels are bold colored, with the diagonal ones underlined.}
\label{fig:edge_additions}
\end{figure}

We then perform the following operations, see Figure~\ref{fig:edge_additions} for examples:
\begin{enumerate}
\item
decrease by $1$ the labels of the special corners around $f$, 
\item
 insert a diagonal $e$ in $f$ connecting the vertex $v$ at $c_i$ and the vertex $v'$ at $c_{i+2}$, which creates a new corner at $v$ (resp. at $v'$), the new corner being considered as the one on the right of $e$, looking from $v$ (resp. from $v'$), 
\item
give label $i$ (resp. $i+2$) to the new corner at $v$ (resp. $v'$).  
\end{enumerate}
 
The conditions of a 4-GS labeling are still satisfied, the crucial point being that the labels of special corners have been decreased in order for Condition (L3) to remain satisfied after adding the diagonal~\footnote{All edge-additions can actually be performed simultaneously, since the status special/non-special of a corner in a quadrangular face is not affected by label changes in other faces.}. 
Moreover, with the notation of Section~\ref{sec:grid-size}, any edge contributing to $\dd_e$ will keep contributing to $\dd_e$ all along subsequent edge-additions, and similarly for $\dd_o$. Hence, the bounds on the grid width and height in Proposition~\ref{prop:grid-bound} are weakly decreasing when adding diagonals by this process (however the actual optimized grid-size may increase).

Let us mention two related previous works. Biedl and Kant have given in~\cite{biedl1997triangulating} a linear-time algorithm that, starting from a \plm with no separating triangle, adds edges incrementally while avoiding the creation of separating triangles, yielding a 4-connected triangulation of the original \plmm. More closely related, in~\cite{fusy2007straight} it is shown that, starting from a simple quadrangulation endowed with a separating decomposition, there is a simple process to add diagonals within faces to yield an irreducible triangulation of the square endowed with a transversal structure. What we describe here can be seen as an extension of this process to adapted 3,4-angulations of the square.



\subsection{Bend minimization}\label{sec:bend}
In this subsection we explain how to minimize the number of bends in the orthogonal drawings obtained by our increasing function algorithm for rooted 3,4-maps. More precisely, given a dual-adapted rooted 3,4-map $G^*$, we want to determine a dual 4-GS structure which leads to a drawing with the least number of bends. 
For this, we will adapt a method pioneered by Tamassia~\cite{tamassia1987embedding}.

\begin{definition}
An SPO drawing of a rooted 3,4-map $G^*$ is \emph{one-bend} if every edge has at most one bend. We use the abbreviation \emph{BSPO} for such a drawing of $G^*$. Two BSPO drawings of $G^*$ are called \emph{equivalent} if they represent the same rooted \plmm, and each half-edge has the same direction in both drawings. The \emph{combinatorial type} of a BSPO drawing is the equivalence class for this relation.
\end{definition}

\begin{figure}[h!]
\begin{center}
\includegraphics[]{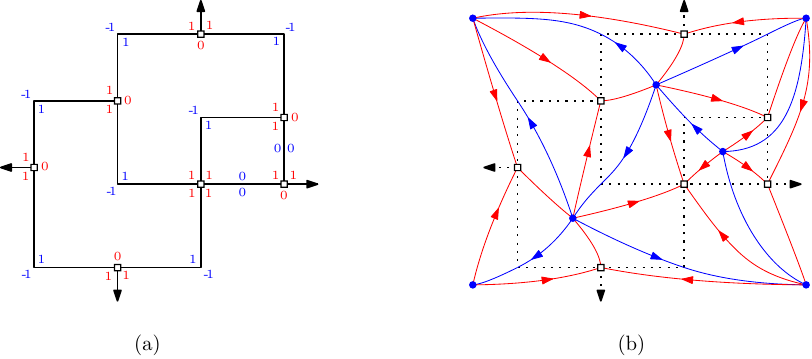}
\end{center}
\caption{(a): a BSPO drawing of a rooted 3,4-map $G^*$, with the encoding of the combinatorial type by a Tamassia assignment. 
(b): the encoding by a valid Tamassia flow: the arrows indicate the arcs of $\Gbox$ for which the flow is~1.}
\label{fig:ortho}
\end{figure}

As shown in Figure~\ref{fig:ortho}(a), and following Tamassia~\cite{tamassia1987embedding}, the combinatorial type of a BSPO drawing can be encoded by assigning values in $\{-1,0,1\}$ at vertex-corners and at edge-sides, with value $1$ (resp. $0$, $-1$) indicating an angle of $\pi/2$ (resp. $\pi$, $3\pi/2$). 
It was shown in~\cite{tamassia1987embedding} that such an assignment encodes a valid combinatorial type (non-empty equivalence class of drawings) if and only if the following conditions are satisfied:
\begin{itemize}
\item
the sum of values at corners around each vertex $v$ is $2\,\mathrm{deg}(v)-4$, 
\item
the sum of edge-side values at each edge is $0$, 
\item
the sum of values in each inner face is $4$,
\item
with the outer contour partitioned by the dangling half-edges, the sum of values in each of the 4 portions is $1$.
\end{itemize}
Note that these conditions are clearly necessary in order to be able to use the assignment to consistently propagate the directions among $\{$west, north, east, south$\}$ to all half-edges after fixing the direction of the 4 dangling half-edges.



We call \emph{Tamassia assignment} such an assignment of values to corners and edge-sides of~$G^*$. A Tamassia assignment for $G^*$ can equivalently be formulated as a flow on a related graph~$\Gbox$, as shown in Figure~\ref{fig:ortho}(b). Let $\Gbox$ be the graph obtained from $G^*$ by
\begin{compactitem}
\item placing a \emph{blue vertex} $v_f$ in each face $f$ of $G^*$ and tracing a \emph{blue edge} (joining 2 blue vertices) across each edge of $G^*$, 
\item adding a \emph{red edge} between blue vertex $v_f$ and each of the corner of the face $f$,
\item deleting the edges of $G^*$.
\end{compactitem}
Given an assignment of a \emph{flow} in $\{0,1\}$ to each arc of $\Gbox$, we call \emph{net flow} of a vertex $v$ the sum of the flows of the outgoing arcs at $v$, minus the sum of the flows of the incoming arcs at $v$. We call \emph{Tamassia flow} an assignment of a \emph{flow} in $\{0,1\}$ to each arc of $\Gbox$ such that 
\begin{itemize}
\item
every vertex $v$ of $G^*$ has net flow $4-2\,\mathrm{deg}(v)$,
\item
every blue vertex $v_f$ has net flow $4$ if $f$ is a non-root face, and net flow 1 otherwise.
\end{itemize} 
It is easy to see that the Tamassia assignments are in bijection with the Tamassia flows such that if an arc of $\Gbox$ has flow 1, then the opposite arc has flow 0; such Tamassia flows are called \emph{valid}. 
We call \emph{min-cost} a Tamassia flow minimizing the sum of the flow values of arcs on blue edges. Note that the min-cost Tamassia flows are valid and they correspond bijectively to the (combinatorial type of) BSPO drawings of $G^*$ minimizing the number of bends. 
Note that $\Gbox$ admits a Tamassia flow (which can be obtained in linear time by applying our drawing algorithm on $G^*$), hence a min-cost one can be found in time $O(n^{3/2}\log(n))$, where $n$ is the number of vertices of $G^*$ (see~\cite{goldberg2015minimum} and references therein).

In order to combine the above framework with our algorithm, we need to establish a correspondence between the combinatorial type of the BSPO of $G^*$ and the dual 4-GS labelings of~$G^*$.

\begin{figure}
\begin{center}
\includegraphics[width=\linewidth]{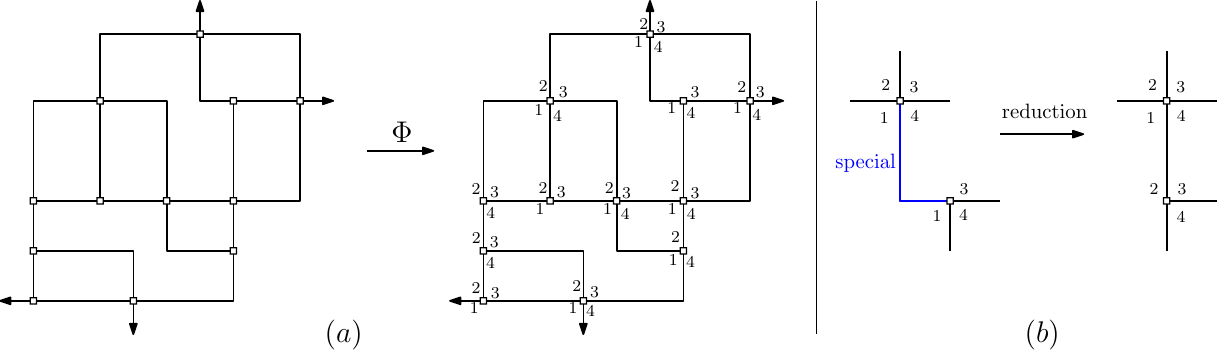}
\end{center}
\caption{(a) Illustrating $\Phi$: A right-chiral BSPO drawing and the corresponding dual 4-GS labeling. (b) Proof that if a dual GS-labeling has a special edge, then it does not correspond to a BSPO type with minimal number of bends.}
\label{fig:Correspond}
\end{figure}

As illustrated in Figure~\ref{fig:Correspond}(a), there is a simple rule to associate to each BSPO drawing of $G^*$ a \emph{4-labeling} of $G^*$, that is, an assignment of a label in $\four$ to each corner of $G^*$: each corner is labeled according to the direction of the half-edge delimiting it on the right side, with labels $1,2,3,4$ corresponding to the directions west, north, east and south respectively. This defines a map $\Phi$ from the set of combinatorial types of the BSPO of $G^*$ to the set of $4$-labelings of $G^*$.

We now consider some particular types of BSPO drawings. 

\begin{definition}
A BSPO drawing of a rooted 3,4-map is called \emph{right-chiral} if it satisfies 
the following properties:
\begin{itemize}
\item[(A1)] no edge leaving a vertex of degree 3 turns left,
\item[(A2)] each dangling half-edge has a corner of angle $\pi/2$ on its right. 
\end{itemize}
\end{definition}

Note that Condition (A2) is close to redundant: it is a consequence of (A1) when there is at least one vertex of degree 4 on the outer contour. 
Moreover, Condition (A1) ensures that the edges connecting two vertices of degree 3 have no bend. 
Note also that when all vertices have degree 4, any BSPO drawing is right-chiral.

\begin{lemma}\label{lem:corresp}
Let $G^*$ be a rooted 3,4-map. 
The mapping $\Phi$ is a bijection between the combinatorial types of right-chiral BSPO drawings of $G^*$ and the dual 4-GS labelings of $G^*$.
\end{lemma}
\begin{proof}
Let $T$ be a combinatorial type of right-chiral BSPO drawing, and let $\cL^*=\Phi(T)$. It is easy to check that the sum of label-jumps in clockwise order around each vertex is $4$, and the sum of label-jumps in counterclockwise order around each edge is $4$. 
Moreover, Property (A1) directly translates to the condition (L3$^*$) at edges. 
It remains to show that there is no label-jump in each of the four outer faces, and that the sum of clockwise label-jumps in each inner face is $4$. 
Let $F_i$ be the set of inner faces, and $F_o$ the set of 4 outer faces. 
For $f\in F_i\cup F_o$, let $J(f)$ be the sum of clockwise label-jumps around $f$, 
thus a multiple of $4$. It is actually a positive multiple of $4$ for $f\in F_i$. 

Moreover, (A2) ensures that the clockwise label-jumps across the dangling half-edges are $+1$. These are called \emph{root label-jumps}. Letting $J_V$ be the sum of 
non-root clockwise label-jumps around vertices, we thus have $J_V=4|V|-4$. Letting $J_E$ be the sum of counterclockwise label-jumps around edges, we have $J_E=4|E|$. 
Each label-jump in counterclockwise around an edge either comes from a non-root label-jump around a vertex, or from a clockwise label-jump around a face. Hence
\[
\sum_{f\in F_i\cup F_o}J(f) = J_E-J_V=4|E|-4|V|+4=4|F_i|,
\]
 by the Euler relation. Since $J(f)\geq 4$ for $f\in F_i$, we conclude that $J(f)=4$ for $f\in F_i$ and $J(f)=0$ for $f\in F_o$. Hence, $\cL^*$ is a dual 4-GS labeling. 

Conversely, for a dual 4-GS labeling $\cL^*$, applying the local rule of $\Phi$ allows to determine the
directions of all half-edges, yielding a valid combinatorial type $T$ of BSPO drawing. Condition~(L3$^*$) then guarantees that (A1) is satisfied, and Condition (A2) is also 
clearly satisfied, so that $T$ is right-chiral. And by construction the two mappings are inverse of each other. 
\end{proof}

As a consequence of Lemma~\ref{lem:corresp}, a rooted 3,4-map $G^*$ admits a right-chiral BSPO drawing if and only if it admits a dual 4-GS labeling, and this holds if and only if $G^*$ is dual-adapted\footnote{It is shown in~\cite{FelsnerKV14} that a graph with max-degree $4$ and a marked vertex $v_\infty$ of degree $4$ admits a (possibly non-planar) one-bend  suspended orthogonal drawing if and only if the following condition (P) is satisfied: for every vertex-subset $S$ not containing  $v_\infty$, there are at most $2|S|-2$ edges with both ends in $S$.  	
		For consistency one can check that, for a   rooted 3,4 map $G^*$, being dual-adapted implies (P). Moreover, if $G^*$ is 4-regular, then being dual-adapted is equivalent to (P).}
		. In this case, one can easily adapt the approach described above in order to find a right-chiral BSPO drawing of $G^*$ having a minimal number of bends. Let us call \emph{right-chiral} a Tamassia flow of $\Gbox$ such that the flow is 0 on every blue arc of $\Gbox$ having a vertex of $G^*$ of degree 3 on its right. Then the combinatorial types of right-chiral BSPO drawings of $G^*$ are in bijection with the right-chiral valid flows of $\Gbox$. In this bijection, the combinatorial types of right-chiral BSPO drawings of $G^*$ having a minimal number of bends 
correspond to the min-cost right-chiral Tamassia flows (as before, min-cost minimizes the sum of the flow values of arcs on blue edges). 
Moreover such a min-cost right-chiral Tamassia flow can be found in time $O(n^{3/2}\log(n))$, where $n$ is the number of vertices of $G^*$.\\

Finally, recall from Section~\ref{sec:prop-edges-dual} that given a dual 4-GS labeling, a non-root arc $a$ is called \emph{special} if 
$$(\text{right-init}(a),\text{right-term}(a),\text{left-term}(a), \text{left-init}(a)) = (k,k,k+2,k+3)$$ for some $k\in\four$. We claim that the 
dual 4-GS labeling corresponding to a min-cost right-chiral Tamassia flow has no special arc. Indeed, if a dual 4-GS labeling $\cL^*$ has a special arc $a$, then one can change one of the labels incident to $a$ and get a new dual 4-GS labeling $\cL^*$ which corresponds (via $\Phi$) to a combinatorial type of right-chiral BSPO with one less bend. This is illustrated in Figure~\ref{fig:Correspond}(b).

Since the dual 4-GS labeling $\cL^*$ corresponding to a min-cost right-chiral Tamassia flow has no special arc, Proposition~\ref{prop:orientation} ensures that the BSPO of $G^*$ obtained using the increasing function drawing algorithm of Section~\ref{sec:increasing_dual} has combinatorial type equal to $\Phi^{-1}(\cL^*)$. Thus the obtained BSPO has minimal number of bends among the left-chiral BSPOs of $G^*$. This is summarized below.

\begin{prop}\label{prop:optimal}
Let $G^*$ be a dual-adapted rooted 3,4-map with $n$ vertices, and 
let $\mathbf{T}(G^*)$ be the set of combinatorial types of right-chiral BSPO drawings of $G^*$. There is an algorithm with time complexity $O(n^{3/2}\log(n))$ 
for finding an element $T\in\mathbf{T}(G^*)$ with minimal number of bends. Moreover, letting $\cL^*=\Phi(T)$ be the associated dual 4-GS labeling, the drawing algorithm of Section~\ref{sec:increasing_dual} applied to $\cL^*$ yields a grid-drawing realization of $T$. 
\end{prop}

\begin{remark}
Proposition~\ref{prop:orientation} ensures that, for any dual 4-GS labeling $\cL^*$, the drawing algorithm of Section~\ref{sec:increasing_dual} applied to $\cL^*$ always yields a right-chiral BSPO drawing of $G^*$ (although its combinatorial type of that BSPO may differ from $\Phi^{-1}(\cL^*)$ if $\cL^*$ has some special edges).  
Hence the optimal 4-GS labeling of $G^*$ in Proposition~\ref{prop:optimal} is also the one for which the algorithm of Section~\ref{sec:increasing_dual} yields the smallest number of bends (among all possible 4-GS labelings of $G^*$).
\end{remark}


\bigskip
\noindent{\bf Conflict of interest statement:} The authors have no conflict of interest to disclose in relation to this article.\\

\smallskip
\noindent{\bf Acknowledgments:} The authors thank the anonymous referees for their many useful suggestions. Olivier Bernardi was partially supported by NSF Grant DMS-2154242. \'Eric Fusy was partially supported by the project ANR19-CE48-011-01 (COMBIN\'E), and the project ANR-20-CE48-0018 (3DMaps).

\bibliographystyle{alpha}
\bibliography{biblio-Schnyder}

\begin{thebibliography}{DFOdMR95}

\bibitem[BETT99]{BattistaETT99}
G.~Di Battista, P.~Eades, R.~Tamassia, and I.G. Tollis.
\newblock {\em Graph Drawing: Algorithms for the Visualization of Graphs}.
\newblock Prentice-Hall, 1999.

\bibitem[BF12]{OB-EF:Schnyder}
O.~Bernardi and {\'E}.~Fusy.
\newblock Schnyder decompositions for regular plane graphs and application to
  drawing.
\newblock {\em Algorithmica}, 62(3):1159--1197, 2012.

\bibitem[BFL24]{OB-EF-SL:Grand-Schnyder}
O.~Bernardi, \'E. Fusy, and S.~Liang.
\newblock Grand {S}chnyder woods.
\newblock {\em Annals of Combinatorics}, 2024.
\newblock arXiv.2303.15630, \url{https://doi.org/10.1007/s00026-024-00729-8}.

\bibitem[BH12]{Barriere-Huemer:4-Labelings-quadrangulation}
L.~Barriere and C.~Huemer.
\newblock 4-labelings and grid embeddings of plane quadrangulations.
\newblock {\em Discrete Mathematics}, 312(10):1722--1731, 2012.

\bibitem[Bie96]{biedl1996optimal}
T.~Biedl.
\newblock Optimal orthogonal drawings of triconnected plane graphs.
\newblock In {\em 5th Scandinavian Workshop on Algorithm Theory, {SWAT}},
  volume 1097 of {\em Lecture Notes in Computer Science}, pages 333--344.
  Springer, Springer, 1996.

\bibitem[BK98]{biedl1998better}
T.~Biedl and G.~Kant.
\newblock A better heuristic for orthogonal graph drawings.
\newblock {\em Computational Geometry}, 9(3):159--180, 1998.

\bibitem[BKK97]{biedl1997triangulating}
T.~Biedl, G.~Kant, and M.~Kaufmann.
\newblock On triangulating planar graphs under the four-connectivity
  constraint.
\newblock {\em Algorithmica}, 19:427--446, 1997.

\bibitem[BLK18]{BhatiaLK18}
S.~Bhatia, K.~Lad, and R.~Kumar.
\newblock Bend-optimal orthogonal drawings of triconnected plane graphs.
\newblock {\em {AKCE} Int. J. Graphs Comb.}, 15(2):168--173, 2018.

\bibitem[DFOdMR95]{de1995bipolar}
H.~De~Fraysseix, P.~Ossona~de Mendez, and P.~Rosenstiehl.
\newblock Bipolar orientations revisited.
\newblock {\em Discrete Applied Mathematics}, 56(2-3):157--179, 1995.

\bibitem[dV03]{CDV:these}
\'E.~Colin de~Verdi\`ere.
\newblock {\em Shortening of curves and decomposition of surfaces}.
\newblock PhD thesis, Universit{\'e} Paris 7, 2003.

\bibitem[Fel01]{Felsner:woods}
S.~Felsner.
\newblock Convex drawings of planar graphs and the order dimension of
  3-polytopes.
\newblock {\em Order}, 18:19--37, 2001.

\bibitem[FHKO11]{FeHuKa}
S.~Felsner, C.~Huemer, S.~Kappes, and D.~Orden.
\newblock {Binary Labelings for Plane Quadrangulations and their Relatives}.
\newblock {\em {Discrete Mathematics \& Theoretical Computer Science}}, {Vol.
  12 no. 3}, January 2011.

\bibitem[FKV14]{FelsnerKV14}
S.~Felsner, M.~Kaufmann, and P.~Valtr.
\newblock Bend-optimal orthogonal graph drawing in the general position model.
\newblock {\em Comput. Geom.}, 47(3):460--468, 2014.

\bibitem[Fus06]{fusy2007straight}
{\'{E}}.~Fusy.
\newblock Straight-line drawing of quadrangulations.
\newblock In {\em Graph Drawing, 14th International Symposium ({GD})}, volume
  4372 of {\em Lecture Notes in Computer Science}, pages 234--239. Springer,
  2006.

\bibitem[Fus09]{Fu07b}
\'E. Fusy.
\newblock Transversal structures on triangulations: {A} combinatorial study and
  straight-line drawings.
\newblock {\em Discrete Math.}, 309:1870--1894, 2009.

\bibitem[GKHT15]{goldberg2015minimum}
A.V. Goldberg, H.~Kaplan, S.~Hed, and R.E. Tarjan.
\newblock Minimum cost flows in graphs with unit capacities.
\newblock In {\em 32nd International Symposium on Theoretical Aspects of
  Computer Science ({STACS})}, volume~30 of {\em LIPIcs}, pages 406--419.
  Schloss Dagstuhl - Leibniz-Zentrum f{\"{u}}r Informatik, 2015.

\bibitem[He93]{He93:reg-edge-labeling}
X.~He.
\newblock On finding the rectangular duals of planar triangulated graphs.
\newblock {\em SIAM Journal on Computing}, 22:1218--1226, 1993.

\bibitem[NR04]{NishizekiR04}
T.~Nishizeki and Md.~S. Rahman.
\newblock {\em Planar Graph Drawing}, volume~12 of {\em Lecture Notes Series on
  Computing}.
\newblock World Scientific, 2004.

\bibitem[RNN98]{RahmanNN98}
Md.~S. Rahman, S.-I. Nakano, and T.~Nishizeki.
\newblock Rectangular grid drawings of plane graphs.
\newblock {\em Comput. Geom.}, 10(3):203--220, 1998.

\bibitem[RNN99]{RahmanNN99}
Md.~S. Rahman, S.{-}I. Nakano, and T.~Nishizeki.
\newblock A linear algorithm for bend-optimal orthogonal drawings of
  triconnected cubic plane graphs.
\newblock {\em J. Graph Algorithms Appl.}, 3(4):31--62, 1999.

\bibitem[RT86]{RosenstiehlT86}
P.~Rosenstiehl and R.E. Tarjan.
\newblock Rectilinear planar layouts and bipolar orientations of planar graphs.
\newblock {\em Discret. Comput. Geom.}, 1:343--353, 1986.

\bibitem[Sch89]{Schnyder:wood1}
W.~Schnyder.
\newblock Planar graphs and poset dimension.
\newblock {\em Order}, 5(4):323--343, 1989.

\bibitem[Sch90]{Schnyder:wood2}
W.~Schnyder.
\newblock Embedding planar graphs on the grid.
\newblock In {\em Proceedings of the First Annual {ACM-SIAM} Symposium on
  Discrete Algorithms ({SODA})}, pages 138--148. {SIAM}, 1990.

\bibitem[Tam87]{tamassia1987embedding}
R.~Tamassia.
\newblock On embedding a graph in the grid with the minimum number of bends.
\newblock {\em SIAM Journal on Computing}, 16(3):421--444, 1987.

\bibitem[TT86]{TamassiaT86}
R.~Tamassia and I.G. Tollis.
\newblock A unified approach a visibility representation of planar graphs.
\newblock {\em Discret. Comput. Geom.}, 1:321--341, 1986.

\bibitem[TT89]{TamassiaTollis89}
R.~Tamassia and I.G. Tollis.
\newblock Tessellation representations of planar graphs.
\newblock In {\em Proc. 27th Allerton Conf. Commun. Control Comput.}, page
  48857, 1989.

\end{thebibliography}

\end{document}